\documentclass[11pt, reqno]{amsart}

\usepackage{etoolbox}

\usepackage{amscd,amsmath}
\usepackage{mathrsfs}
\usepackage{amsfonts}
\usepackage{centernot}

\usepackage{amssymb, bm, xspace}
\usepackage{enumerate}

\usepackage{cite}

\usepackage[OT2,OT1]{fontenc}

\usepackage{datetime}

\usepackage{enumitem}

\usepackage[dvipsnames]{xcolor}

\usepackage[pagebackref=true, colorlinks=true, citecolor=blue]{hyperref}

\usepackage[margin=1.2in]{geometry}

\usepackage{mathtools}

\mathtoolsset{centercolon}

\usepackage{cleveref}

\usepackage{stackengine}
\stackMath

\newbool{HaveBBM}
\booltrue{HaveBBM}

\ifbool{HaveBBM}{
	\usepackage{bbm}
    }
    {
    }

\renewcommand{\mathbf}[1]{\bm{#1} \textbf{ *** Use bm instead of mathbf ***}}

\newcommand{\eqn}{\begin{eqnarray}}
\newcommand{\een}{\end{eqnarray}}

\newtheorem{theorem}{Theorem}[section]
\newtheorem*{theorem*}{Theorem}				
\newtheorem{prop}[theorem]{Proposition}
\newtheorem{lemma}[theorem]{Lemma}
\newtheorem{cor}[theorem]{Corollary}

\newtheorem{definition}[theorem]{Definition}
\newtheorem{remark}[theorem]{Remark}
\newtheorem*{remark*}{Remark}

\numberwithin{equation}{section}

\newcommand{\BoldTau}
    {{\mbox{\boldmath $\tau$}}}
    
\newcommand{\bomega}
    {{\mbox{\boldmath $\omega$}}}
\newcommand{\bOmega}
    {{\mbox{\boldmath $\Omega$}}}
\newcommand{\bmu}
    {{\mbox{\boldmath $\mu$}}}

\newcommand{\BB}[1]{\ensuremath{\mathbb{#1}}}
\newcommand{\R}{{\ensuremath{\BB{R}}}}

\newcommand{\iny}{\ensuremath{\infty}}
\newcommand{\grad}{\ensuremath{\nabla}}

\newcommand{\CharFunc}{
    \ifbool{HaveBBM}{
        \ensuremath{\mathbbm{1}}
        }
        {
        \ensuremath{\bm{1}}
        }
    }

\DeclareMathOperator{\dv}{div} %
\DeclareMathOperator{\curl}{curl} %

\DeclareMathOperator{\dist}{dist} %
\DeclareMathOperator{\diam}{diam} %
\DeclareMathOperator{\supp}{supp} %
\DeclareMathOperator{\cond}{cond} %
\DeclareMathOperator{\lincond}{lincond} %
\DeclareMathOperator{\Dom}{Dom}

\newcommand{\prt}{\ensuremath{\partial}}
\newcommand{\brac}[1]{\ensuremath{\left[ #1 \right]}}

\newcommand{\pr}[1]{\ensuremath{\left( #1 \right)}}

\DeclarePairedDelimiter{\set}{\{}{\}}

\DeclarePairedDelimiterX{\norm}[1]{\lVert}{\rVert}{#1}

\DeclarePairedDelimiterX{\abs}[1]{\lvert}{\rvert}{#1}

\newcommand{\cN}{{c_{\Cal{N}}}}

\DeclarePairedDelimiter{\pbrac}{\lbrack}{\rbrack}
\DeclarePairedDelimiter{\Spanpr}{\lparen}{\rparen}

\newcommand{\n}{{\bm{n}}}

\newcommand{\olUPlus}{\ol{U}_+ \setminus \set{0} \times \Gamma_+}

\newcommand{\FMOC}{{\widetilde{F}}}

\renewcommand{\epsilon}{\varepsilon}
\newcommand{\eps}{\ensuremath{\varepsilon}}

\newcommand{\Cal}[1]{\ensuremath{\mathcal{#1}}}
\newcommand{\al}{\ensuremath{\alpha}}

\newcommand{\diff}[2]{\frac{ d#1}{d#2}}

\newcommand{\ol}{\overline}

\renewcommand{\le}{\leqslant}
\renewcommand{\ge}{\geqslant}

\newcommand{\MOC}{modulus of continuity\xspace}

\newcommand{\Holder}
    {H\"{o}lder\xspace}

\newcommand{\Poincares}
	{Poincar\'{e}'s\xspace}
	
\newcommand{\Gronwalls}
    {Gr\"{o}nwall's\xspace}

\newcommand{\LinearPaper}[1]{#1 of \cite{Paper1}\xspace}
\newcommand{\LinearPaperCap}[1]{#1 of \cite{Paper1}\xspace}

\crefname{cor}{Corollary}{Corollaries} 
									   
\crefname{lemma}{Lemma}{Lemmas}	       

\crefname{section}{Section}{Sections}
\Crefname{section}{Section}{Sections}

\crefname{appendix}{Appendix}{Appendices}
\Crefname{appendix}{Appendix}{Appendices}

\crefname{theorem}{Theorem}{Theorems}
\Crefname{theorem}{Theorem}{Theorems}

\crefname{prop}{Proposition}{Propositions}
\Crefname{prop}{Proposition}{Propositions}

\crefname{conj}{Conjecture}{Conjectures}
\Crefname{conj}{Conjecture}{Conjectures}

\crefname{definition}{Definition}{Definitions}
\Crefname{definition}{Definition}{Definitions}

\crefname{remark}{Remark}{Remarks}
\Crefname{remark}{Remark}{Remarks}

\crefname{assumption}{Assumption}{Assumptions}
\Crefname{assumption}{Assumption}{Assumptions}

\creflabelformat{enumi}{(#1)}
\creflabelformat{enumii}{(#1)}

\crefformat{equation}{(#2#1#3)}
\crefrangeformat{equation}{(#3#1#4) through (#5#2#6)}
\crefmultiformat{equation}
    {(#2#1#3)}%
    { and~(#2#1#3)}
    {, (#2#1#3)}
    { and~(#2#1#3)}

\newcommand{\aaa}{{\bm{\mathrm{a}}}}
\newcommand{\bbb}{{\bm{\mathrm{b}}}}

\newcommand{\f}{\bm{\mathrm{f}}}

\newcommand{\g}{{\bm{\mathrm{g}}}}
\newcommand{\G}{\bm{\mathrm{G}}}

\renewcommand{\H}{{\bm{\mathrm{H}}}}

\newcommand{\uu}{{\bm{\mathrm{u}}}}   

\newcommand{\vv}{{\bm{\mathrm{v}}}}

\newcommand{\ww}{{\bm{\mathrm{w}}}}
\newcommand{\x}{{\bm{\mathrm{x}}}}
\newcommand{\y}{{\bm{\mathrm{y}}}}

\newcommand{\w}{{\bm{\mathrm{w}}}}
\newcommand{\W}{{\bm{\mathrm{W}}}}

\newcommand{\z}{{\bm{\mathrm{z}}}}
\newcommand{\zz}{{\bm{\mathrm{z}}}}

\newcommand{\UU}{{\bm{\mathrm{U}}}}
\newcommand{\VV}{{\bm{\Cal{V}}}}

\newcommand{\Our}{\cref{e:OurNonlinear}\xspace}
\newcommand{\OurOne}{\cref{e:OurNonlinear}$_1$\xspace}

\newcommand{\OurBCsOnly}{\cref{e:OurNonlinear}$_{4, 5}$\xspace}
\newcommand{\OurWeakBCs}{\cref{e:OurNonlinear}$_{1\text{-}4}$\xspace}
\newcommand{\OurInflow}{\cref{e:OurNonlinear}$_5$\xspace}

\newcommand{\OurLinear}{\cref{e:OurLinear}\xspace}
\newcommand{\OurLinearOne}{\cref{e:OurLinear}$_1$\xspace}

\newcommand{\uSolSpaceZero}{{S^{1, \al}}}
\newcommand{\uSolSpaceN}{{S^{N + 1, \al}}}

\newcommand{\uInitSpaceN}{{C_\sigma^{N + 1, \al}(\Omega)}}

\newcommand{\uInputSpaceZero}{\uSolSpaceZero}
\newcommand{\uInputSpaceN}{\uSolSpaceN}

\newcommand{\USpaceN}{{C^{N + 2, \al}(Q)}}

\newcommand{\USpaceNInf}{{C^{N + 2, \al}(Q_\iny)}}

\newcommand{\vortSpace}{{C^\al(Q)}}
\newcommand{\vortSpaceN}{{C^{N, \al}(Q)}}

\newcommand{\vortSpaceNInit}{{C^{N, \al}(\Omega)}}

\newcommand{\fspaceNInf}{{C^{N + 1, \al}(Q_\iny) \cap C([0, \iny); H_0)}}

\newcommand{\gspaceNInf}{{C^{N, \al}(Q_\iny)}}

\newcommand{\WeakSTwoArgs}[2]{{\mathring{S}^{#1, #2}}}
\newcommand{\WeakSArg}[1]{{\WeakSTwoArgs{#1}{\beta}}}
\newcommand{\WeakS}{{\WeakSArg{N + 1}}}
\newcommand{\WeakSalpha}{{\WeakSTwoArgs{N + 1}{\al}}}

\renewcommand{\time}{{\tau}}
\newcommand{\pos}{{\mbox{\boldmath $\gamma$}}}

\newcommand{\wprt}{\widetilde{\prt}_t}

\newcommand{\PullMarginsIn}[1]%
{%
\begin{minipage}[c]{0.27\textwidth}%
\phantom{x}%
\end{minipage}%
\begin{minipage}[c]{0.46\textwidth}%
\begin{center}
#1%
\end{center}
\end{minipage}%
\begin{minipage}[c]{0.27\textwidth}%
\phantom{x}%
\end{minipage}%
}

\newcommand{\Part}[1]
{
\bigskip
\begin{center}
{\Large \textsc{#1}}
\phantomsection   
\addcontentsline{toc}{chapter}{#1}
\end{center}
}

\allowdisplaybreaks

\begin{document}
\newdateformat{mydate}{\THEDAY~\monthname~\THEYEAR}

\title
	[{E}uler with inflow, outflow and vorticity]
	{The 3{D} {E}uler equations with inflow, outflow and vorticity boundary
	conditions}

\author[G.-M. Gie, J. Kelliher, and A. Mazzucato]
{Gung-Min Gie$^{1}$, James P. Kelliher$^{2}$, and Anna L. Mazzucato$^{\ast,3}$}
\address{$^1$ Department of Mathematics, University of Louisville, Louisville, KY 40292}
\address{$^2$ Department of Mathematics, University of California, Riverside, 900 University Ave., Riverside, CA 92521}
\address{$^3$ Department of Mathematics, Penn State University, University Park, PA 16802}
\email{gungmin.gie@louisville.edu}
\email{kelliher@math.ucr.edu}
\email{alm24@psu.edu}



\begin{abstract}
	The 3D incompressible Euler equations in a bounded domain are most often
	supplemented with impermeable boundary conditions,
	which constrain the fluid to neither enter nor leave the domain.
	We establish well-posedness with inflow, outflow of velocity when either
	the full value of the
	velocity is specified on inflow, or only the normal component is specified
	along with the vorticity (and an additional constraint). We derive
	compatibility conditions to
	obtain regularity in a \Holder space with prescribed arbitrary index,
	and allow multiply connected domains.
	Our results apply as well to impermeable boundaries,
	establishing higher regularity of solutions in \Holder spaces.
\end{abstract}

\maketitle

\markleft{G.-M. GIE, J. KELLIHER, AND A. MAZZUCATO}

\vspace{-0.8em}

\begin{center}
Compiled on {\dayofweekname{\day}{\month}{\year} \mydate\today} at \currenttime \\
\end{center}

\bigskip

\vspace{-0.8em}

\centerline{
\begin{minipage}{0.8\textwidth}
\begin{center}
\renewcommand\contentsname{}
\setcounter{tocdepth}{1}		
\tableofcontents
\end{center}
\end{minipage}
}
\newpage

%
%
\Part{Part I: Overview}

\vspace{-0.25em}

\section{Introduction}\label{S:Introduction}

\noindent Let $\Omega$ be a bounded domain in $\R^3$, possibly multiply connected, having a boundary that is at least $C^2$ regular. We define $\n$ to be the outward unit normal vector to the boundary, $\Gamma := \prt \Omega$, and follow the convention that for any vector field $\vv$,
\begin{align}\label{e:ComponentConvention}
	v^\n := \vv \cdot \n, \quad
	\vv^\n := v^\n \n, \quad
	\vv^\BoldTau := \vv - \vv^\n
		\text{ on } \Gamma.
\end{align}

Fixing $T > 0$, the Euler equations on $Q := (0, T) \times \Omega$ can be written,
\begin{align}\label{e:EOnly}
	\begin{cases}
		\prt_t \uu + \uu \cdot \grad \uu + \grad p = \f
			& \text{in } Q, \\
		\dv \uu = 0
			& \text{in } Q, \\
		\uu(0) = \uu_0
			& \text{on } \Omega.
	\end{cases}
\end{align}
Here, $\uu$ is the velocity field of a constant-density incompressible fluid, $p$ its scalar pressure, $\f$ the divergence-free external force tangential to the boundary, and $\uu_0$ the initial velocity.

To complete the system of equations in \cref{e:EOnly} we impose inflow, outflow boundary conditions in the spirit of \cite{AKM}. We partition the boundary $\Gamma$ into three portions, $\Gamma_+$, $\Gamma_-$, and $\Gamma_0$, corresponding to inflow, outflow, and impermeability, respectively. Each portion consists of a finite number of components (with $\Gamma_0 = \emptyset$ or $\Gamma_0 = \Gamma$ allowed---see \cref{R:ImpermeableOnly}). We fix a vector field $\UU$ on $[0, T] \times \Omega$ and assume that
\begin{align}\label{e:UBoundCond}
	U^\n < 0 \text{ on } \Gamma_+, \qquad
	U^\n > 0 \text{ on } \Gamma_-, \qquad
	U^\n = 0 \text{ on } \Gamma_0.
\end{align}
We then define inflow, outflow boundary conditions as
\begin{align}\label{e:OurBCs}
	\begin{cases}
		u^\n = U^\n
			& \text{on } [0, T] \times \Gamma, \\
		\uu = \UU
			& \text{on } [0, T] \times \Gamma_+.
	\end{cases}
\end{align}
We also impose on $\UU$ the constraint that
$
	\int_{\Gamma_+} U^\n = - \int_{\Gamma_-} U^\n
$,
required to allow $\dv \uu = 0$.

We choose to impose inflow, outflow boundary conditions in terms of a vector field $\UU$ defined on all of $\Omega$ because it will be productive for us to view $\UU$ as a background flow as done in \cite{TemamWang2002,Pet06,GHT11}. We will also choose $\UU$ to be divergence-free as, shown can be done in \cite{GHT11}, as this will be convenient, though not strictly necessary.

Defining the vorticity,
\begin{align*}
	\bomega := \curl \uu,
\end{align*}
applying $\curl$ to both sides of \cref{e:EOnly}$_1$ yields the vorticity equation,
\begin{align}\label{e:vortEq}
	\prt_t \bomega + \uu \cdot \grad \bomega - \bomega \cdot \grad \uu = \g := \curl \f.
\end{align}
It follows from \cref{e:vortEq} that the vorticity is transported and stretched (pushedforward) by the flow map for $\uu$ (when $\g \equiv 0$).

In particular, the vorticity is brought into the domain from the inflow boundary, making inflow, outflow substantially more difficult to treat than impermeable boundaries: the mechanism for generating vorticity on the inflow boundary must be understood and controlled. This is a key reason for using \Holder spaces, as there is no loss of regularity of the trace of the vorticity on the boundary over that in the domain.

Higher regularity solutions for inflow, outflow boundary conditions are employed, for instance, in Prandtl-type boundary layer expansions (such as \cite{TemamWang2002,GHT11} and work in progress of the authors). The validity of such expansions for inflow, outflow boundary conditions results from a stability mechanism of injection, suction in boundary layers. These applications were the original motivation for this work: because of this, in \cref{A:CompCondsSpecialCase} we give the explicit form of the compatibility conditions for those works.

The system of equations we study, then, are \cref{e:EOnly} with \cref{e:OurBCs}:
\begin{align}\label{e:OurNonlinear}
	\begin{cases}
		\prt_t \uu + \uu \cdot \grad \uu + \grad p = \f
			& \text{in } Q, \\
		\dv \uu = 0
			& \text{in } Q, \\
		\uu(0) = \uu_0
			&\text{on } \Omega, \\
		u^\n = U^\n
			& \text{on } [0, T] \times \Gamma, \\
		\uu = \UU
			& \text{on } [0, T] \times \Gamma_+.
	\end{cases}
\end{align}

We can state the main result of this paper informally as follows, where \textbf{throughout, we fix} $\bm{\al \in (0, 1)}$:

\begin{theorem*}[Informal statement of main result]\label{T:WellPosednessRough}
	Assume that for some integer $N \ge 0$, $\uu_0$ is a divergence-free vector field
	in the classical \Holder space $C^{N + 1, \al}(\Omega)$, satisfies \cref{e:OurBCs},
	and satisfies a compatibility condition to be described below.
	There is a $T > 0$ such that there exists a unique solution
	to \Our with $\curl \uu(t) \in C^{N , \al}(\Omega)$ for
	all $t \in [0, T]$.
\end{theorem*}

We state our main result rigorously in \cref{T:WellPosedness}, but to do so, we must define the function spaces in which we will work, determine proper conditions on the forcing, and determine the required compatibility conditions. It will be helpful, however, to first explain how the boundary conditions in \OurBCsOnly arise.

\subsection*{Possible boundary conditions}
Being motivational, we will argue somewhat heuristically. Some of what we observe will echo observations in \cite{OligerSundstrom1978}---in particular, the comments on an ``open boundary'' in Section 2 of \cite{OligerSundstrom1978} on the linearized compressible Euler equations and in Section 3 of \cite{OligerSundstrom1978} on the linearized incompressible Euler equations.

By taking the divergence of \cref{e:EOnly}$_1$, the pressure can be recovered from the velocity field by
\begin{align}\label{e:TruePressureEq}
	\begin{cases}
		\Delta p = - \grad \uu \cdot (\grad \uu)^T
			&\text{in } \Omega, \\
		\grad p \cdot \n = \prt_t \uu \cdot \n - (\uu \cdot \grad \uu) \cdot \n
			&\text{on } \Gamma.
	\end{cases}
\end{align}
But also, starting from the Gromeka-Lamb form of the Euler equations, one can easily show (see \cref{P:omegaOnGammaPlus}) that any $(\uu, p)$ that satisfies \cref{e:EOnly} must satisfy, on $\Gamma$, the identity,
\begin{align}\label{e:pTau}
	\begin{split}
	u^\n \bomega^\BoldTau
		&= \brac{
				-\prt_t \uu^\BoldTau - \grad_\Gamma \pr{p + \frac{1}{2} \abs{\uu}^2}
			 	 	+ \f}^\perp
				+ \curl_\Gamma \uu^\BoldTau \uu^\BoldTau, \qquad
	\omega^\n
		= \curl_\Gamma \uu^\BoldTau.
	\end{split}
\end{align}
Here, $\vv^\perp = \n \times \vv$ is the tangential vector field $\vv$ on $\Gamma$ rotated 90 degrees counterclockwise around the normal vector $\n$ when viewed from outside $\Omega$, $\grad_\Gamma$ is the tangential derivative, and $\curl_\Gamma$ is the curl operator on the boundary. (\cref{A:BoundaryCalculations} gives details.)

If we impose impermeable boundary conditions, $u^\n \equiv 0$ on $\Gamma$, then the vorticity term disappears in \cref{e:pTau} and there is no constraint on the vorticity. But on portions of the boundary where $\uu \cdot \n$ does not vanish, \cref{e:pTau} gives a relation among $\bomega$, $\grad_\Gamma p$, $\uu \cdot \n$, and $\uu^\BoldTau$ on the boundary. At the same time, \cref{e:TruePressureEq} gives a (global) relation between $\uu \cdot \n$ (via its time derivative) and $p$. At the risk of oversimplifying, together, \cref{e:TruePressureEq,e:pTau} give two relations among four quantities, so we must have an independent means of determining two of them so as to obtain the value of the other two.

To better understand the consequences of \cref{e:pTau}, we turn to the vorticity equation, \cref{e:vortEq}. For impermeable boundary conditions, one can express a Lagrangian solution to the vorticity equation by introducing the flow map, $\eta(t_1, t_2; \x)$,  for $\uu$. This flow map gives the position that a particle at $\x \in \ol{\Omega}$ at time $t_1$ will be as it moves, forward or backward, along the flow line to time $t_2$. Given the flow map, $\bomega(t, \x) := \grad \eta(0, t; \eta(t, 0; \x)) \bomega_0(\eta(t, 0; \x))$ is a Lagrangian solution. (In 2D, it would be  $\omega(t, \x) := \omega_0(\eta(t, 0; \x))$.) This works, because $\eta$ maps any point in $\ol{\Omega}$ to another point in $\ol{\Omega}$, so one can always evaluate $\bomega_0(\eta(t, 0; \x))$.

At points on the boundary at which $\uu \cdot \n < 0$, however, the flow lines enter the domain, and we must have a way of determining, or \textit{generating}, the vorticity so that it can be transported into the domain. From \cref{e:pTau}, we have at such inflow points,
\begin{align}\label{e:omegaNonVanishingPoint}
	\begin{split}
	\bomega^\BoldTau
		&:= \frac{1}{u^\n}
			\brac{
				-\prt_t \uu^\BoldTau - \grad_\Gamma \pr{p + \frac{1}{2} \abs{\uu}^2}
			 	 	+ \f}^\perp
				+ \frac{1}{u^\n}
					\curl_\Gamma \uu^\BoldTau \uu^\BoldTau, \\
	\omega^\n
		&:= \curl_\Gamma \uu^\BoldTau.
	\end{split}
\end{align}

There is, however, another constraint: As vorticity is generated on the boundary and pushed forward into the domain, the resulting vorticity must lie in the range of the curl; that is, the vector field that results must actually itself be the vorticity of some divergence-free vector field. In 2D, this is automatic, because vorticity is simply a scalar field. But in 3D, vorticity is in the range of the curl only if it is divergence-free and has vanishing fluxes across each boundary component.

Taking the divergence of \cref{e:vortEq} leads, after some calculation, to the conclusion that $\prt_t \dv \bomega + \uu \cdot \grad \dv \bomega = 0$; that is, the divergence of the vorticity is transported by the flow. Since $\dv \bomega_0 = \dv \curl \uu_0 = 0$, we need only show that $\dv \bomega = 0$ at inflow points on the boundary. But another calculation gives that on $\Gamma$,
\begin{align*}
	u^\n \dv \bomega
		&= \g \cdot \n - \prt_t \bomega \cdot \n 
			- \dv_\Gamma [\omega^\n \uu^\BoldTau - \uu^\n \bomega^\BoldTau],
\end{align*}
where $\dv_\Gamma$ is the divergence operator on the boundary (see \cref{A:BoundaryCalculations}). This leads to the constraint,
\begin{align}\label{e:VorticityConstraint}	
	\prt_t \bomega \cdot \n 
		+ \dv_\Gamma [\omega^\n \uu^\BoldTau - \uu^\n \bomega^\BoldTau]
		- \g \cdot \n = 0
\end{align}
at inflow points. These calculations are all formal, but are worked out rigorously in detail in \LinearPaper{Section 6}.

Insuring that \cref{e:VorticityConstraint} holds at inflow points is an issue that must be addressed regardless of the manner in which vorticity is generated on the boundary. (Dealing with the external fluxes vanishing is relatively straightforward, and is also treated in \LinearPaper{Section 6}.)

Now, without prescribing at least the sign of $\uu \cdot \n$ on the boundary, we would have to determine the regions of inflow dynamically. To avoid this considerable difficulty, we impose Dirichlet conditions for $\uu \cdot \n$ on all of $\Gamma$. It remains to select a second condition that allows the constraints in \cref{e:TruePressureEq,e:pTau} to be met.
Using the value of $\uu^\BoldTau$ at inflow points as the second condition leads to \OurBCsOnly.

It is not obvious, but specifying the full velocity field at inflow points on the boundary and generating the vorticity at inflow points via \cref{e:omegaNonVanishingPoint} automatically gives \cref{e:VorticityConstraint}, as we show in \cref{P:H0CondSatisfied}. Hence, \cref{e:VorticityConstraint} does not represent an additional constraint.

Another possibility, which we also treat in \cref{T:VorticityBCs}, is to specify the value of $\uu \cdot \n$ and, at inflow points, the value of $\bomega$---so-called \textit{vorticity boundary conditions}. This, however, does not lead to the constraint in \cref{e:VorticityConstraint} being automatically satisfied; rather, we must impose a restriction on our choice of $\bomega$. It is not clear how to do this in greatest generality, but by requiring that the prescribed vorticity be tangential to the inflow boundary, we obtain well-posedness nearly for free from the technology we develop to handle inflow, outflow boundary conditions.

(One could also choose to use an independent relation between $\bomega$ and $\uu^\BoldTau$ for the second condition. This was done by Chemetov and Antontsev \cite{ChemetovAntontsev2008} for 2D weak solutions in vorticity form, without uniqueness, for Navier friction boundary conditions.)

Once we have points on the boundary at which $\uu \cdot \n < 0$, we must have other points at which $\uu \cdot \n > 0$ else the fluid could not be incompressible. Hence, we must have
$
	0
		= \int_\Gamma \uu \cdot \n
		= \int_\Omega \dv \uu
$.
Reflecting upon \cref{e:omegaNonVanishingPoint}, it would be very difficult to handle $\uu \cdot \n$ vanishing at a point or, even worse, changing sign, especially to obtain classical solutions with higher regularity, which is our intent. To avoid this, each boundary component must have $\uu \cdot \n$ strictly negative (inflow), strictly positive (outflow), or vanish identically. If a component has inflow, then at least one other component must have outflow.

\subsection*{Contrast with the analytic setting}
Motivated in part by the results in \cite{Paper1}, the authors of \cite{KOS} have recently established well-posedness of the Euler equations with inflow and outflow in the analytic category. More precisely, they show existence and uniqueness of solutions in certain spaces of functions that are analytic in the tangential direction near the boundary and otherwise belong to a Sobolev space with sufficiently high index of regularity.

By working with analytic norms they can absorb the loss of a derivative at the boundary in the weak formulation of the Euler equation, this loss arising from the non-homogeneous inflow and outflow boundary conditions. Because of this, they can employ a sequence of approximating solutions based on velocity and pressure with no need to treat the behavior of the vorticity. This also allows them to directly recover the pressure from the velocity in their analytic spaces. As a consequence, they need to prescribe only the normal component of the velocity at the boundary: in their setting, no compatibility conditions are necessary, as the relation in \cref{e:omegaNonVanishingPoint}, while it must hold, never enters into the analysis.

At the same time, working with Sobolev norms they can derive suitable {\em a priori} estimates by higher energy estimates. Then, the unique solution is obtained directly from the velocity-pressure formulation of the Euler equations via a Picard iteration.

By contrast, we use the transport of vorticity to establish existence as in \cite{AKM}, bypassing the loss of derivative at the boundary. The vorticity generated at the inflow boundary, however, requires both knowledge of the pressure and of the tangential component of the velocity field, as we can see in \cref{e:omegaNonVanishingPoint}. Hence, we must impose an additional constraint, satisfied by imposing full inflow boundary conditions as well as compatibility conditions among the data.

\subsection*{Function spaces}
Returning to stating our main result rigorously, we define the function spaces for our solutions.
For any $N \ge 0$ we define the affine spaces,
\begin{align}\label{e:PrimarySpaces}
	\begin{split}
	\uInitSpaceN
		&:= \set{\uu \in C^{N + 1, \al}(\Omega) \colon \dv \uu = 0, \uu \cdot \n = U^\n(0)
					\text{ on } \Gamma}, \\
	\uSolSpaceN
		&:= \set{\uu \in C^{N, \al}(Q) \colon \curl \uu \in C^{N, \al}(Q), \,
                        \prt_t^{N + 1} \uu \in L^\iny([0, T]; C^\al(\Omega)), \\
                &\qquad
                        \dv \uu = 0, \uu \cdot \n = U^\n
					\text{ on } [0, T] \times \Gamma}, \\
        \norm{\uu}_\uSolSpaceN
            &:= \norm{\uu}_{C^{N, \al}(Q)}
                + \norm{\curl \uu}_{C^{N, \al}(Q)}
                + \norm{\prt_t^{N + 1} \uu}_{L^\iny([0, T]; C^\al(\Omega))}.
	\end{split}
\end{align}
Although these spaces depend on $\UU$, for notational simplicity, we drop the $\UU$, as it is fixed.

We also define
\begin{align}\label{e:HSpace}
	H
		&:= \set{\uu \in L^2(\Omega)^3 \colon \dv \uu = 0, \, \uu \cdot \n = 0
				\text{ on } \Gamma}
		= H_0 \oplus H_c,
\end{align}
where
\begin{align}\label{e:HcH0}
	H_c &:= \set{\vv \in H: \curl \vv = 0}, \quad
	H_0 := H_c^\perp.
\end{align}
Let $P_{H_c}$ be the projection operator from $H$ to $H_c$. We call $P_{H_c} \uu$ the \textit{harmonic} component of $\uu$.

We define the boundary values (via $\UU$) and the forcing $\f$ for all time on $Q_\iny := [0, \iny) \times \Omega)$. We will prove existence only for short time.

\begin{definition}\label{D:NReg}
	We say that the data has regularity $N$ for an integer $N \ge 0$ if
	\begin{itemize}
		\item
			$\Gamma$ is $C^{N + 2, \al}$,
			$\f \in \fspaceNInf$;
			
		\item
			$\UU \in \USpaceNInf$, $\dv \UU = 0$, and \cref{e:UBoundCond} holds;
		
		\item
			$U_{min} := \min \set{\abs{U^\n(t, \x)} \colon (t, \x)
				\in [0, \iny) \times \Gamma_+} > 0$;

		\item
			$\uu_0 \in \uInitSpaceN$, $\uu_0^\BoldTau = \UU_0^\BoldTau$ on $\Gamma_+$.
	\end{itemize}
\end{definition}

We assumed that $\UU$ has one more derivative than $\uu$, as explained in \cref{R:HHigherRegularity}.

\subsection*{Compatibility conditions}

The vorticity generated at the inflow boundary is carried by the flow into the interior; at the same time, the flow pushes the initial vorticity forward in time. The interaction between these two sources of vorticity may potentially lead to a singularity. The main thrust of this work is to show that it is possible to avoid such singularities, at least for short time, by imposing suitable conditions on the data. We refer to these conditions as \textit{compatibility conditions}, satisfying two primary principles:
\begin{enumerate}
	\item
		They depend only upon the initial data, $\UU$, and $\f$.
		
	\item
		They are compatible with being a solution to \Our; that is, a solution to
		\Our could, in principle, satisfy them.
\end{enumerate}
The conditions we develop will ensure regularity of the solution for short time. It remains an open question whether a regular solution persists for all time even in 2D.

Given $\uu$ with data regularity $N$ for some $N \ge 0$, we define the $N^{th}$ compatibility condition,
\begin{align}\label{e:condN}
	\begin{split}
	&\cond_{-1}:
		\uu_0^\BoldTau = \UU_0^\BoldTau
			\text{ on } \Gamma_+, \\
	&\cond_N: \cond_{N - 1} \text{ and } 
		\prt_t^{N + 1} \UU^\BoldTau|_{t = 0}
			= \wprt^{N + 1} \uu_0^\BoldTau
			\text{ on } \Gamma_+.
	\end{split}
\end{align}
For integers $n \ge 0$, we define $\wprt^n \uu_0$ inductively by setting $\wprt^0 \uu_0 = \uu_0$, while for $n \ge 1$, we take the time derivative of $\wprt^{n - 1} \uu$ at time zero and replace each instance of $\prt_t \uu$ in the resulting expression by $- \uu_0 \cdot \grad \uu_0 - \grad p^0 + \f(0)$. Here, $p^0$ is the value the pressure would have at time zero if $\uu$ actually solved \Our; that is, $p^0$ is the solution to the system in \cref{e:TruePressureEq} at time zero:
\begin{align}\label{e:p0}
	\begin{cases}
		\Delta p^0 = - \dv (\uu_0 \cdot \grad \uu_0)
			&\text{in } \Omega, \\
		\grad p^0 \cdot \n = -\prt_t U^\n(0) - \uu_0 \cdot \grad \uu_0
			&\text{on } \Gamma.
	\end{cases}
\end{align}
We give a more complete account of the compatibility conditions in \cref{S:CompConds}.

For $N = 0$, \cref{e:condN} is the compatibility condition in (1.10), (1.11) of Chapter 4 of \cite{AKM}:
\begin{align*}
	&\cond_0: 
		\prt_t \UU^\BoldTau|_{t = 0}
			= [- \uu_0 \cdot \grad \uu_0 - \grad p^0 + \f(0)]^\BoldTau
			\text{ on } \Gamma_+.	
\end{align*}

\subsection*{Main result}
We can now rigorously state the main result of this paper as follows:

\begin{theorem}\label{T:WellPosedness}
	Assume the data has regularity $N$ for some integer $N \ge 0$ as in \cref{D:NReg}
	and satisfies $\cond_N$ of \cref{e:condN}. There is a $T > 0$ such that there exists a solution
	$(\uu, p)$ to \Our with
	$\uu \in \uSolSpaceN$ and $\grad p$ in $L^\iny([0, T]; C^{N, \al}(\Omega))$,
	which is unique up to an additive constant for the pressure.
	If $N \ge 1$, $\grad p$ is also in $C^{N - 1, \al}(Q)$.
\end{theorem}

\begin{remark}\label{R:TEstimate}
It follows from the proof of \cref{T:WellPosedness} that $T$ is bounded below by a continuous, increasing function of norms of $(U^\n)^{-1}$, $\UU$, $\f$, and $\uu_0$ in the spaces appearing in \cref{D:NReg}. The explicit form of the estimate is, however, involved and may not be optimal. 
It also follows from the proof that the higher regularity of $\UU$ and so $\uu$ on the inflow boundary leads to higher time regularity of the pressure near the inflow boundary that need not, however, extend through the whole domain.
\end{remark}

\subsection*{Vorticity boundary condition}

We also consider solutions  $(\uu, p, \zz)$ to the Euler equations with vorticity boundary conditions, where the value of the vorticity on the inflow boundary is given by a function $\H$ on $[0, T] \times \Gamma_+$:
\begin{align}\label{e:OurVorticity}
	\begin{cases}
		\prt_t \uu + \uu \cdot \grad \uu + \grad p = \f + \zz
			& \text{in } Q, \\
		\dv \uu = 0
			& \text{in } Q, \\
		\uu(0) = \uu_0
			&\text{on } \Omega, \\
		u^\n = U^\n
			& \text{on } [0, T] \times \Gamma, \\
		\curl \uu = \H
			& \text{on } [0, T] \times \Gamma_+.
	\end{cases}
\end{align}
Here, $\z \in H_c$ is an harmonic vector field.

As stated, this system is not yet complete, as we must have a means of determining the harmonic component $P_{H_c}$ of $\uu$. We can do this either (1) directly, by setting $P_{H_c} \uu = \uu_c$ for some harmonic vector field $\uu_c$ and letting $\z$ be obtained as part of the solution or (2) indirectly, by prescribing the value of $\z$ and obtaining $P_{H_c} \uu$ as part of the solution. In \cref{T:VorticityBCs} we choose (1), as it allows for the uniqueness of solutions.

\begin{theorem}\label{T:VorticityBCs}
	Fix $\uu_c \in C^{N + 1, \al}(Q) \cap C([0, T]; H_c)$.
	Assume that the data has regularity $N$ for some integer $N \ge 0$ as in \cref{D:NReg},
	that $\cond_N$ holds, and that $\uu_c(0) = P_{H_c} P_H \uu_0$.
	Also assume that $\H \in C^{\max\set{N, 1}, \al}([0, T] \times \Gamma_+)$
	and
	\begin{align}\label{e:RangeOfCurlCondSpecial}
		H^\n
			= 0, \quad
		\dv_\Gamma[U^\n \H^\BoldTau] + \curl \f \cdot \n
			= 0
			\text{ on } [0, T] \times \Gamma_+.
	\end{align}
	There is a $T > 0$ such that there exists a solution
	$(\uu, p, \zz)$ in $\uSolSpaceN \times L^\iny([0, T]; C^{N, \al}(\Omega))
	\times (C^{N + 1, \al}(Q) \cap C([0, T]; H_c)$ to \cref{e:OurVorticity}
	supplemented with the condition that
	\begin{align*}
		P_{H_c} \uu = \uu_c \text{ on } [0, t] \times \Omega.
	\end{align*}
	If $N \ge 1$,
	$\grad p$ is also in $L^\iny([0, T]; C^{N, \al}(\Omega))$	
	and the solution is unique up to an additive constant for the pressure.
	In addition, $\z(0) = 0$.
\end{theorem}

The condition in \cref{R:VortCondExplained} reflects the constraint in  \cref{e:VorticityConstraint}.

\subsection*{Approaches to well-posedness of Euler equations}

There are many proofs of well-posedness of the Euler equations taking different approaches. Most such proofs in \Holder spaces in a 3D domain with boundary, including this paper, and many in the whole space or a periodic domain, follow in the tradition of McGrath \cite{McGrath1966,McGrath1967} and Kato \cite{Kato1967}, in which the solution is obtained as a fixed point of an operator $A$ derived from a linearization of the Euler equations, employing Schauder's fixed point theorem. A notable exception is \cite{BardosFrisch1976}, which establishes well-posedness in $C^{N + 1, \al}(\Omega)$ for impermeable boundary conditions (finite time for 3D) using a more direct iteration scheme, obtaining a contraction mapping.

For inflow, outflow boundary conditions, the Schauder's fixed point theorem approach was taken in Chapter 4 of \cite{AKM}, which establishes \cref{T:WellPosedness} for $N = 0$ and simply connected domains. The operator $A$ is derived from a linearization of the vorticity equation \cref{e:vortEq} with prescribed values on the inflow boundary (see \cref{D:A}). This leads to linear compatibility conditions based on vorticity, whereas the nonlinear boundary conditions are based on the velocity. In fact, one challenge is to ensure that the nonlinear compatibility conditions at the level of the velocity imply the linear ones at the level of the vorticity.

To handle inflow, outflow boundary conditions, the authors of \cite{AKM} make many adaptations to the Kato-McGrath approach, but we would identify their two key innovations as the following:
\begin{itemize}
	\item
		They obtain estimates on the operator $A$
		under the simple linear compatibility condition that on the inflow boundary,
		the vorticity matches the prescribed inflow vorticity at time zero
		(akin to the Rankine-Hugoniot condition).

	\item
		They show how to achieve the needed regularity of the inflow vorticity
		from the pressure.
\end{itemize}

\subsection*{What is novel in our approach}

For $N \ge 1$, several complications arise. We can still use the same operator $A$ as in \cite{AKM}, but now the linear compatibility conditions becomes more involved (see \cref{e:lincond01,e:lincondN}). This linear issue was resolved in \cite{Paper1}, but deriving and relating the nonlinear compatibility conditions to the linear ones remained a significant challenge, which we address here. (In \cref{R:Paper1}, we say a few words about the analysis in \cite{Paper1}, after we have introduced a number of the constructs involved.)

Moreover, while for $N = 0$, the linear compatibility condition implies the nonlinear compatibility condition, this is no longer the case for $N \ge 1$. To address this, we must restrict the domain of the operator $A$ by imposing an additional condition on the time derivative of the initial velocity (as in \cref{e:DomNA}) and show that, in fact, the resulting domain is nonempty.

The estimates on the operator $A$ that result become much more complex for the higher regularity we treat here. This is in contrast to proving existence in the full space or a periodic domain, where one can bootstrap as in Section 4.4 of \cite{MB2002}, which takes advantage of the simple form of Biot-Savart kernel for the full space. And in 2D, where the vorticity equation has no stretching term, one can bootstrap as Marchioro and Pulvirenti do in \cite{MP1994} (which originates in their earlier text \cite{MP1984}). 
Instead, we must obtain existence directly in the higher-regularity spaces: the resulting estimates  are therefore much more intricate than the $N = 0$ case.

\subsection*{Other Prior work} 

In addition to \cite{AKM}, we also drew ideas from \cite{Koch2002}, which proves well-posedness of the 3D Euler equations for impermeable boundary conditions in \Holder spaces (the equivalent of our $N = 0$ regularity). We mention also the work of Petcu \cite{Pet06}, who presents a version of the argument in Chapter 4 of \cite{AKM}, specializing it to a 3D channel with a constant $\UU$, which simplifies and makes clearer some of the arguments in \cite{AKM}. 

Section 1.4 of \cite{Mamontov2010} contains an extensive survey of results, both 2D and 3D, related to the problem we are studying here. We also mention the 2D work of Boyer and Fabrie \cite{BoyerFabrie2006,BoyerFabrie2013} and the recent works \cite{BravinSueur2021,NoisetteSueur2021}.

Vorticity boundary conditions were studied in 2D by Yudovich in \cite{Yudovich1964V}. We refer in addition to the historical comments in Section 1.4 of \cite{Mamontov2010} concerning partial results in 3D. Both vorticity and pressure-based boundary conditions with inflow, outflow were treated in \cite{Zaj1980,Zaj1982,Zaj1984a,Zaj1984b}.

Finally, we mention the works \cite{TemEuler75, KatoLai1984}, which give well-posedness of solutions to the Euler equations with impermeable boundaries in Sobolev spaces.

\subsection*{Structure of this paper}

This paper consists of three parts, along with three appendices.

In Part I, following this introduction, we begin in \cref{S:Paper1} by summarizing results from \cite{Paper1} on the linearization of \Our, a key tool at the heart of all of our arguments. In \cref{S:InflowVorticity} we detail how vorticity is generated on the inflow boundary and then define the operator $A$. In \cref{S:CompConds}, we explore in-depth the nonlinear compatibility conditions $\cond_N$ as they apply to \Our and their counterparts for the linearized equations. We then give the proof of our main result, \cref{T:WellPosedness}, in \cref{S:Proof}. This proof relies upon three propositions, \crefrange{P:InvariantSet}{P:WeakVorticityIsVelocity}: the rest of the paper is devoted to proving these propositions.

In Part II, we develop some properties related to the function space $\uSolSpaceN$, then present identities and estimates on the flow map, on the vorticity generated on the boundary, and on the pressure.

In Part III, we use results primarily from the second part to prove \cref{P:InvariantSet}, then leverage it to obtain \cref{P:AContinuous}. We also give the proof of \cref{P:WeakVorticityIsVelocity}. In the final section of this part, we describe how \cref{T:VorticityBCs} follows from the estimates obtained in Part II.

\cref{A:HolderLemmas} contains a number of estimates in \Holder spaces, some very standard, some specific to this paper. In \cref{A:BoundaryCalculations} we construct a convenient coordinate system in an $\eps$-neighborhood of $\Gamma_+$. We use this system to develop properties of the operators $\grad_\Gamma$, $\dv_\Gamma$, and $\curl_\Gamma$ we use in the body of the paper. This allows us to treat the various calculations on the boundary in a coordinate-free manner, which makes the calculations more transparent. Finally, in \cref{A:CompCondsSpecialCase}, we discuss the compatibility conditions in the special case in which $\UU^\BoldTau \equiv 0$ and $U^\n$ is constant along $\Gamma_+$ (as occurs in \cite{TemamWang2002,GHT11}).

We have structured this paper so as to allow the reader to grasp the overall structure of the proof of \cref{T:WellPosedness} without it being obscured by the many technical details. It is possible to read only Part I and get the gist of the proof. A more in-depth reading would involve at least examining the key pressure estimates in \cref{S:Pressure} and a reading of \cite{Paper1}, to understand how the linear compatibility conditions arise.

\subsection*{On notation} Our notation, while fairly standard, has a few subtleties. If $M$ is a matrix, $M^i_n$ refers to the entry in row $i$ of $M$, column $n$; $v^i$ refers to the $i^{th}$ entry in the vector $\vv$, which we always treat as a column vector for purposes of multiplication. If $M$ and $N$ are matrices of the same dimensions then $M \cdot N := M^i_n N^i_n$, where here, and elsewhere, we use implicit sum notation. If $\uu$ and $\vv$ are vectors then the matrix $\uu \otimes \vv$ has components $[\uu \otimes \vv]^i_n = u^i v^n$.

We define the divergence of a  matrix row-by-row, so $\dv M$ is the column vector with components $[\dv M]^i = \prt_\ell M^i_\ell$. Hence, $[\dv [\uu \otimes \vv]]^i = \dv [\uu \otimes \vv]^i = \prt_\ell (u^i v^\ell)$, where $\prt_\ell$ is the derivative with respect to the $\ell^{th}$ spatial variable. The notation $\grad$ means the gradient with respect to the spatial variables only; by $D$ we mean the gradient with respect to all variables, time and
space. When applied to the flow map $\eta(t_1, t_2, \x)$, we write $\prt_{t_1} \eta$, $\prt_{t_2} \eta$ to mean the derivative with respect to the first, second time variable. Finally, for vector fields $\uu$ and $\vv$, we will interchangeably write $\uu \cdot \grad \vv$ and $\grad \vv \, \uu$, as they both are vectors with $i^{th}$ component $u^m \prt_m v^i$.

For any tangent vector field $\vv$ on $\Gamma$, $\vv^\perp = \n \times \vv$ is the tangent vector field $\vv$ on $\Gamma$ rotated 90 degrees counterclockwise around the normal vector $\n$ when viewed from outside $\Omega$.
We write the gradient, divergence, and curl operators on the boundary as $\grad_\Gamma$, $\dv_\Gamma$, and $\curl_\Gamma$, as defined in \cref{A:BoundaryCalculations}.

When we write that a function is in a \Holder space $C^{k, \al}(U)$ (defined in \cref{S:BSLaw}) we mean not just that it has the given regularity but that its norm is finite. Since a function in $C^{k, \al}(U)$ extends uniquely to a function in $C^{k, \al}(\ol{U})$, this will rarely have an impact.

%
%
\section{The linearized problem}\label{S:Paper1}

\noindent
The linearized Euler equations corresponding to the vorticity form of \OurOne are
\begin{align}\label{e:OurLinear}
	\begin{cases}
		\prt_t \ol{\bomega} + \uu \cdot \grad \ol{\bomega} - \ol{\bomega} \cdot \grad \uu = \g
			&\text{in } Q, \\
		\ol{\bomega} = \H
			&\text{on } [0, T] \times \Gamma_+, \\
		\ol{\bomega}(0) = \ol{\bomega}_0
			&\text{on } \Omega.
	\end{cases}
\end{align}
Here, $\H$ is given on $[0, T] \times \Gamma_+$, $\ol{\bomega}_0$ is given on $\Omega$, $\uu$ and $\g$ are \textbf{given} on $Q$, and \OurLinear is to be solved for $\ol{\bomega}$. In application, we will set $\ol{\bomega}_0 = \bomega_0 := \curl \uu(0)$, though then $\ol{\bomega}(t) \ne \curl \uu(t)$ in general for $t > 0$.

We employ the following three types of solution to \cref{e:OurLinear}:

\begin{enumerate}[itemsep=0.5em]
	\item
		\textit{Classical Eulerian} or simply \textit{classical} solutions to \OurLinear,
		by which we mean that \OurLinearOne holds pointwise, and each term is continuous
		on $\ol{Q}$.
		
	\item
		\textit{Weak Eulerian solutions}, defined as follows:
		
		\vspace{-0.3em}
		
		\begin{definition}
			\label{D:Weak}
			We say that $\ol{\bomega} \in C(\ol{Q})$ is a weak (Eulerian) solution to \OurLinear if
			$\ol{\bomega} = \H$ on $[0, T] \times \Gamma_+$ pointwise,
			$\ol{\bomega}(0) = \ol{\bomega}_0$, and
			$\prt_t \ol{\bomega} + \dv(\ol{\bomega} \otimes \uu) - \ol{\bomega} \cdot \grad \uu = \g$
			in $\Cal{D'}(Q)$.
		\end{definition}
		\noindent Note that $\ol{\bomega}$ has sufficient
		time and boundary regularity that we do not need to enforce the initial and boundary
		conditions weakly. Also, $\ol{\bomega} \otimes \uu$ is a regular distribution,
		so $\dv(\ol{\bomega} \otimes \uu)$ is a distribution even for $N = 0$.
		
	\item
		\textit{Lagrangian solutions},
		are obtained by pushing forward the vorticity by the flow map, including
		the vorticity generated on the inflow boundary. Because we must
		first introduce some concepts related to this inflow, we defer to
		\cref{D:LagrangianSolution}.
\end{enumerate}

To allow $C^{N, \al}$ solutions to \OurLinear, we must impose linear compatibility condition, $\lincond_N$, defined for $N = 0, 1, \dots$, as follows:
\begin{align}\label{e:lincond01}
	\begin{array}{ll}
		\lincond_0: &\H(0) = \bomega_0 \text{ on } \Gamma_+, \\
		\lincond_1: &\lincond_0 \text{ and }
			\prt_t \H|_{t = 0} 
				=  \bomega_0 \cdot \grad \uu_0 - \uu_0 \cdot \grad \bomega_0 + \g(0)
				\text{ on } \Gamma_+,
	\end{array}
\end{align}
where $\uu_0 := \uu(0)$. In $\lincond_1$, we formally replaced $\prt_t \ol{\bomega}(0)$ with the value it would have were $\ol{\bomega}$ an actual classical solution to \OurLinear at time zero. Continuing this process inductively on higher derivatives defines a formal operator $\wprt^n$ (see \cref{D:wprt} for the details). Then we set, for all $N \ge 2$,
\begin{align}\label{e:lincondN}
	\begin{array}{ll}
		\lincond_N: &\lincond_{N - 1} \text{ and }
			\prt_t^N \H|_{t = 0}
				= \wprt^N \bomega_0
				\text{ on } \Gamma_+.
	\end{array}
\end{align}

\begin{theorem}[\cite{Paper1}]\label{T:LinearExistence}
	Assume that the data has regularity $N$ for some $N \ge 0$ and that
	\begin{itemize}
		\item
			$\g := \curl \f$,
			
		\item
			$\uu \in \uInputSpaceN$,
			
		\item
			$\H \in C^{\max \set{N, 1}, \al}([0, T] \times \Gamma_+)$,
			
		\item
			$\lincond_N$ holds,
			
		\item
			$\bomega_0$ is in the range of the curl, by which we mean that
			$\bomega_0 \in \curl (H^1(\Omega)^3)$,
			
		\item
                The following condition on $\H$ holds on $(0, T] \times \Gamma_+$:
                \begin{align}\label{e:RangeOfCurlCond}
                	\prt_t H^\n
                		+ \dv_\Gamma[H^\n \uu^\BoldTau  - U^\n \H^\BoldTau]
                		- \g \cdot \n
                		= 0.
                \end{align}
	\end{itemize}
	There exists a
	solution $\ol{\bomega}$ to \OurLinear in $C^{N, \al}(Q)$, such that
	$\ol{\bomega}(t)$ is in the range of the curl for all $t \in [0, T]$.
	When $N \ge 1$, the solution is classical Eulerian and unique.
	When $N = 0$, the solution is Lagrangian and is also
	the unique weak Eulerian solution as in \cref{D:Weak} for 
	which $\ol{\bomega}(t)$ is in the range of the curl for all $t \in [0, T]$.

	Moreover, there exists a unique $\vv \in \uSolSpaceN$ with $\curl \vv = \ol{\bomega}$
	and $\vv(0) = \uu_0$, 
	and a unique mean-zero pressure field $\pi$ with
	$\grad \pi$ in $L^\iny([0, T]; C^{N, \al}(\Omega))$ and, if $N \ge 1$,
	also in $C^{N - 1, \al}(Q)$,
	satisfying
	\begin{align}\label{e:LinearVelEq}
		\prt_t \vv + \uu \cdot \grad \vv - \uu \cdot (\grad \vv)^T
		      &=\prt_t \vv + \uu \cdot \ol{\bOmega}
			= -\grad \pi + \f,
	\end{align}
    where the antisymmetric matrix $\ol{\bOmega} = \grad \vv - (\grad \vv)^T$,
    whose nonzero components form the vector $\ol{\bomega}$.
\end{theorem}

\begin{remark}\label{R:LinearFixes}
    In \cite{Paper1}, the stronger condition that $\uu \in C^{N + 1, \al}(Q)$ was assumed, but $\uu \in \uSolSpaceN$ suffices. This is because for the linearized problem it is only the regularity of the flow map for $\uu$ that is needed, and as we will see in \cref{S:FlowMap}, $\uu \in \uSolSpaceN$ is sufficient. It was concluded in \cite{Paper1} that $\vv \in C^{N + 1, \al}(Q)$, but the argument there only gives $\vv \in \uSolSpaceN$ (see \cref{C:AuSolspaceBound}). Also, \cite{Paper1} states the requirement that $\g$ have, in effect, one more derivative of regularity in space than in time, but in fact does not use that additional regularity. Hence, we need only assume that $\g \in C^{N, \al}(Q)$, which follows from \cref{D:NReg}.
\end{remark}

\begin{remark}\label{R:VortCondExplained}
	The condition in \cref{e:RangeOfCurlCond} is the linear analog of \cref{e:VorticityConstraint},
	required to insure that $\ol{\bomega}(t)$
	lies in the range of the curl.
	As applied to the solution of the linearized problem given by \cref{T:LinearExistence},
	\cref{e:RangeOfCurlCond} is a condition on the data, not
	on the solution, since $\uu$ is given. Applied to the fully nonlinear problem
	with $H = \curl \uu$, however,
	the appearance of $\uu^\BoldTau$ in \cref{e:RangeOfCurlCond}
	makes \cref{e:RangeOfCurlCond} a condition on the solution.
	Eliminating the term involving $\uu^\BoldTau$ by requiring that the normal component
	of the vorticity on inflow vanish gives \cref{e:RangeOfCurlCondSpecial},
	which is a condition on the data: $\uu_0$, $\f$, $\UU$, and $\H$.
\end{remark}

In what follows, we will use $\ol{\bomega}$ as a Lagrangian solution, but we will need to estimate $\vv$, which is obtained from the Eulerian solution. Hence, it is crucial that Eulerian and Lagrangian solutions agree.

%
%
\section{Vorticity on inflow, and the operator $A$}\label{S:InflowVorticity}

\noindent At the end of this section, we will use a solution to \OurLinear to define an operator $A$ whose fixed point is a solution to the fully nonlinear Euler equations, \Our. To do this, we first show in \cref{P:omegaOnGammaPlus} how vorticity is generated on the inflow boundary if $(\uu, p)$ is a classical solution to \Our.

\begin{prop}\label{P:omegaOnGammaPlus}	
	Assume that $(\uu, p)$ satisfies \OurOne in a classical sense and let
	$\bomega := \curl \uu$.
	Then on $[0, T] \times \Gamma$,
	\begin{align*}
		\begin{split}
		u^\n \bomega^\BoldTau
			&= \brac{
					-\prt_t \uu^\BoldTau - \grad_\Gamma \pr{p + \frac{1}{2} \abs{\uu}^2}
			 	 		+ \f}^\perp
				 + \pr{\curl_\Gamma \uu^\BoldTau} \uu^\BoldTau, \qquad
		\omega^\n
			= \curl_\Gamma \uu^\BoldTau.
		\end{split}
	\end{align*}
	Here, $\grad_\Gamma$ is the tangential derivative, and $\curl_\Gamma$ is the curl operator
	on the boundary. (See \cref{A:BoundaryCalculations}.)
\end{prop}
\begin{proof}
    As on p. 155 of \cite{AKM}, we start with the Gromeka-Lamb form of
	the Euler equations,
	\begin{align}\label{e:Lamb}
		\prt_t \uu + \grad \pr{p + \frac{1}{2} \abs{\uu}^2}
			- \uu \times \bomega - \f = 0.
	\end{align}
	The equivalence of \cref{e:Lamb} and \OurOne follows from the identity,
	\begin{align}\label{e:NonlinerIdentity}
		\uu \cdot \grad \uu
			&= - \uu \times \bomega + \frac{1}{2} \grad \abs{\uu}^2,
	\end{align}
	which holds as long as $\bomega = \curl \uu$.

	From \cref{L:CrossProduct}
	\begin{align*}
		[\uu \times \bomega]^\BoldTau
			&= u^\n [\bomega^\BoldTau]^\perp - \omega^\n [\uu^\BoldTau]^\perp,
	\end{align*}
	so restricting \cref{e:Lamb} to $[0, T] \times \Gamma_+$, we have
	\begin{align*}
		\prt_t \uu^\BoldTau + \grad_\Gamma \pr{p + \frac{1}{2} \abs{\uu}^2}
			- u^\n [\bomega^\BoldTau]^\perp + \omega^\n [\uu^\BoldTau]^\perp - \f^\BoldTau  = 0.
	\end{align*}
	Hence, since $(\vv^\perp)^\perp = - \vv$ for any tangent vector $\vv$,
	\begin{align*}
		u^\n \bomega^\BoldTau
			&= \brac{-\prt_t \uu^\BoldTau - \grad_\Gamma \pr{p + \frac{1}{2} \abs{\uu}^2}
			 	 + \f^\BoldTau }^\perp + \omega^\n \uu^\BoldTau.
	\end{align*}
	The proof is completed by observing that $\omega^\n = \curl_\Gamma \uu^\BoldTau$
	by \cref{e:curlRestricted}.
\end{proof}

We see from \cref{P:omegaOnGammaPlus} that for a solution to \OurWeakBCs with $\bomega := \curl \uu$, we have
\begin{align}\label{e:omegaIsW}
	\bomega = \W[\uu, p]
		\text{ on } [0, T] \times \Gamma_+,
\end{align}
where, using that $\uu^n = \UU^n \ne 0$ on $\Gamma_+$, $\W[\uu, p]$ is defined on $[0, T] \times \Gamma_+$ by
\begin{align}\label{e:W}
	\begin{split}
	\W^\BoldTau[\uu, p]
		&:= \frac{1}{U^\n}
			\brac{-\prt_t \uu^\BoldTau - \grad_\Gamma \pr{p + \frac{1}{2} \abs{\uu}^2}
					+ \f^\BoldTau}^\perp
			+ \frac{1}{U^\n} \curl_\Gamma \uu^\BoldTau\uu^\BoldTau, \\
	W^\n[\uu, p]
		&:= \curl_\Gamma \uu^\BoldTau.
	\end{split}
\end{align}

Now let $\uu$ be any element of $\uInputSpaceN$, not necessarily a solution of \Our. We seek to define a function $\H$ in $C^{N, \al}([0, T] \times \Gamma_+)$ as a modification of the expression for $\W[\uu, p]$ in such a way that when the data has regularity $N$, at least the following two properties hold:
\begin{itemize}[itemsep=0.2em]			
	\item[(P1)]\label{i:HInitOnly}
		$\H$ at time zero can be defined in terms of $\uu_0$, $\f$,
		and $\UU$ only.
				
	\item[(P2)]\label{i:HHWup}
		If $(\uu, p)$ solves \OurWeakBCs and
		the vorticity $\bomega = \curl \uu$ satisfies $\bomega = \H$
		on $[0, T] \times \Gamma_+$ then
		$(\uu, p)$ satisfies \OurInflow as well---and so solves \Our.
\end{itemize}

We define the function $\H$ for all $N \ge 0$ as done in \cite{AKM} for $N = 0$. First construct a ``regularized pressure'' $p_r$ from $\uu$ as the unique mean-zero solution to
\begin{align}\label{e:pr}
	\begin{cases}
		\Delta p_r = - \dv (\uu \cdot \grad \uu)
			&\text{in } \ol{Q}, \\
		\grad p_r \cdot \n = -\prt_t U^\n - N[\uu]
			&\text{on } [0, T] \times \Gamma,
	\end{cases}
\end{align}
where
\begin{align}\label{e:NuDef}
	\begin{split}
	N[\uu]
		&:=
		\begin{cases}
			(\uu \cdot \grad \uu) \cdot \n
					&\text{on } [0, T] \times (\Gamma_- \cup \Gamma_0), \\
			(\uu \cdot \grad \uu) \cdot \n
				+ \dv_\Gamma (U^\n (\uu^\BoldTau - \UU^\BoldTau))
					&\text{on } [0, T] \times \Gamma_+.
		\end{cases}
	\end{split}
\end{align}

Though not so evident here, using $N[\uu]$ in place of $(\uu \cdot \grad \uu) \cdot \n$ regularizes the pressure, as we will explain in \cref{S:Boundary}. But it is clear from its definition that if $\uu$ satisfies \OurInflow, that is, $\uu^\BoldTau = \UU^\BoldTau$ on $\Gamma_+$, then $N[\uu] = (\uu \cdot \grad \uu) \cdot \n$ on $[0, T] \times \Gamma$, so that $\grad p_r = \grad p$ on $\ol{Q}$, where $p$, given by the system \cref{e:TruePressureEq}, is the ``true pressure.''

Finally, define $\H$ on $[0, T] \times \Gamma_+$ by replacing $\uu^\BoldTau$ with $\UU^\BoldTau$ in all terms in the expression for $\W[\uu, p]$ having a derivative on $\uu^\BoldTau$. This gives
\begin{align}\label{e:H}
	\begin{split}
	\H^\BoldTau
		&:= \frac{1}{U^\n}
			\brac{
				-\prt_t \UU^\BoldTau - \grad_\Gamma \pr{p_r + \frac{1}{2} \abs{\UU}^2}
			 	 	+ \f^\BoldTau }^\perp
				+ \frac{1}{U^\n}
					\curl_\Gamma \UU^\BoldTau \uu^\BoldTau, \\
	H^\n
		&:= \curl_\Gamma \UU^\BoldTau,
	\end{split}
\end{align}
and we see that property $(P1)$ of $\H$ holds. We will show property (P2) in \cref{P:WeakVorticityIsVelocity}.

\begin{remark}\label{R:HHigherRegularity}
    Because we assumed that $\UU$ has higher regularity than $\uu$, the function $\H$ has one more derivative than $\W[\uu, p]$ in \cref{e:W}. This higher regularity will be needed to obtain estimates on $p_r$ in \cref{S:Pressure};
    it is needed as well to solve the linearized problem in \cref{T:LinearExistence}, though only for $N = 0$.
\end{remark}

The next proposition shows that our choice of $\H$ does, in fact, satisfy the constraint in \cref{e:RangeOfCurlCond} that is necessary to ensure $\curl \uu = \bomega$.

\begin{prop}\label{P:H0CondSatisfied}
	Assume that the data has regularity $0$ as in \cref{D:NReg}.
	For $\uu \in \uInputSpaceZero$,
	the condition in \cref{e:RangeOfCurlCond} is satisfied on $(0, T] \times \Gamma_+$.
\end{prop}
\begin{proof}
	From \cref{e:H} and using that $\curl_\Gamma \UU^\BoldTau = H^\n$ we have
	\begin{align*}
		U^\n &\H^\BoldTau - H^\n \uu^\BoldTau
			= \brac{
				-\prt_t \UU^\BoldTau - \grad_\Gamma \pr{p_r + \frac{1}{2} \abs{\UU}^2}
			 	 	+ \f^\BoldTau}^\perp.
	\end{align*}
	By \cref{e:curlRestricted}, $\dv_\Gamma \vv = - \dv_\Gamma ((\vv^\perp)^\perp)
	= \curl_\Gamma \vv^\perp$ for any tangent vector $\vv$. Hence,
	\begin{align*}
		\prt_t H^\n
				+ &\dv_\Gamma[H^\n \uu^\BoldTau  - U^\n \H^\BoldTau]
				- \g \cdot \n
				= \prt_t H^\n
					+ \curl_\Gamma[(H^\n \uu^\BoldTau - U^\n \H^\BoldTau)^\perp]
					- \g \cdot \n \\
			&= \prt_t \curl_\Gamma \UU^\BoldTau - \prt_t \curl_\Gamma \UU^\BoldTau
					+ \g \cdot \n - \g \cdot \n
			= 0,
	\end{align*}
	where $\curl_\Gamma \f^\BoldTau = (\curl \f) \cdot \n = \g \cdot \n$ by
	\cref{e:curlRestricted}. This gives \cref{e:RangeOfCurlCond}.
\end{proof}

\subsection*{The operator $A$} Our proof of \cref{T:WellPosedness} in \cref{S:Proof} will involve establishing the existence of a fixed point of an operator, which we denote $A$. We first define its domain:

\begin{definition}[Domain of $A$]\label{D:DomNA}
    Assume that the data has regularity $N$ as in \cref{D:NReg} and fix $\uu_0 \in \uInitSpaceN$ satisfying $\cond_N$. Define
    \begin{align}\label{e:DomNA}
    \begin{split}
    	\Dom_N(A)
    		:= \set{&\uu \in \uSolSpaceN \colon
    			\uu(0) = \uu_0, \,
    			\prt_t^n \uu|_{t = 0} = \wprt^n \uu_0
    			\text{ on } \ol{\Omega}, \,
    		0 \le n \le N},
    	\end{split}
    \end{align}
    where $\wprt^n \uu_0$, which appears in $\cond_N$ of \cref{e:condN}, will be defined
    in detail in \cref{D:wprt}.
\end{definition}

\begin{remark}\label{R:DomACond}
    The condition in $\Dom_N(A)$ that $\prt_t^n \uu(0) = \wprt^n \uu_0$ on $\ol{\Omega}$ for all $1 \le n \le N$
arises from the need to show that the nonlinear compatibility conditions imply the linear compatibility conditions for velocities in  $\Dom_N(A)$---as we will show in  \cref{P:condImplieslincond}.
\end{remark}

\begin{definition}[Operator $A$]\label{D:A}
Assume that the data has regularity $N$ as in \cref{D:NReg} and fix $\uu_0 \in \uInitSpaceN$ satisfying $\cond_N$. Let $\uu \in \Dom_N(A)$ and define $\H$ as in \cref{e:H}. Let $\ol{\bomega} \in C^{N, \al}(Q)$ be the unique solution to \OurLinear with $\ol{\bomega}_0 = \bomega_0 = \curl \uu_0$ given by \cref{T:LinearExistence} (see \cref{R:A}), with $\vv, \pi$ the corresponding velocity field $\vv \in \uSolSpaceN$ and pressure $\pi$ with $\curl \vv = \ol{\bomega}$ satisfying \cref{e:LinearVelEq}. Define
\begin{align}\label{e:A}
	A \uu := \vv,
\end{align}
and define also
\begin{align}\label{e:Lambda}
	\Lambda \uu
		:= \curl A \uu
		= \ol{\bomega}.
\end{align}
\end{definition}

\begin{remark}\label{R:A}
In \cref{D:A}, we use that $\lincond_N$ is satisfied for any $\uu \in \Dom_N(A)$ as in \cref{R:DomACond}, and that \cref{e:RangeOfCurlCond} is satisfied by \cref{P:H0CondSatisfied}. This allows us to apply \cref{T:LinearExistence} to obtain $\ol{\bomega}$, $\vv$, and $\pi$ in the given spaces.
\end{remark}

%
%
\section{Compatibility conditions}\label{S:CompConds}

\noindent 
If $(\uu, p)$ is a classical solution to \OurWeakBCs and $\bomega := \curl \uu$, then, of course,
\begin{align}\label{e:RealSolAtInitTime}
	\begin{split}
	\prt_t \bomega(0)
		&= \bomega_0 \cdot \grad \uu_0 - \uu_0 \cdot \grad \bomega_0
			+ \g(0), \\
	\prt_t \uu(0)
		&= - \uu_0 \cdot \grad \uu_0 - \grad p(0) + \f(0),
	\end{split}
\end{align}
where $\g := \curl \f$. This simple observation is behind both $\cond_N$ of \cref{e:condN} and $\lincond_N$ of \cref{e:lincondN}, which are based upon applying $\prt_t$, $n$ times, each time replacing $\prt_t \uu$ or $\prt_t \bomega$ at time zero with the relation in \cref{e:RealSolAtInitTime}, thereby replacing all time derivatives with spatial derivatives. This produces expressions, $\wprt^n \uu_0$ and $\wprt^n \bomega_0$, which are equal to $\prt_t^n \uu|_{t = 0}$ and $\prt^n \bomega|_{t = 0}$, respectively, for any actual solution to the Euler equations having sufficient regularity.

We make this process of constructing $\wprt^n \uu_0$ and $\wprt^n \bomega_0$ precise in \cref{D:wprt}. We stress that in this definition, we assume of $\uu$ only that it lies in $\uInputSpaceN$ with $\uu(0) = \uu_0$.

\begin{definition}\label{D:wprt}
Assume that the data has regularity $N \ge 0$ as in \cref{D:NReg} and let $\uu \in \uInputSpaceN$ with $\uu(0) = \uu_0$.
In accord with \cref{e:RealSolAtInitTime}, we define
\begin{align*}
	&\wprt \uu_0
		:= - \uu_0 \cdot \grad \uu_0 - \grad p^0 + \f(0), \quad
	\wprt \bomega_0
		:= - \uu_0 \cdot \grad \bomega_0 + \bomega_0 \cdot \grad \uu_0 + \g(0),
\end{align*}
where $p^0$ satisfies \cref{e:p0}.
We then define (recall that $\g := \curl \f$)
\begin{align}\label{e:wprt2}
	\begin{split}
	\wprt^2 \uu_0
		&:= - \wprt (\uu_0 \cdot \grad \uu_0) - \grad \wprt p^0 + \prt_t \f|_{t = 0}, \\
	\wprt^2 \bomega_0
		&:= - \wprt \uu_0 \cdot \grad \bomega_0 - \uu_0 \cdot \grad \wprt \bomega_0
			+ \wprt \bomega_0 \cdot \grad \uu_0 + \bomega_0 \cdot \grad \wprt \uu_0 + \prt_t \g|_{t = 0},
	\end{split}
\end{align}
where
\begin{align}\label{e:wprtu0gradu0}
	\wprt (\uu_0 \cdot \grad \uu_0)
		:= \wprt \uu_0 \cdot \grad \uu_0 + \uu_0 \cdot \grad \wprt \uu_0,
\end{align}
and define $\wprt p^0$ to be the unique mean-zero solution to (see \cref{R:wprtqSolvable}, below)
\begin{align*}
	\begin{cases}
		\Delta \wprt p^0
			= - \dv \wprt (\uu_0 \cdot \grad \uu_0)
			&\text{in } \Omega, \\
		\grad \wprt p^0 \cdot \n = -\prt_t^2 U^\n|_{t = 0} - \wprt (\uu_0 \cdot \grad \uu_0) \cdot \n
			&\text{on } \Gamma.
	\end{cases}
\end{align*}
We note that
\begin{align*}
	\wprt^2 \uu_0
		= - (- \uu_0 \cdot &\grad \uu_0 - \grad p^0 + \f(0)) \cdot \grad \uu_0
				- \uu_0 \cdot \grad (- \uu_0 \cdot \grad \uu_0 - \grad p^0 + \f(0)) \\
			&- \grad \wprt p^0 + \prt_t \f|_{t = 0}.
\end{align*}

For $\wprt^n$, we repeat this process inductively, up to order $N + 1$ for $\wprt \uu$ and order $N$ for $\wprt \bomega_0$.
\end{definition}

\begin{remark}\label{R:wprtqSolvable}
In the inductive extension of $\wprt^n p^0$ in \cref{D:wprt}, $\wprt^n p^0$ is the unique mean-zero solution to
\begin{align}\label{e:wprtqSystem}
	\begin{cases}
		\Delta \wprt^n p^0
			= - \dv \wprt^n (\uu_0 \cdot \grad \uu_0)
			&\text{in } \Omega, \\
		\grad \wprt^n p^0 \cdot \n = -\prt_t^{n + 1} U^\n|_{t = 0}
			- \wprt^n (\uu_0 \cdot \grad \uu_0) \cdot \n
			&\text{on }\Gamma.
	\end{cases}
\end{align}
Then
\begin{align*}
	\int_\Gamma \brac{\prt_t^n U^\n(0) + \wprt^n (\uu_0 \cdot \grad \uu_0) \cdot \n}
		= \int_\Omega \dv \wprt^n (\uu_0 \cdot \grad \uu_0),
\end{align*}
since $\dv \UU = 0$. Hence, \cref{e:wprtqSystem} is solvable.
Also, from the manner in which $\wprt^n p^0$ was defined,
$
	\wprt^n \uu_0 \cdot \n = \prt_t^n U^\n(0)
		\text{ on } \Gamma.
$
\end{remark}

In \cref{D:wprt}, $\wprt^n$ does not represent a derivative. Rather, it is a shorthand notation to properly account for the combinatorial nature of $\lincond_N$ and $\cond_N$.

\cref{P:curlwprt} shows that, in effect, $\wprt^n \curl \uu = \curl \wprt^n \uu$.
		 
\begin{prop}\label{P:curlwprt}
	Let $\uu$ and $\wprt^n$ be as in \cref{D:wprt}.
	Then
	\begin{align}\label{e:wprtcurl}
		\begin{array}{ll}
		\dv \wprt^n \uu_0 = 0
			&\text{for all } 0 \le n \le N + 1, \\
		\wprt^n \bomega_0
			= \curl \wprt^n \uu_0
			&\text{for all } 0 \le n \le N.
		\end{array}
	\end{align}
\end{prop}
\begin{proof}
	We constructed the pressure $p^0$ in \cref{D:wprt}
	so that $\dv \wprt^n \uu_0 = 0$. Then,
	for $n = 1$, \cref{e:wprtcurl}$_2$ follows from the identity,
	$\curl (\uu_0 \cdot \grad \uu_0 + \grad p^0)
		= \uu_0 \cdot \grad \bomega_0 - \bomega_0 \cdot \grad \uu_0$.
	
	For $n = 2$, we will use two vector identities that hold for vector fields
	$\aaa, \bbb \in C^2(\ol{\Omega})$ and follow from direct calculation:
	\begin{align}\label{e:TwoVectorIds}
		\begin{split}
		&\aaa \cdot (\grad \bbb)^T + \bbb \cdot (\grad \aaa)^T
			= \grad (\aaa \cdot \bbb), \\
		&\curl(\aaa \cdot \grad \bbb + \bbb \cdot (\grad \aaa)^T)
			= \aaa \cdot \grad \curl \bbb - \curl \aaa \cdot \grad \bbb
				\quad (\text{if } \dv \aaa = 0).
		\end{split}
	\end{align}
	Then,
	\begin{align*}
			\wprt^2 \bomega_0
		&= - \wprt \uu_0 \cdot \grad \bomega_0 - \uu_0 \cdot \grad \curl (\wprt \uu_0)
			+ \curl (\wprt \uu_0) \cdot \grad \uu_0 + \bomega_0 \cdot \grad \wprt \uu_0
				+ \g(0) \\
		&= \curl (\wprt \uu_0) \cdot \grad \uu_0 - \uu_0 \cdot \grad \curl (\wprt \uu_0)
			+ \bomega_0 \cdot \grad \wprt \uu_0 - \wprt \uu_0 \cdot \grad \bomega_0 + \g(0)
			\\
		&= - \curl (\wprt \uu_0 \cdot \grad \uu_0 + \wprt \uu_0 \cdot (\grad \uu_0)^T)
				- \curl (\uu_0 \cdot \grad \wprt \uu_0 + \uu_0 \cdot (\grad \wprt \uu_0)^T)
					+ \g(0) \\
		&= \curl (- \wprt \uu_0 \cdot \grad \uu_0 - \uu_0 \cdot \grad \wprt \uu_0
				+ \prt_t \f(0))
			- \curl(\wprt \uu_0 \cdot (\grad \uu_0)^T + \uu_0 \cdot (\grad \wprt \uu_0)^T) \\
		&= \curl (- \wprt \uu_0 \cdot \grad \uu_0 - \uu_0 \cdot \grad \wprt \uu_0
			- \grad \wprt p^0 + \prt_t \f(0))
		= \curl \wprt^2 \uu_0.
	\end{align*}
	The first equality is \cref{e:wprt2}$_2$, the second equality is a rearrangement of it.
	The third equality used \cref{e:TwoVectorIds}$_2$, the fourth equality rearranges
	terms. The fifth equality uses \cref{e:TwoVectorIds}$_2$ and $\curl \grad = 0$.
	The final equality uses \cref{e:TwoVectorIds}$_1$.
	
	Equality in \cref{e:wprtcurl}$_2$ follows inductively for higher values of $n$.
\end{proof}

For the linearized problem in \OurLinear, $\uu$ is given, so $\lincond_N$ is a condition on the data, $\H$, along with the forcing $\g$. For the nonlinear problem in \Our, however, $\H$ must be generated on the inflow boundary from the solution itself. This means that our nonlinear compatibility condition $\cond_N$ must be such that the function $\H$ given by \cref{e:H} satisfies the linear compatibility condition $\lincond_N$. We begin with the $N = 0$ case.

Using \cref{L:CrossProduct} along with $\curl_\Gamma \UU^\BoldTau = H^\n$, on all of $[0, T] \times \Gamma_+$, we have
\begin{align*}
	- [[\uu \times \H]^\BoldTau]^\perp
		&= U^\n \H^\BoldTau - H^\n \uu^\BoldTau
		= \brac{-\prt_t \UU^\BoldTau - \grad_\Gamma \pr{p_r + \frac{1}{2} \abs{\UU}^2}
			 	 	+ \f^\BoldTau}^\perp,
\end{align*}
so,
\begin{align}
	\begin{split}
	[\uu \times \H]^\BoldTau
		&= \prt_t \UU^\BoldTau + \grad_\Gamma \pr{p_r + \frac{1}{2} \abs{\UU}^2}
			 	 	- \f^\BoldTau \\
		&= \prt_t \UU^\BoldTau + \grad_\Gamma \pr{p_r + \frac{1}{2} \abs{\uu}^2}
			 	 	- \f^\BoldTau
					+ \frac{1}{2} \grad_\Gamma \pr{\abs{\UU}^2 - \abs{\uu}^2}.
	\end{split}
\end{align}

Then using the vector identity in \cref{e:NonlinerIdentity},
\begin{align}\label{e:FromNonlinearIdentity}
	\grad_\Gamma &\pr{p_r + \frac{1}{2} \abs{\uu}^2} - \f^\BoldTau
		= [\uu \cdot \grad \uu + \grad p_r - \f]^\BoldTau
			+ [\uu \times \bomega]^\BoldTau.
\end{align}
At $t = 0$, $p_r = p^0$ so $[\uu \cdot \grad \uu + \grad p_r - \f]^\BoldTau = -\wprt \uu_0^\BoldTau$, and on $\Gamma_+$, 
\begin{align}\label{e:FromHwithwprt}
	[\uu \times \H]^\BoldTau_{t = 0}
		&= \prt_t \UU^\BoldTau|_{t = 0} -\wprt \uu_0^\BoldTau 
					+ \frac{1}{2} \grad_\Gamma \pr{\abs{\UU(0)}^2 - \abs{\uu_0}^2}
					+ [\uu_0 \times \bomega_0]^\BoldTau.
\end{align}

\begin{prop}\label{P:condImplieslincond0}
	Assume the data has regularity $0$, $\uu \in \uInputSpaceZero$, and
	$\cond_0$ in \cref{e:condN} holds. Then
	$\lincond_0$ in \cref{e:lincondN} holds.	
\end{prop}
\begin{proof}
\textbf{All the calculations in this proof apply at time zero on $\Gamma_+$.}

We have $\prt_t \UU^\BoldTau -\wprt \uu_0^\BoldTau = 0$ by $\cond_0$. Since also $\uu(0) = \UU(0)$ on $\Gamma_+$, we know that $\grad_\Gamma \abs{\UU}^2 = \grad_\Gamma \abs{\uu}^2$, and \cref{e:FromHwithwprt} reduces to $[\UU \times \H]^\BoldTau = [\UU \times \bomega]^\BoldTau$, or,
\begin{align*}
	[\UU \times (\H - \bomega)]^\BoldTau = 0.
\end{align*}
Also from \cref{e:H}$_2$, $H^\n = \curl_\Gamma \UU^\BoldTau = \curl_\Gamma \uu^\BoldTau = \omega^\n$. Then, since $H^\n = \omega^\n$ and only $(\H - \bomega)^\BoldTau$ contributes to $\n \times (\H - \bomega)$, we can apply the vector identity, $A \times (B \times C )=(A \cdot C )B -(A \cdot B )C$ to give
\begin{align*}
	0
		&= \n \times [\UU \times (\H - \bomega)]^\BoldTau
		= \n \times [\UU \times (\H - \bomega)] \\
		&= [\n \cdot (\H - \bomega)] \UU
			- [\n \cdot \UU] (\H - \bomega)
		= - U^\n (\H - \bomega).
\end{align*}
Since $U^\n$ never vanishes on $\Gamma_+$, we conclude that $\H = \bomega$ on $\set{0} \times \Gamma_+$, meaning that $\lincond_0$ is satisfied.
\end{proof}

If $\uu \in \uSolSpaceN$ for $N \ge 1$ , however, $\cond_N \centernot \implies \lincond_N$ unless we restrict $\uu$ to the subspace $\Dom_N(A)$ of $\uSolSpaceN$. To show this, we will find it convenient to extend our definition of the $\wprt^n$ ``operator'' to apply to all of $[0, T] \times \ol{\Omega}$
by replacing $p^0$ with $p_r$ in \cref{D:wprt} and not restricting the calculations to $t = 0$. Since $p_r|_{t = 0} = p^0$, the definitions of $\wprt \uu_0$ and $\wprt \bomega_0$ are unchanged in the sense that
\begin{align*}
	\wprt \uu_0 = \wprt \uu|_{t = 0}, \quad \wprt \bomega_0 = \wprt \bomega|_{t = 0},
\end{align*}
where we have used the same symbol $\wprt$ for both versions of the ``operator.''
We then define
\begin{align}\label{e:wprt2ut0}
	\wprt^2 \uu
		= - (- \uu\cdot &\grad \uu - \grad p_r + \f) \cdot \grad \uu
				- \uu \cdot \grad (- \uu \cdot \grad \uu - \grad p_r + \f)
			- \grad \wprt p_r + \prt_t \f
\end{align}
on $[0, T] \times \ol{\Omega}$, where
\begin{align*}
	\begin{cases}
		\Delta \wprt p_r
			= - \dv \wprt (\uu \cdot \grad \uu)
			&\text{in } \ol{Q}, \\
		\grad \wprt p_r \cdot \n = -\prt_t U^\n - \wprt N[\uu]
			&\text{on } [0, T] \times \Gamma,
	\end{cases}
\end{align*}
with
\begin{align*}
	\begin{split}
	\wprt N[\uu]
		&:=
		\begin{cases}
			\wprt (\uu \cdot \grad \uu) \cdot \n
					&\text{on } [0, T] \times (\Gamma_- \cup \Gamma_0), \\
			\wprt (\uu \cdot \grad \uu) \cdot \n
				+ \dv_\Gamma (U^\n (\wprt \uu^\BoldTau - \prt_t \UU^\BoldTau))
					&\text{on } [0, T] \times \Gamma_+.
		\end{cases}
	\end{split}
\end{align*}
Then $\wprt^n \uu$ is defined inductively for all $n > 2$, and $\wprt^n \bomega$ is defined similarly. An additional assumption on $u$ is required, however, to have the two definitions of $\wprt^n \uu$ at time zero to agree for $n \ge 2$, as we see in \cref{P:wprtInduct}.

\begin{prop}\label{P:wprtInduct}
	Assume that the data has regularity $N \ge 0$, $\cond_N$ holds,
	and $\prt_t \uu|_{t = 0} = \wprt \uu_0$ on $\ol{\Omega}$.
	Then
	\begin{itemize}
		\item[(1)]
			$\wprt^n p_r|_{t = 0} = \wprt^n p^0$
				on $\ol{\Omega}$ for all $0 \le n \le N$;\\[-4pt]
	
		\item[(2)]
			$\wprt^n \uu|_{t = 0} = \wprt^n \uu_0$
				on $\ol{\Omega}$ for all $0 \le n \le N + 1$;\\[-4pt]

		\item[(3)]			
			$\prt_t \wprt^n \uu|_{t = 0} = \wprt^{n + 1} \uu_0$ 
				on $\ol{\Omega}$ for all $0 \le n \le N$.
	\end{itemize}
\end{prop}
\begin{proof}
	By $\cond_N$, we have $\prt_t^n \UU^\BoldTau = \wprt^n \uu^\BoldTau$
	on $\Gamma_+$, so $\wprt^n N[\uu] = \wprt^n (\uu \cdot \grad \uu)$ and hence
	$\wprt^n p_r|_{t = 0} = \wprt^n p^0$.
	This gives (1) and then (2) follows directly from (1).
	For (3), observe that $\wprt^n \uu$ fully expanded
	(as in \cref{e:wprt2ut0} for $n = 2$) contains no time dervatives of $\uu$.
	Using the product rule to obtain $\prt_t \wprt^n \uu|_{t = 0}$ and using
	that $\prt_t \uu|_{t = 0} = \wprt \uu_0$ we obtain, using also (1), the same expression
	as $\wprt^{n + 1} \uu_0$.
\end{proof}

\begin{prop}\label{P:condImplieslincond}
	Assume that the data has regularity $N \ge 0$, $\cond_N$ holds,
	and $\uu \in \Dom_N(A)$.
	Then $\lincond_N$ in \cref{e:lincondN} holds.
\end{prop}
\begin{proof}
Let $N = 1$.
With our extended definition of $\wprt$, \cref{e:FromNonlinearIdentity} holds beyond time zero; that is,
\begin{align*}
	[\uu \times \H]^\BoldTau
		&= \prt_t \UU^\BoldTau -\wprt \uu^\BoldTau 
					+ \frac{1}{2} \grad_\Gamma \pr{\abs{\UU}^2 - \abs{\uu}^2}
					+ [\uu \times \bomega]^\BoldTau
			\text{ on }
\end{align*}
on all of $[0, T] \times \Gamma_+$ (this does not require any compatibility conditions).
Differentiating both sides in time gives
\begin{align}\label{e:prttH0}
	\begin{split}
	[\prt_t \uu \times \H]^\BoldTau + [\uu \times \prt_t \H]^\BoldTau
		= \prt_{tt} &\UU^\BoldTau - \prt_t \wprt \uu^\BoldTau 
					+ \frac{1}{2} \grad_\Gamma \prt_t \pr{\abs{\UU}^2 - \abs{\uu}^2} \\
		&
			+ [\prt_t \uu \times \bomega]^\BoldTau + [\uu \times \prt_t \bomega]^\BoldTau
	\end{split}
\end{align}
on $[0, T] \times \Gamma_+$. We know from \cref{P:condImplieslincond0} that if $\cond_0$ holds then $\H = \bomega$ on $\set{0} \times \Gamma_+$, so two terms above cancel, leaving, at time zero,
\begin{align}\label{e:uprtHprtomega}
	[\uu \times \prt_t \H]^\BoldTau
		= \pbrac[\bigg]{\prt_{tt} &\UU^\BoldTau - \prt_t \wprt \uu_0^\BoldTau  
					+ \frac{1}{2} \grad_\Gamma \prt_t \pr{\abs{\UU}^2 - \abs{\uu}^2}}
			+ [\uu \times \prt_t \bomega]^\BoldTau
			\text{ on } \set{0} \times \Gamma_+.
\end{align}

	From \cref{P:curlwprt},
	$\wprt \bomega_0 = \curl \wprt \uu_0 = \curl \prt_t \uu(0)
	= \prt_t \curl \uu(0) = \prt_t \bomega(0)$,
	and from \cref{P:wprtInduct} we know that
	$\prt_t \wprt \uu|_{t = 0} = \wprt^2 \uu_0$.
	Also,
	\begin{align}\label{e:prttU2u2}
		\begin{split}
		\prt_t \pr{\abs{\UU}^2 - \abs{\uu}^2}|_{t = 0}
			&= 2 (\UU \cdot \prt_t \UU - \uu \cdot \prt_t \uu)|_{t = 0}
			= 2 (\UU(0) \cdot \prt_t \UU|_{t = 0} - \UU(0) \cdot \wprt \uu_0) \\
			&= 2 (\UU \cdot \prt_t \UU - \UU \cdot \prt_t \UU)|_{t = 0}
			= 0,
		\end{split}
	\end{align}
	where we used that $\uu \in \Dom_N(A)$ in the second equality
	and $\cond_1$ with $\wprt \uu_0 \cdot \n = \prt_t U^\n(0)$
	on $\Gamma$ as in \cref{R:wprtqSolvable} in the third equality.
 
	Thus, the term in the brackets in \cref{e:uprtHprtomega} vanishes
	because of $\cond_1$,
	and we are left with
	\begin{align*}
		[\uu_0 \times \prt_t \H|_{t = 0}]^\BoldTau
			= [\uu_0 \times \wprt \bomega_0]^\BoldTau
			\text{ on } \Gamma_+.
	\end{align*}
	Also, $\prt_t H^\n|_{t = 0} = \curl_\Gamma \prt_t \UU^\BoldTau|_{t = 0}
	= \curl_\Gamma \wprt \uu_0^\BoldTau = \wprt \omega_0^\n$,
	so arguing as in the proof of \cref{P:condImplieslincond0}, we see that
	$\prt_t \H|_{t = 0} = \wprt \bomega_0$, which is $\lincond_1$.
	
	The result for $N \ge 2$ follows inductively, where we note that,
	as in \cref{e:prttU2u2}, showing that
	$\prt_t^N \pr{\abs{\UU}^2 - \abs{\uu}^2}|_{t = 0} = 0$
	uses that $\prt_t^n \uu|_{t = 0} = \wprt^n \uu_0$ on $\Gamma_+$ for all
	$0 \le n \le N$, since $\uu \in \Dom_N(A)$.
\end{proof}

\begin{remark}
	In  the proof of \cref{P:condImplieslincond}, we only required of $\uu$
	that $\uu(0) = \uu_0$ satisfy $\cond_N$,
	$\prt_t \uu|_{t = 0} = \wprt \uu_0$ on $\ol{\Omega}$ (because we applied
	\cref{P:wprtInduct}),
	and $\prt_t^n \uu_{t = 0} = \wprt^n \uu_0$ on $\Gamma_+$.
	The full conditions on $\Dom_N(A)$ will be required
	shortly, however, in \cref{P:DomatoDomA}.
\end{remark}

%
%
\section{Proof of well-posedness with inflow, outflow}\label{S:Proof}

\noindent In this section, we present the three key propositions on which the proof of \cref{T:WellPosedness} relies, then give the proof  of \cref{T:WellPosedness} itself.

\begin{prop}\label{P:DomatoDomA}
	$A$ maps $\Dom_N(A)$ to itself.
\end{prop}
\begin{proof}
	Let $\uu \in \Dom_N(A)$ and let $\vv = A \uu$.
	\cref{T:LinearExistence} shows that $\vv \in  \uSolSpaceN$ and $\vv(0) = \uu_0$,
	so it remains only to show that $\prt_t^n \vv|_{t = 0} = \wprt^n \uu_0$
	for $1 \le n \le N$.
	
	Suppose $N = 1$. Then since $\vv(0) = \uu(0)$, \cref{e:LinearVelEq} gives
	\begin{align*}
		\prt_t \vv|_{t = 0}
			&= -\uu_0 \cdot \grad \uu_0 + \uu_0 \cdot (\grad \uu_0)^T - \grad \pi(0) + \f(0).
	\end{align*}
	But $\uu_0 \cdot (\grad \uu_0)^T  = (1/2) \grad \abs{\uu_0}^2$, so we have
	\begin{align*}
		\prt_t \vv|_{t = 0}
			&= -\uu_0 \cdot \grad \uu_0 - \grad r + \f(0)
	\end{align*}
	for some ``pressure'' $r$. But $r$ is recovered in the same manner as $p$,
	which is the same as $p_r$ at time zero. We see, then, that
	$\prt_t \vv|_{t = 0} = \wprt \uu_0$.
	
	For $N = 2$, a time derivative of \cref{e:LinearVelEq} yields,
	\begin{align*}
		\prt^2_t \vv
			+ \prt_t \uu \cdot \grad \vv
			+ \uu \cdot \grad \prt_t \vv
			- [\prt_t \uu \cdot (\grad \vv)^T
				+ \uu \cdot (\grad \prt_t \vv)^T]
			= -\grad \prt_t \pi + \prt_t \f.
	\end{align*}
	From the $N = 1$ result, $\prt_t \uu|_{t = 0} = \prt_t \vv|_{t = 0} = \wprt \uu_0$ 
	and we have
	\begin{align*}
		\prt^2_t \vv|_{t = 0}
			+ \wprt \uu_0 \cdot \grad \uu_0
			+ \uu_0 \cdot \grad \wprt \uu_0
			- [\wprt \uu_0 \cdot (\grad \uu_0)^T
				+ \uu_0 \cdot (\grad \wprt \uu_0)^T]
			= -\grad \prt_t \pi|_{t = 0}  + \prt_t \f|_{t = 0}.
	\end{align*}
	By \cref{e:TwoVectorIds}$_1$, $\wprt \uu_0 \cdot (\grad \uu_0)^T
	+ \uu_0 \cdot (\grad \wprt \uu_0)^T$ is a gradient, so we see, also using
	\cref{e:wprtu0gradu0}, that
	\begin{align*}
		\prt^2_t \vv|_{t = 0}
			+ \wprt (\uu_0 \cdot \grad \uu_0)
			= \grad q + \prt_t \f|_{t = 0}
	\end{align*}
	for some $q$.
	But from \cref{e:wprt2}$_1$,
	\begin{align*}
		\wprt^2 \uu_0
			+ \wprt (\uu_0 \cdot \grad \uu_0)
			= \grad \wprt p^0 + \prt_t \f|_{t = 0}.
	\end{align*}
	Hence, $\ww := \prt^2_t \vv|_{t = 0} - \wprt^2 \uu_0$ is a gradient
	in the space $H$, since $\dv \ww = 0$ with $\ww \cdot \n = 0$ on $\Gamma$; hence,
	$\ww = 0$, giving $\prt^2_t \vv|_{t = 0} = \wprt^2 \uu_0$.
	
	The result for $N > 2$ follows inductively. We note that because it involves
	differentiating in time \cref{e:LinearVelEq} $N - 1$ times, it requires that
	$\prt_t^n \uu|_{t = 0} = \wprt^n \uu_0$ on $\Omega$ for all  $1 \le n \le N$,
	so the full condition on $\uu$ in $\Dom_N(A)$ is required.	
\end{proof}

We will also show in \cref{L:DomNA} that $\Dom_N(A)$ is a nonempty, convex subset of $\uInputSpaceN$.

We will apply Schauder's fixed point theorem to obtain the existence of a fixed point of $A$, but this requires that $A$ be continuous. Results in \cite{Paper1} would give that $A$ is bounded as a map from $\Dom_N(A)$ to $\Dom_N(A)$ in the $\uInputSpaceN$ norm, as long as we can obtain sufficient control of the pressure so as to control $\H$. But $A$, which is nonlinear, need not be continuous from $\uInputSpaceN$ to itself. To ensure continuity, we need to work with a weaker topology, which we introduce next.

\begin{definition}
    Fixing
    $
        \beta \in (0, \al]
    $,
    we define the (affine) space $\WeakS$ to be all vector fields in $\uSolSpaceN$ endowed with the norm,
    \begin{align*}
        \norm{\uu}_\WeakS
            &= \norm{\uu}_{C^{N, \beta}(Q)} + \norm{\curl \uu}_{C^{N, \beta}(Q)}.
    \end{align*}
\end{definition}

We note that $\WeakS$ is a locally convex topological affine space. Because $C^{N, \al}(Q)$ is compactly embedded in $C^{N, \beta}(Q)$ for $\beta < \al$, we see that $\WeakSalpha$ is compactly embedded in $\WeakS$ for $\beta < \al$. Like $\uSolSpaceN$, which is also an affine space, we will often apply the $\WeakS$ norm to the difference of two elements in $\WeakS$, even though that difference does not lie in the space. In particular, we do this in \cref{P:AContinuous}. Finally, observe that $\norm{\uu}_\uSolSpaceN = \norm{\uu}_{\WeakSalpha} + \norm{ \prt_t^{N + 1} \uu}_{L^\iny([0, T]; C^\al(Q))}$.

In outline, our proof of \cref{T:WellPosedness} is as follows: We show that $A$ maps a nonempty convex set $K \subseteq \Dom_N(A)$ into itself (an $A$-invariant set), that $K$ is compact in the $\WeakS$ norm, and that $A$ is continuous on $K$ in the $\WeakS$ norm. Applying Schauder's fixed point theorem gives the existence of a fixed point. We show a posteriori that the full inflow, outflow boundary conditions in \OurBCsOnly are satisfied, and, finally, prove uniqueness.

These steps are detailed in \crefrange{P:InvariantSet}{P:WeakVorticityIsVelocity}, followed by the proof itself. To streamline the presentation, we defer the proofs of these technical propositions to later sections.

\begin{prop}\label{P:InvariantSet}
    For all $M$ larger than a value that depends only upon the data, there exists $T > 0$ for which the set
    \begin{align}\label{e:KBall}
        \Cal{K}
        	= \Cal{K}_{M, T}
            := \set{\uu \in \Dom_N(A) \colon
                \norm{\uu}_\uSolSpaceN \le M}
    \end{align}
    is invariant under $A$. That is, $\uu \in \Dom_N(A)$ with $\norm{\uu}_\uSolSpaceN \le M$ implies that $A \uu \in \Dom_N(A)$ with $\norm{A \uu}_\uSolSpaceN  \le M$.
\end{prop}
\begin{proof}
	Given in \cref{S:InvariantSet}. We note here only that $\Cal{K}$ depends on
	$T$ because each $\uu \in \Dom_N(A)$ is defined on $Q = Q_T$.
	When $M$ and $T$ are fixed, we will generally refer to the set simply as $\Cal{K}$.
\end{proof}

\begin{prop}\label{P:AContinuous}
	With $\Cal{K}$ as in \cref{e:KBall},
	for any $\beta \in (0, \al)$,
	$A \colon \Cal{K} \to \Cal{K}$ is continuous in the $\WeakS$ norm.
\end{prop}
\begin{proof}
	Given in \cref{S:AContinuous}.
\end{proof}

\begin{prop}\label{P:WeakVorticityIsVelocity}
	Assume that $(\uu, \grad p_r) \in \uSolSpaceZero \times C^\al(Q)$ and $(u, p_r)$
	solves \OurWeakBCs (with $p_r$ in place of $p$) and that
        $\curl \uu = \H$ on $[0, T] \times \Gamma_+$,
	with $\H$ given in \cref{e:H}.
	Then \OurInflow also holds.
\end{prop}
\begin{proof}
	Given in \cref{S:FullBCs}.
\end{proof}

\subsection*{Proof of well-posedness}
\cref{T:WellPosedness} we now see is a consequence of \crefrange{P:InvariantSet}{P:WeakVorticityIsVelocity}:

\begin{proof}[\textbf{Proof of \cref{T:WellPosedness}}]
\phantom{x}
Let $M > 0$ depending on the initial data, $T > 0$, and $\Cal{K} = \Cal{K}_{M, T}$ be given by \cref{P:InvariantSet}.

Choose any $\beta \in (0, \al)$. Because $C^{N, \al}$ is compactly embedded in $C^{N, \beta}$, and also using \cref{L:DomNA}, below, we see that $\Cal{K}$ is a nonempty convex compact subset of $\WeakS$, and $A \colon \Cal{K} \to \Cal{K}$ by \cref{P:InvariantSet}. By \cref{P:AContinuous}, $A$ is continuous as a map from $\Cal{K}$ to $\Cal{K}$ in the $\WeakS$ norm, and so has a fixed point $\uu$ by Schauder's Fixed Point Theorem. It follows that $A \uu = \uu$ with $\uu \in \WeakS$. Since $\uu \in \Dom_N(A)$ it follows that, in fact, $\uu \in \uSolSpaceN$.

Since $\vv := A \uu = \uu$, \cref{T:LinearExistence} implies that
$
	\prt_t \uu + \uu \cdot \grad \uu + \grad p
		= \f
$
for some pressure $p$. Hence, $(\uu, p)$ is a solution to \OurWeakBCs.
But since $\uu = A \uu$, we have $\bomega := \curl \uu = \H$ on $[0, T] \times \Gamma_+$. Thus, \cref{P:WeakVorticityIsVelocity} gives that \OurInflow holds, so $(\uu, p)$ is a solution to \Our.

To prove uniqueness, let $(\uu_1, p_1)$, $(\uu_2, p_2)$ be two solutions to \Our with the same initial velocity in $C^{1, \al}$ (so we prove uniqueness for $N = 0$ and it then follows for all $N \ge 0$). Letting $\ww = \uu_1 - \uu_2$, subtracting \OurOne for $(\uu_2, p_2)$ from \OurOne for $(\uu_1, p_1)$, 
\begin{align*}
	\prt_t \ww + \uu_1 \cdot \grad \ww + \ww \cdot \grad \uu_2 + \grad (p_1 - p_2) = 0.
\end{align*}
Multiplying by $\ww$ and integrating over $\Omega$, we obtain
\begin{align*}
	\frac{1}{2} \diff{}{t} \norm{\ww}^2
		&= - \int_\Omega (\ww \cdot \grad \uu_2) \cdot \ww
			- \frac{1}{2} \int_\Omega \uu_1 \cdot \grad \abs{\ww}^2
		\le \norm{\grad \uu_2}_{L^\iny(Q)} \norm{\ww}^2
			- \frac{1}{2} \int_\Omega \uu_1 \cdot \grad \abs{\ww}^2.
\end{align*}
But,
\begin{align*}
	-\int_\Omega \uu_1 \cdot \grad \abs{\ww}^2
		&= - \int_\Gamma (\uu_1 \cdot \n) \abs{\ww}^2
		= - \int_{\Gamma_-} (\uu_1 \cdot \n) \abs{\ww}^2
		\le 0,
\end{align*}
since $\ww = 0$ on $\Gamma_+$, $\uu_1 \cdot \n = 0$ on $\Gamma_0$, and
$\uu_1 \cdot \n > 0$ on $\Gamma_-$. Hence,
\begin{align*}
	\diff{}{t} \norm{\ww}^2
		&\le 2 \norm{\grad \uu_2}_{L^\iny(Q)} \norm{\ww}^2,
\end{align*}
and we conclude that $\ww = 0$ by \Gronwalls Lemma, giving the uniqueness in \cref{T:WellPosedness}.
\end{proof}

When $\Gamma_0 = \Gamma$---that is, when classical impermeable boundary conditions are imposed on the entire boundary---\cref{T:WellPosedness} gives well-posedness of the 3D Euler equations in $\uInputSpaceN$ for any $N \ge 0$. The proof simplifies in this case, as we discuss briefly in \cref{R:ImpermeableOnly}.

%
%
\Part{Part II: Preliminary Estimates}

\subsection*{Organization of Part II}
We introduce in \cref{S:Conventions} some conventions that we will use throughout the remainder of this paper to streamline the presentation. In \cref{S:BSLaw} we develop some properties related to the function space $\uSolSpaceN$, and in \cref{S:FlowMap} derive some properties of the flow map. We describe the generation of vorticity on the boundary in \cref{S:Boundary}, and obtain critical estimates on the pressure in \cref{S:Pressure}.

%
%
\section{Some conventions}\label{S:Conventions}
\noindent 

\subsection*{Constants}
To simplify notation, we write $M$ as a universal but unspecified bound on $\norm{\uu}_\uSolSpaceN$. Thus, we assume that
\begin{align}\label{e:M}
	\norm{\uu}_\uSolSpaceN \le M
		\text{ for some } M \ge 1
\end{align}
in what follows.
(Having $M \ge 1$ simplifies the form of some estimates.)

\begin{definition}\label{D:cc}
We define the following three types of positive ``constant'':
\begin{align*}
	c_0 &= c_0(\norm{\uu_0}_\uInitSpaceN,
			U_{min}^{-1}, \norm{\UU}_\USpaceNInf, \norm{\curl \f}_{\gspaceNInf}), \\
	c_X &= c_X(c_0, M), \\
    \cN &= \cN(c_0, M, T),
\end{align*}
where $U_{min}$ is as in \cref{D:NReg}.
Each of $c_0$, $c_X$, and $\cN$ are continuous, increasing functions of each of their arguments, and each appearance of $c_0$, $c_X$, and $\cN$  may have different values, even within the same expression. Moreover, we require of $\cN$ that for all $M > 0$,
\begin{align}\label{e:cN0}
    \cN(c_0, M, 0) = c_0.
\end{align}
\end{definition}

The property in \cref{e:cN0} is critical in establishing the existence of an invariant set for the operator $A$ in the proof of \cref{P:InvariantSet} in \cref{S:InvariantSet}. All the necessary quantities will be bounded by a $\cN$ constant, and all those bounds ultimately derive from the key bound on $\norm{\uu}_{C^{N, \al}(Q)}$ in \cref{P:SConsequences} for $\uu \in \Dom_N(A)$, which we establish in the next section.

In the process of obtaining constants $c_0$, $c_X$, or $\cN$ it will be clear that they increase with their arguments. It is sometimes clearer to write the expression for a constant $\cN$ more explicitly; typical examples are $c_0 + c_X T^\al$ and $e^{M T}$.

Because we imposed the restrictions on $\UU$ and $\f$ in \cref{D:NReg} for all time, a $c_0$ constant has no dependence on $T$.

\begin{remark}\label{R:TPowers}
	Many of our estimates contain factors
	of the form $C_1 T^{e_1} + C_2 T^{e_2} + C_3 T^{e_3}$, $0 < e_1 < e_2 < e_3$,
	where $C_1$, $C_2$, and $C_3$ may depend upon the norms of the data or
	the solution, but have no explicit dependence on time. To simplify matters,
	we will assume that $T \le T_0$ for some fixed but arbitrarily large $T_0 > 0$.
	Then
	\begin{align*}
		&C_1 T^{e_1} + C_2 T^{e_2} + C_3 T^{e_3}
			\le C_1 T^{e_1} + C_2 T^{e_1} T_0^{e_2 - e_1} + C_3 T^{e_1} T_0^{e_3 - e_1}
			\le C' T^{e_1}, \\\
		&C' := (1 + T_0^{e_2 - e_1} + T_0^{e_3 - e_1}) \max \set{C_1, C_2, C_3}.
	\end{align*}	
	Hence, in the final forms of estimates, we will only keep the lowest exponents
	of $T$ and, similarly, of $\abs{t_1 - t_2}$ for $t_1, t_2 \in [0, T]$.
\end{remark}

%
%
\section{Function spaces and the Biot-Savart law}\label{S:BSLaw}

\noindent In this section we give some basic properties of \Holder spaces that we will use throughout most of this paper and summarize those properties of the Biot-Savart law that we will need. We use these properties to establish the key estimate on velocity fields in \cref{P:SConsequences}, and then to prove \cref{L:DomNA}, showing that $\Dom_N(A)$ is nonempty and convex.

\subsection*{\Holder spaces}

Let $k \ge 0$ be an integer and $U$ be an open subset $U$ of $\R^d$, $d \ge 1$. We define $C^k(U)$ to be the space of all $k$-times continuously differentiable functions with the norm
\begin{align*}
    \norm{f}_{C^k(U)}
        &:= \sum_{\abs{\gamma} \le k} \norm{D^\gamma f}_{L^\iny(U)}.
\end{align*}
Letting $r \in (0, 1)$ we define the \Holder space, $C^{k, r}(U)$, to be the space of all $f \in C^k(U)$ for which
\begin{align}\label{e:HolderNormDef}
	\begin{split}
    \norm{f}_{C^{k, r}(U)}
        &:= \norm{f}_{C^k(U)}
            + \sum_{\abs{\gamma} = k} \norm{D^\gamma f}_{\dot{C}^r(U)} < \iny, \\
    \norm{g}_{\dot{C}^r(U)}
        &:= \sup_{x \ne y \in U} \frac{\abs{g(x) - g(y)}}{\abs{x - y}^r}.
	\end{split}
\end{align}

Now consider a time-space domain of the form $Q = [0, T] \times \Omega$. For any $f$ continuous on $Q$ and $r \in (0, 1]$, define
\begin{align}\label{e:Crtx}
	\begin{split}
    \norm{f}_{\dot{C}^r_t(Q)}
        &:= \norm{f}_{\dot{C}^r([0, T]; L^\iny(\Omega))}
        = \sup_{\x \in \Omega} \norm{f(\cdot, \x)}_{\dot{C}^r([0, T])},
        \\
    \norm{f}_{\dot{C}^r_x(Q)}
        &:= \norm{f}_{L^\iny([0, T]; \dot{C}^r(\Omega))}
        = \sup_{t \in [0, T]} \norm{f(t, \cdot)}_{\dot{C}^r(\Omega)}.
	\end{split}
\end{align}

\begin{lemma}\label{L:HolderEquiv}
    For any integer $k \ge 0$ and $r \in (0, 1)$,
    \begin{align*}
        \norm{f}_{C^N(Q)} + \norm{f}_{\dot{C}^r_t(Q)} + \norm{f}_{\dot{C}^r_t(Q)}
    \end{align*}
    is equivalent to the $C^{N, r}(Q)$ norm.
\end{lemma}

\subsection*{Biot-Savart law}
We need a few facts from \cite{Paper1} related to the Biot-Savart law, which we present now. We use the spaces $H$, $H_c$, and $H_0$ of \cref{e:HSpace,e:HcH0}.

\begin{lemma}\label{L:TamingHc}
	Assume that $\Gamma$ is $C^{n, \al}$-regular and
	let $X$ be any function space contained in $C^{n, \al}(\Omega)^3$.
	For any $\vv \in H$,
	$
		\norm{P_{H_c} \vv}_X \le C(X) \norm{\vv}_H
	$.
\end{lemma}

\begin{lemma}\label{L:VortGivesUniqueness}
    If $\uu \in H$ with $\curl \uu = 0$ and $P_H \uu = 0$ then $\uu = 0$.
\end{lemma}

For any $\bomega$ in the range of the curl, $\curl (H^1(\Omega)^3)$, there exists a unique $\uu = K[\bomega] \in H_0 \cap H^1(\Omega)^3$ for which $\curl \uu = \bomega$. The operator $K$, which recovers the unique divergence-free vector field in $H_0$ having a given curl, encodes the Biot-Savart law.

There exists a vector field $\VV$ as regular as $\UU$ with $\dv \VV = 0$, $\curl \VV = 0$, and $\VV \cdot \n = U^\n$ on $[0, T] \times \Gamma$. We define
\begin{align}\label{e:KUDef}
	K_{U^\n}[\bomega]
		:= K[\bomega] + \VV.
\end{align}

\begin{lemma}\label{L:BSLaw}
    Assume $\UU \in \uSolSpaceN$. Let $\bomega$ be a divergence-free vector field on $\Omega$ having vanishing external fluxes. Let $\uu_c \in H_c$ and set $\uu = K_{U^\n}[\bomega] + \uu_c$. For all $t \in [0, T]$ and all integers $k$ with $0 \le k \le N$,
    \begin{align*}
        \norm{\uu(t)}_{W^{k + 1, p}(\Omega)}
            &\le C \norm{\bomega(t)}_{W^{k, p}(\Omega)}
                + \norm{\UU(t)}_{W^{k + 1, p}(\Omega)}
                + \norm{\uu_c(t)}_{W^{k + 1, p}(\Omega)}, \\
        \norm{\uu(t)}_{\uInitSpaceN}
            &\le C \norm{\bomega(t)}_{\vortSpaceNInit} + \norm{\UU(t)}_{C^{N + 1, \al}(\Omega)}
                + \norm{\uu_c(t)}_{C^{N + 1, \al}(\Omega)}
    \end{align*}
    for all $p \in (1, \iny)$, whenever the norms on the right-hand side are finite. In each case, $\UU$ can be replaced by $\VV$ and the final term can be replaced by $C \norm{\uu}_H$.
\end{lemma}
\begin{proof}
    For the two inequalities see, for instance, \cite{Paper1}. \cref{L:TamingHc} allows us to replace each of the final terms by $C \norm{\uu}_H$.
\end{proof}

\subsection*{A key property of $\bm{S^{N + 1, \al}}$} The purpose of this subsection is to prove the following:

\begin{prop}\label{P:SConsequences}
    Assume that $\uu \in \uSolSpaceN$. Then
    \begin{align*}
        \norm{\uu}_{C^{N, \al}(Q)}
            &\le C \sum_{j = 0}^N \norm{\prt_t^j \uu|_{t = 0}}_{C^{N - j, \al}(\Omega))}
                + C \brac{\norm{\uu}_\uSolSpaceN + \norm{\UU}_{S^{N+ 1, \al}}}
                    \max \set{T^\al, T^{1 - \al}}.
    \end{align*}
    Moreover, if $\uu \in \Dom_N(A)$ then
    \begin{align*}
        \norm{\uu}_{C^{N, \al}(Q)}
            &\le c_0 + c_X \max \set{T^\al, T^{1 - \al}}  \le \cN,
    \end{align*}
    where $c_0$, $c_X$, and $\cN$ are as in \cref{D:cc}.
\end{prop}

To prove \cref{P:SConsequences}, we will make use of the following space:
\begin{definition}\label{D:Rkal}
For an integer $k \ge 0$ and $\al \in (0, 1)$, define the space,
\begin{align*}
    R^{k, \al}
        &:= \set{f \in C^{k - 1}(Q):
            \prt_t^j f \in L^\iny([0, T]; C^{k - j, \al}
            (\Omega)), 0 \le j \le k}, \\
    \norm{f}_{R^{k, \al}}
        &= \sum_{j = 0}^k \norm{\prt_t^j f}_{{L^\iny([0, T]; C^{k - j, \al}(\Omega))}}.
\end{align*}
We allow $R^{k, \al}$ to apply to scalar-, vector-, or matrix-valued functions.
\end{definition}

\begin{lemma}\label{L:ROneLess}
    Let $f \in R^{k, \al}$, $k \ge 1$. Then
    \begin{align*}
        \norm{f}_{R^{k - 1, \al}}
            &\le  \sum_{j = 0}^{k - 1} \norm{\prt_t^j f|_{t = 0}}_{C^{k -1 - j, \al}(\Omega)}
                + \norm{f}_{R^{k, \al}} T.     
    \end{align*}
\end{lemma}
\begin{proof}
    Let $j$ be an integer with $0 \le j \le k - 1$. Because 
    \begin{align*}
        \prt_t^j f(t, \x)
            &= \prt_t^j f(t, \x)|_{t= 0}
                + \int_0^t \prt_t^{j + 1} f(s, \x) \, ds,
    \end{align*}
    we have
    \begin{align*}
        \norm{\prt_t^j D^\beta f}_{L^\iny([0, T]; C^{k -1 - j, \al}(\Omega))}
            &\le \norm{\prt_t^j D^\beta f|_{t= 0}}_{C^{k -1 - j, \al}(\Omega)}
                + \norm{\prt_t^{j + 1} D^\beta f}_{L^\iny([0, T]; C^{k - j, \al}(\Omega))} T.
    \end{align*}
    Summing over $j$ from $0$ to $k - 1$ gives the result.
\end{proof}

\begin{lemma}\label{L:Rkal}
    If $f \in R^{k, \al}$ for $k \ge 1$ then $f \in C^{k - 1, \al}(Q)$ with
    \begin{align}\label{e:RkGrowthBound}
        \norm{f}_{C^{k - 1, \al}(Q)}
            &\le \sum_{j = 0}^{k - 1} \norm{\prt_t^j f|_{t = 0}}_{C^{k -1 - j, \al}(\Omega)}
                + 4 \norm{f}_{R^{k, \al}} \max \set{T^{1 - \al}, T}.
    \end{align}
\end{lemma}
\begin{proof}
    First suppose that $k = 1$.
    Since $\prt_t f \in L^\iny([0, T]; C^\al(\Omega))$ and $f(t, \cdot) \in C(\ol{\Omega})$ for all $t \in [0, T]$, we have that for all $(t, \x) \in [0, T] \times \Omega$,
    \begin{align*}
        f(t, \x)
            &= f(0, \x) + \int_{0}^t \prt_s f(s, \x) ds.
    \end{align*}
    Then for all $t_1, t_2 \in [0, T]$ and $x_1, x_2 \in \Omega$,
    \begin{align}\label{e:fDoubleDiff}
        \begin{split}
        \abs{f(t_1, \x_1) - f(t_2, \x_2)}
            &\le \abs{f(t_1, \x_1) - f(t_1, \x_2)}
                + \abs{f(t_1, \x_2) - f(t_2, \x_2)} \\
            &\le \norm{f(t_1, \cdot)}_{\dot{C}^\al} \abs{\x_1 - \x_2}^\al
                + \norm{\prt_t f}_{L^\iny([0, T]; C^\al(\Omega))} \abs{t_1 - t_2}.
        \end{split}
    \end{align}
    Dividing both sides by $\abs{(t_1, \x_1) - (t_2, \x_2)}^\al$, which we note is greater than both $\abs{\x_1 - \x_2}^\al$ and $\abs{t_1 - t_2}$, we see that, in fact,  $f \in C^\al(Q)$ with $\norm{f}_{C^\al(Q)} \le \norm{f}_{R^{1, \al}}$.

    Moreover, we can estimate the term $\abs{f(t_1, \x_1) - f(t_1, \x_2)}$ in two other ways. First, because $f \in L^\iny([0, T]; C^{1, \al}(\Omega))$, we have
    \begin{align*}
        \abs{f(t_1, \x_1) - f(t_1, \x_2)}
            &\le a_1
            := \norm{f}_{L^\iny([0, T]; C^{1, \al}(\Omega))} \abs{\x_1 - \x_2}.
    \end{align*}
    Second, we have
    \begin{align*}
        \abs{f(t, \x_1) - f(t, \x_2)}
            &= \abs[\Bigg]{\int_0^t \prt_s f(s, \x_1) \, ds - \int_0^t \prt_s f(s, \x_2) \, ds}
            \le a_2
            := 2 \norm{\prt_t f}_{L^\iny([0, T] \times \Omega)} T.
    \end{align*}
    Hence,
    \begin{align*}
        \abs{f(t, \x_1) - f(t, \x_2)}
            &= \abs{f(t, \x_1) - f(t, \x_2)}^\al \abs{f(t, \x_1) - f(t, \x_2)}^{1 - \al}
            \le a_1^\al a_2^{1 - \al} \\
            &\le \norm{f}_{L^\iny([0, T]; C^{1, \al}(\Omega)}^\al
                2^{1 - \al} \norm{\prt_t f}_{L^\iny([0, T] \times \Omega)}^{1 - \al}
                \abs{\x_1 - \x_2}^\al T^{1 - \al} \\
            &\le 2 \norm{f}_{R^{1, \al}}
            \abs{\x_1 - \x_2}^\al T^{1 - \al}.
    \end{align*}
    Then, as in \cref{{e:fDoubleDiff}}, and using that
    $\abs{t_1 - t_2} \le \abs{t_1 - t_2}^\al T^{1 - \al}$,
    \begin{align*}
       \abs{f(t_1, \x_1) - f(t_2, \x_2)}
            &\le 2 \norm{f}_{R^{1, \al}}
            \abs{\x_1 - \x_2}^\al T^{1 - \al}
                + \norm{\prt_t f}_{L^\iny([0, T]; C^\al(\Omega))} \abs{t_1 - t_2}^\al T^{1 - \al}.
    \end{align*}
    Dividing both sides by $\abs{(t_1, \x_1) - (t_2, \x_2)}^\al$ and taking the supremum over all $(t_1, \x_1) \ne (t_2, \x_2)$ yields
    \begin{align*}
        \norm{f}_{\dot{C}^\al}
            &\le 3 \norm{f}_{R^{1, \al}} T^{1 - \al}.
    \end{align*}

    Also, for any $(t, \x) \in Q$,
    \begin{align*}
        \abs{f(t, \x)}
            &\le \abs{f(0, \x)} + \abs{f(t, \x) - f(0, \x)}
            \le \abs{f(0, \x)} + \norm{\prt_t f}_{L^\iny([0, T]; C^\al(\Omega))} \abs{t}
    \end{align*}
    so
    \begin{align*}
        \norm{f}_{L^\iny([0, T] \times \Omega)}
            &\le \norm{f(0)}_{L^\iny(\Omega)}
                + \norm{\prt_t f}_{L^\iny([0, T]; C^\al(\Omega))} T.
    \end{align*}
    Combined, these estimates yield \cref{e:RkGrowthBound} for $k = 1$.

    The result for $k > 1$ follows from applying the above argument to $D^\beta f$ for any $\abs{\beta} = k$, and controlling all the lower-order derivatives via \cref{L:ROneLess}.
\end{proof}
\begin{proof}[\textbf{Proof of \cref{P:SConsequences}}]
    First we prove that $\uu \in R^{N + 1, \al}$ with
    \begin{align}\label{e:PreliminaryRBound}
         \norm{\uu}_{R^{N + 1, \al}}
                    \le C \norm{\uu}_\uSolSpaceN + \norm{\UU}_{S^{N+ 1, \al}}.
    \end{align}
    We have,
    \begin{align*}
        \prt_t^j \uu
            &= \prt_t^j K_{U^\n}[\bomega] + \prt_t^j \uu_c
            =  K_{U^\n}[\prt_t^j \bomega] + \prt_t^j \uu_c,
    \end{align*}
    so by \cref{L:BSLaw}, for $j \le N$ and all $t \in [0, T]$,
    \begin{align*}
        \norm{\prt_t^j \uu(t)}_{C^{N + 1- j, \al}(\Omega)}
            &\le C \norm{\prt_t^j \bomega(t)}_{C^{N - j, \al}(\Omega)}
            + \norm{\UU(t)}_{C^{N + 1 - j, \al}(\Omega)}
            + C \norm{\prt_t^j \uu(t)}_H \\
            &\le C \norm{\bomega}_{C^{N, \al}(Q)}
                + C \norm{\uu}_{C^{N, \al}(Q)}
                + \norm{\UU(t)}_{C^{N, \al}(Q)} \\
            &\le C \norm{\uu}_{S^{N+ 1, \al}}
                + \norm{\UU(t)}_{C^{N, \al}(Q)},
    \end{align*}
    where we used \cref{L:TamingHc}. For $j = N + 1$, the $L^\iny([0, T]; C^\al(\Omega))$ norm is included as part of both the $R^{N +1, \al}$ and $\uSolSpaceN$ norms, and we see that \cref{e:PreliminaryRBound} follows.
    
    From \cref{L:Rkal} and \cref{e:PreliminaryRBound}, then,
    \begin{align*}
        \norm{\uu}_{C^{N, \al}(Q)}
            &\le \sum_{j = 0}^N \norm{\prt_t^j \uu|_{t = 0}}_{C^{N - j, \al}(\Omega)}
                + 4 \norm{\uu}_{R^{N + 1, \al}} \max \set{T^{1 - \al}, T}. \\
            &\le \sum_{j = 0}^N \norm{\prt_t^j \uu|_{t = 0}}_{C^{N - j, \al}(\Omega)}
                + C \brac{\norm{\uu}_\uSolSpaceN + \norm{\UU}_{S^{N+ 1, \al}}}
                    \max \set{T^\al, T^{1 - \al}},
    \end{align*}
    giving the first bound on $\norm{\uu}_{C^{N, \al}(Q)}$. If $\uu \in \Dom_N(A)$, then
    \begin{align*}
        \sum_{j = 0}^N \norm{\prt_t^j \uu|_{t = 0}}_{C^{N - j, \al}(\Omega)}
            &= \sum_{j = 0}^N \norm{\wprt^j \uu_0}_{C^{N - j, \al}(\Omega)}
            = c_0,
    \end{align*}
    giving the second bound on $\norm{\uu}_{C^{N, \al}(Q)}$.
\end{proof}

We now have the tools needed to prove \cref{L:DomNA}:

\begin{lemma}\label{L:DomNA}
	Assuming $\cond_N$ holds,
	$\Dom_N(A)$ is a nonempty, convex subset of $\uInputSpaceN$.
\end{lemma}
\begin{proof}
	We first show that $\Dom_N(A)$ is convex.
	Let $a, b \in [0, 1]$ with $a + b = 1$, let $\vv$, $\ww$ be in $\Dom_N(A)$,
	and let $\uu = a \vv + b \ww$. Then $\uu(0) = a \uu_0 + b \uu_0 = \uu_0$,
	and so also $\cond_N$ is satisfied.
	Similarly,
	$\prt_t^n \uu|_{t = 0}
	= a \prt_t^n \vv|_{t = 0} + b \prt_t^n \ww|_{t = 0} 
	= a \wprt^n \uu_0 + b \wprt^n \uu_0
	= \wprt^n \uu_0$. It follows that $\Dom_N(A)$ is convex.	
	
	To show that $\Dom_N(A)$ is nonempty, let $\bomega_0 := \curl \uu_0$ and define
	\begin{align*}
		\bomega(t)
			:= \bomega_0 + \sum_{n = 1}^N \frac{t^n}{n!} \wprt^n \bomega_0,
	\end{align*}
	so that for all $0 \le n \le N$, $\prt_t^n \bomega(0) = \wprt^n \bomega_0$.
	Because $\bomega(t)$ is in the range of the curl for all $t \in [0, T]$
	by \cref{P:curlwprt}, we can define
	\begin{align*}
		\uu(t)
			:= K_{U^\n}[\bomega] + \sum_{n = 0}^N \frac{t^n}{n!} P_{H_c} \wprt^n \uu_0,
	\end{align*}
	which we note lies in $\uSolSpaceN$.
	Then $\uu(0)$ and $\uu_0$ have the same $\curl$ and same harmonic component,
        and $\uu(0) - \uu_0 \in H$, so $\uu(0) = \uu_0$ by \cref{L:VortGivesUniqueness}.
        Moreover, for $1 \le n \le N$,
	\begin{align*}
		\curl \prt_t^n \uu(0)
			 &= \prt_t^n \bomega(0)
			 = \wprt^n \bomega_0
			 = \curl \wprt^n \uu_0
	\end{align*}
	by \cref{P:curlwprt}. Also, $P_{H_c} \prt_t^n \uu(0) = P_{H_c} \wprt^n \uu_0$.
	That is, $\prt_t^n \uu(0)$ and $\wprt^n \uu_0$ have the same curl and same
        harmonic component, while
        $\prt_t^n \uu(0) \cdot \n= \wprt^n \uu_0 \cdot \n$ on $\Gamma$. Hence, it follows from \cref{L:VortGivesUniqueness} that
	$\prt_t^n \uu(0) = \wprt^n \uu_0$, and we see that $\uu \in \Dom_N(A)$, demonstrating
	that $\Dom_N(A)$ is nonempty.
\end{proof}

%
%
\section{Flow map estimates}\label{S:FlowMap}

\noindent The pushforward of the initial vorticity by the flow map meets, along a hypersurface $\Cal{S}$ in $Q$, the pushforward of the vorticity generated on the inflow boundary. This requires some analysis at the level of the flow map. For the most part, the analysis in \cite{Paper1}, which we summarize here, suffices. The coarse bounds developed on the flow map in \cite{Paper1}, however, would only be sufficient for us to obtain small data existence of solutions: for the short time result for general data that we desire, we will require more explicit and refined bounds, which we develop in \cref{L:QuickFlow}.

We assume throughout this section that
$
	\UU \in \USpaceN,
	\uu \in \uInputSpaceN
$
for some $N \ge 0$. As in \cite{Paper1}, we extend $\uu$ to be defined on all of $\R \times \R^3$ using an extension operator like that in Theorem $5'$, chapter VI of \cite{Stein1970}. This extension need not be divergence-free, and is used only as a matter of convenience in stating results; it is only the value of $\uu$ on $\ol{Q}$ that ultimately concerns us.

We define $\eta \colon \R \times \R \times \R^3 \to \R^3$ to be the unique flow map for $\uu$, so that
$\prt_{t_2} \eta(t_1, t_2; \x) = \uu(t_2, \eta(t_1, t_2; \x))$. Then $\eta(t_1, t_2; \x)$ is the position that a particle starting at time $t_1$ at position $\x \in \R^3$ will be at time $t_2$ as it moves under the action of the velocity field $\uu$.

For any $(t, \x) \in \ol{Q}$ let
\begin{itemize}
	\item
		$\pos(t, \x)$ be the point on $\Gamma_+$ at which the flow line through $x$ at time $t$
		intersects $\Gamma_+$;
		
	\item
		$\time(t, \x)$ be the time at which that intersection occurs.
\end{itemize}
For our purposes, we can leave $\time$ and $\pos$ undefined if the flow line never intersects with $\Gamma_+$.

\begin{remark}\label{R:Brevity}
	We will often drop the $(t, \x)$ arguments on $\time$ and $\pos$
	for brevity.
\end{remark}

We define the hypersurface,
\begin{align*}
	\Cal{S} := \set{(t, \x) \in \ol{Q} \colon \time(t, \x) = 0},
\end{align*}
which is nonempty since it contains at least $\Gamma_+ \times \set{0}$,
and the open sets $U_\pm \subseteq Q$,
\begin{align*}
	U_-
		&:= \set{(t, \x) \in Q
				\colon (t, \x) \notin \text{ domain of } \time, \pos}, \\
	U_+
		&:= \set{(t, \x) \in Q
				\colon \time(t, \x) > 0}.
\end{align*}
Then $\Cal{S}$ is of class $C^{N + 1, \al}$ as a hypersurface in $Q$ and $\Cal{S}(t) := \set{\x \in \Omega \colon (t, \x) \in \Cal{S}}$ is of class $C^{N + 1, \al}$ as a surface in $\Omega$. 

The estimates on the flow map in \cref{L:QuickFlow} are more explicit than in \cite{Paper1}, where we required only coarse estimates.  We note that $\eta$ has one more derivative in both time variables than has $\uu$, which we can see in the explicit estimates.
\begin{lemma}\label{L:QuickFlow}
	The flow map
	$\eta \in C^{N + 1, \al}([0, T]^2 \times \R^3)$.
	Define $\mu \colon U_+ \to [0, T] \times \Gamma_+$ by
	\begin{align*}
		\mu(t, \x) = (\time(t, \x), \pos(t, \x))
	\end{align*}
    and let $M := \norm{\uu}_{S^{1, \al}}$.
	The functions $\time$, $\pos$, $\mu$ lie in $C^{N + 1, \al}(\olUPlus)$.
	Moreover,
	\begin{align}\label{e:etaBounds}
		\begin{split}
		\norm{\prt_{t_1} \eta(t_1, t_2; \x)}_{L^\iny_\x}
			&\le \norm{\uu}_{L^\iny(Q)} h(t_1, t_2),
				\\
		\norm{\grad \eta(t_1, t_2; \x)}_{L^\iny_\x}
			&\le h(t_1, t_2),
				\\
		\norm{\grad \eta(0, t_2; \x)}_{\dot{C}^\al_{t_2}(Q)}
			&\le \norm{\grad \uu}_{L^\iny(Q)} h(0, T) T^{1 - \al},
				\\
		\norm{\grad \eta(0, t_2; \x)}_{\dot{C}^\al_\x(Q)}
			&\le h(0, T)^{1 + 2 \al} \int_0^T 
				\norm{\grad \uu(s)}_{\dot{C}^\al} \, ds,
					\\
		\norm{\grad \eta(0, T; \x)}_{\dot{C}^\al(Q)}
			&\le e^{(1 + 2 \al) M T} M T^{1 - \al},
		\end{split}
	\end{align}
	where
	\begin{align*}
		h(t_1, t_2)
			:= \exp \abs[\Bigg]{\int_{t_1}^{t_2} \norm{\grad \uu(s)}_{L^\iny} \, ds}
			\le e^{M T}.
	\end{align*}
	Also,
	\begin{align}\label{e:KeyBound}
		\norm{D \mu}_{L^\iny(U_+)}
			&\le C U_{min}^{-1} [1 + \norm{\uu}_{L^\iny(Q)}^2] h(0, T),
	\end{align}
	where $U_{min}$ is as in \cref{D:NReg}.

	More generally, for any $N \ge 0$, defining $\exp^n$ to be $\exp$
	composed with itself $n$ times,
	\begin{align}\label{e:etaBoundsN}
		\begin{split}
		\norm{\prt_{t_1}^{N + 1} \eta(t_1, t_2; \x)}_{L^\iny([0, T]^2 \times \Omega)}
			&\le C \norm{\uu}_{C^N(Q)} \exp^{N + 1}(M T),
				\\
		\norm{\grad^{N + 1} \eta(t_1, t_2; \x)}_{L^\iny([0, T]^2 \times \Omega)}
			&\le \exp^{N + 1}(M T),
				\\
		\norm{\grad^{N + 1} \eta(0, t_2; \x)}_{\dot{C}^\al_{t_2}(Q)}
			&\le \norm{\grad^{N + 1} \uu}_{L^\iny(Q)} \exp^{N + 1}(M T) T^{1 - \al},
				\\
		\norm{\grad^{N + 1} \eta(0, t_2; \x)}_{\dot{C}^\al_\x(Q)}
			&\le \exp^{N + 1} (C M T) \int_0^T 
				\norm{\grad^{N + 1} \uu(s)}_{\dot{C}^\al} \, ds,
					\\
		\norm{\grad^{N + 1} \eta(0, T; \x)}_{\dot{C}^\al(Q)}
			&\le \exp^{N + 1} (C M T)  M T^{1 - \al}, \\
		\norm{D^{N + 1} \mu}_{L^\iny(U_+)}
			&\le c_0 [1 + \norm{\uu}_{C^N(Q)}^{2(N + 1)}] \exp^{N + 1}(MT).
		\end{split}
	\end{align}
\end{lemma}
\begin{proof}
	We will apply \cref{L:HolderComp}
	multiple times without explicit reference.
	
	Taking the gradient of the integral expression in \LinearPaper{(3.1)},
	\begin{align}\label{e:gradeta}
		\grad \eta(t_1, t_2; \x)
			= I + \int_{t_1}^{t_2} \grad \uu(s, \eta(t_1, s; \x)) \grad \eta(t_1, s; \x) \, ds.
	\end{align}
	Thus,
	\begin{align*}
		\norm{\grad \eta(t_1, t_2; \x)}_{L^\iny_\x}
			&\le 1 + \abs[\Bigg]{\int_{t_1}^{t_2} \norm{\grad \uu(s)}_{L^\iny}
				\norm{\grad \eta(t_1, s; \x)}_{L^\iny_\x}
				\, ds}.
	\end{align*}
	\Gronwalls Lemma, applied with fixed $t_1$, gives \cref{e:etaBounds}$_2$.
	\LinearPaperCap{Lemma 3.1} gives
	$\prt_{t_1} \eta(t_1, t_2; \x) = - \uu(t_1, \x) \cdot \grad \eta(t_1, t_2; \x)$,
	from which \cref{e:etaBounds}$_1$ follows.
	
	It also follows from \cref{e:gradeta} that
	\begin{align*}
		\norm{\grad \eta(0, t_2; \x)}_{\dot{C}(Q)^\al_{t_2}}
			&\le \sup_{t_2 \ne t_2'}
				\frac{\norm{\grad \uu}_{L^\iny(Q)}
				\norm{\grad \eta}_{L^\iny(Q)}}
				{\abs{t_2 - t_2'}^\al}
				\abs{t_2 - t_2'} \\
			&\le \norm{\grad \uu}_{L^\iny(Q)} h(0, T) T^{1 - \al},
	\end{align*}
	giving \cref{e:etaBounds}$_3$.
	
	Returning once more to \cref{e:gradeta},
	\begin{align*}
		\norm{\grad \eta(t_1, t_2; \x)}_{\dot{C}^\al_\x}
			&\le \int_0^{t_2} \norm{\grad \uu(s, \eta(t_1, s; \x))
				\grad \eta(t_1, s; \x)}_{\dot{C}^\al_\x} \, ds.
	\end{align*}
	But, using \cref{L:HolderProd},
	\begin{align*}
		&\norm{\grad \uu(s, \eta(t_1, s; \x))
				\grad \eta(t_1, s; \x)}_{\dot{C}^\al_\x} \\
			&\qquad
			\le \norm{\grad \uu(s, \eta(t_1, s; \x))}_{\dot{C}^\al_\x}
					\norm{\grad \eta(t_1, s; \x)}_{L^\iny_\x}
				+
				\norm{\grad \uu(s, \eta(t_1, s; \x))}_{L^\iny_\x}
					\norm{\grad \eta(t_1, s; \x)}_{\dot{C}^\al_\x}
				\\
			&\qquad
			\le \norm{\grad \uu(s)}_{\dot{C}^\al} 
				\norm{\eta(t_1, s; \x)}_{Lip_\x}^\al
				\norm{\grad \eta(t_1, s; \x)}_{L^\iny_\x}
				+
				\norm{\grad \uu(s)}_{L^\iny}
					\norm{\grad \eta(t_1, s; \x)}_{\dot{C}^\al_\x}
				\\
			&\qquad
			\le \norm{\grad \uu(s)}_{\dot{C}^\al}
				h(t_1, s)^{2 \al}
				+
				\norm{\grad \uu(s)}_{L^\iny}
					\norm{\grad \eta(t_1, s; \x)}_{\dot{C}^\al_\x},
	\end{align*}
	so
	\begin{align*}
		&\norm{\grad \eta(0, t_2; \x)}_{\dot{C}^\al_\x} \\
			&\qquad
			\le \int_0^{t_2} 
				\norm{\grad \uu(s)}_{\dot{C}^\al} h(0, s)^{2 \al} \, ds
				+ \int_0^{t_2} \norm{\grad \uu(s)}_{L^\iny(\Omega)}
				\norm{\grad \eta(0, s; \x)}_{\dot{C}^\al_\x}
				\, ds \\
			&\qquad
			\le h(0, t_2)^{2 \al} \int_0^{t_2} 
				\norm{\grad \uu(s)}_{\dot{C}^\al} \, ds
				+ \int_0^{t_2} \norm{\grad \uu(s)}_{L^\iny(\Omega)}
				\norm{\grad \eta(0, s; \x)}_{\dot{C}^\al_\x}
				\, ds.
	\end{align*}
	Applying \Gronwalls Lemma gives
	\begin{align*}
		\begin{split}
		\norm{\grad \eta(0, t_2; \x)}_{\dot{C}^\al_\x}
			&\le \brac{h(0, t_2)^{2 \al} \int_0^{t_2} 
				\norm{\grad \uu(s)}_{\dot{C}^\al} \, ds}
				\exp \int_0^{t_2} \norm{\grad \uu(s)}_{L^\iny(\Omega)}
				\, ds \\
			&= h(0, t_2)^{1 + 2 \al} \int_0^{t_2} 
				\norm{\grad \uu(s)}_{\dot{C}^\al} \, ds,
		\end{split}
	\end{align*}
	which is \cref{e:etaBounds}$_4$.

	From \LinearPaper{Lemma 3.5},
	\begin{align}\label{e:timeposidentities}
		\arraycolsep=0.5em\def\arraystretch{1.2}
		\begin{array}{ll}
		\prt_t \time 
			= - U^\n(\time, \pos)^{-1} \prt_{t_1} \eta(t, \time; \x) \cdot \n(\pos),
		&\grad \time
			= - U^\n(\time, \pos)^{-1} (\grad \eta(t, \time; \x))^T \n(\pos),
				\\
		\prt_t \pos
			= 
			\prt_{t_1} \eta(t, \time; \x)
			+ \prt_t \time \uu(\time, \pos),
		&\grad \pos
			= \uu(\time, \pos) \otimes \grad \time + \grad \eta(t, \time; \x).
		\end{array}
	\end{align}
	We use these expressions to calculate,
	\begin{align*}
		\norm{\prt_t \time}_{L^\iny(U_+)}
			&\le C U_{min}^{-1} \norm{\prt_{t_1} \eta}_{L^\iny(Q)}
			\le C U_{min}^{-1} \norm{\uu}_{L^\iny(Q)} h(0, T),
				\\
		\norm{\grad \time}_{L^\iny}
			&\le C U_{min}^{-1} \norm{\grad \eta}_{L^\iny(Q)}
			\le C U_{min}^{-1} h(0, T),
				\\
		\norm{\prt_t \pos}_{L^\iny}
			&\le C U_{min}^{-1} \norm{\prt_{t_1} \eta}_{L^\iny(Q)}
				+ \norm{\uu}_{L^\iny(Q)} \norm{\prt_t \time}_{L^\iny} \\
			&\le C U_{min}^{-1} [\norm{\uu}_{L^\iny(Q)} + \norm{\uu}_{L^\iny(Q)}^2] h(0, T),
				\\
		\norm{\grad \pos}_{L^\iny}
			&\le \norm{\uu}_{L^\iny(Q)}
					\norm{\grad \time}_{L^\iny}
					+ \norm{\grad \eta}_{L^\iny(Q)}
			\le [1 + C U_{min}^{-1} \norm{\uu}_{L^\iny(Q)}] h(0, T).
	\end{align*}
	Summing these estimates gives the bound on $D \mu = (\prt_t \mu, \grad \mu)$.
	
	The bounds for higher $N$ follow from inductive extension of these arguments.
\end{proof}

\begin{remark}\label{R:cNFlowBounds}
    The exact bounds in \cref{L:QuickFlow} are not so important, but it is important that $M$ only appear in the exponentials. Because of that and \cref{P:SConsequences}, we see that for $\uu \in \Dom_N(A)$, each of the bounds in \cref{e:etaBounds,e:KeyBound,e:etaBoundsN} is of the form $\cN$ of \cref{D:cc}.
    Similarly, $\time$ and $\pos$ can be bounded in $C^{N, \al}(U_+)$ by $\cN$.
\end{remark}

We are now in a position to give the definition of a Lagrangian solution to \OurLinear, as it appears in \cite{Paper1}. For this purpose, define
\begin{align}\label{e:pos0}
	\pos_0
		= \pos_0(t, \x)
		:= \eta(t, 0; \x).
\end{align}
As with $\time$ and $\pos$ (see \cref{R:Brevity}) we will often drop the $(t, \x)$ arguments on $\pos_0$.

\begin{definition}[Lagrangian solution to \OurLinear]\label{D:LagrangianSolution}
Define $\ol{\bomega}_\pm$ and $\G_\pm$ on $U_\pm$ by
\begin{align}\label{e:LagrangianForm}
	\begin{split}
	\ol{\bomega}_-(t, \x)
		&= \grad \eta(0, t; \pos_0) \ol{\bomega}_0(\pos_0) + \G_+(t, \x), \\
	\ol{\bomega}_+(t, \x)
		&= \grad \eta(\time, t; \pos) \H(\time, \pos) + \G_-(t, \x), \\
	\G_-(t, \x)
		&:= \int_0^t \grad \eta(s, t; \eta(t, s; \x))
				\g(s, \eta(t, s; \x)) \, ds, \\
	\G_+(t, \x)
		&:= \int_{\time(t, \x)}^t \grad \eta(s, t; \eta(t, s; \x))
				\g(s, \eta(t, s; \x)) \, ds.
	\end{split}
\end{align}
Then $\ol{\bomega}$ defined by $\ol{\bomega}|_{U_\pm} = \ol{\bomega}_\pm$ is called a Lagrangian solution to \OurLinear.
\end{definition}

\begin{remark}\label{R:Paper1}
	A few words are appropriate here about the treatment of \cref{e:LagrangianForm} in \cite{Paper1}.
	
	The outflow vorticity, $\ol{\bomega}_-$, behaves and can be analyzed much like
	the full vorticity in the classical
	setting of an impermeable boundary: the initial vorticity is pushed forward from time zero
	by the flow map $\eta$ for the given velocity field, $\uu$, and Duhamel's principle is used
	to treat the forcing term, $\G_-$. Such an analysis yields
	$\ol{\bomega}_- \in C^{N, \al}(\Omega_-)$.
	
	The inflow vorticity, $\ol{\bomega}_+$, is somewhat more complicated since the
	vorticity is pushed off the 2D inflow boundary $\Gamma_+$ into the 3D domain $\Omega$,
	producing an inflow component $\Omega_+(t)$ expanding in time. Also, the regularity
	of $\time(t, \x)$ and $\pos(t, \x)$ must be accounted for, and the time $t$ enters
	into the Duhamel integral in both limits.  Nevertheless, 
	$\ol{\bomega}_+ \in C^{N, \al}(\Omega_+)$ holds.
	
	The key difficulty, however, lies not with the inflow or outflow vorticity individually,
	but rather
	with insuring that they meet across the hypersurface $\Cal{S}$ in a manner that allows
	the full vorticity $\ol{\bomega}$ to be regular enough to lie in $C^{N, \al}(\Omega)$.
	That $\lincond_N$ is the right condition to insure this is natural and is easy to show
	for $N = 0$, primarily because no derivatives are involved to obtain $C^\al(Q)$ regularity.
	The situation for $N > 0$, is much more involved.
	
	To obtain regularity across $\Cal{S}$, the $N = 1$ case formally reduces to the $N = 0$
	case.
	But for $N \ge 2$, such a reduction to the $N - 1$ case can be obtained,
	and allows an induction argument to be made to reduce the problem to the $N = 1$
	case. This leaves the $N = 1$ case, which requires a delicate analysis.
	
	Complicating the argument slightly is that for $N \ge 1$, the two terms making up
	$\ol{\bomega}_\pm$ 
	in \cref{e:LagrangianForm} need not be $C^{N, \al}$-continuous across $\Cal{S}$, though their
	sum is.
	But as long as $\lincond_N$ holds, $C^{N, \al}$ estimates on each of the four terms
	making up $\omega^+$ and $\omega^-$ can be combined 
	to give estimates on $\ol{\bomega}$ in  $\vortSpaceN$.

\end{remark}

%
%
\section{The nonlinear term on the boundary}\label{S:Boundary}

\noindent
\cref{P:ugradunOnGammaPlus} gives coordinate-free expressions for $(\uu \cdot \grad \uu) \cdot \n$ on $\Gamma$. The proof of \cref{P:ugradunOnGammaPlus} is most readily obtained  using the boundary coordinates introduced in \cref{A:BoundaryCalculations}, so we defer it to that appendix.

\begin{definition}\label{D:perp}
	For any tangent vector field $\vv$ on $\Gamma$, define $\vv^\perp$
	to be $\vv$ rotated 90 degrees counterclockwise around the normal vector
	when viewed from outside $\Omega$
	(so $\vv^\perp = \n \times \vv$).
\end{definition}

We write the gradient and divergence on the boundary as $\grad_\Gamma$ and $\dv_\Gamma$, as in \cref{A:BoundaryCalculations}.

\begin{prop}\label{P:ugradunOnGammaPlus}
	Assume that $\Gamma$ is $C^2$.
	Let $\uu$ be a divergence-free differentiable vector field,
	let $u^\n = \uu \cdot \n$,
	and, as in \cref{e:ComponentConvention}, let $\uu^\BoldTau = \uu - u^\n \n$.
	Let $\kappa_1, \kappa_2$ be the principal curvatures on $\Gamma$.
	On $[0, T] \times \Gamma$, we have
	\begin{align}\label{e:ugradun}
		(\uu \cdot \grad \uu) \cdot \n
		&= -u^\n \dv_\Gamma \uu^\BoldTau + \uu^\BoldTau \cdot \grad_\Gamma u^\n
			- (\kappa_1 + \kappa_2) (u^\n)^2
			- \uu^\BoldTau \cdot \Cal{A} \uu^\BoldTau.
	\end{align}	
	Here, $\Cal{A}$ is the shape operator on the boundary: for any tangential vector field,
	$\Cal{A} \vv$ is the directional derivative of $\n$ in the direction of $\vv$,
	which is also a tangential vector field. 
\end{prop}

The nonlinear term on the boundary is key to recovering the pressure, as we will see in the next section. It was for these purposes that we used $N[\uu]$ given in \cref{e:NuDef} to define the regularized pressure in \cref{e:pr}.
Using that $\uu^\n = \UU^\n$, substituting the expression in \cref{e:ugradun} for $(\uu \cdot \grad \uu) \cdot \n$, and using \cref{e:dvfv}, we see that on $\Gamma_+$,
\begin{align}\label{e:NuNice}
	\begin{split}
		N[\uu]
			= -U^\n \dv_\Gamma \UU^\BoldTau + \uu^\BoldTau \cdot \grad_\Gamma U^\n
			- (\kappa_1 + \kappa_2) (U^\n)^2
			- \uu^\BoldTau \cdot \Cal{A} \uu^\BoldTau,
	\end{split}
\end{align}
so $N[\uu]$ has no derivatives on $\uu^\BoldTau$. Nonetheless, integrating \cref{e:NuDef}$_2$ by parts along each boundary component using \cref{L:IBPOnBoundary}, we see that
\begin{align}\label{e:NuDefuComp}
	\int_\Gamma N[\uu]
		= \int_\Gamma (\uu \cdot \grad \uu) \cdot \n.
\end{align}
Hence, replacing $(\uu \cdot \grad \uu) \cdot \n$ with $N[\uu]$ does not alter the compatibility condition for recovering the pressure, as in \cref{S:Pressure}.

%
%
\section{Pressure Estimates}\label{S:Pressure}

\noindent We can determine the pressure from the velocity as in \cref{e:TruePressureEq}.
On $\Gamma_0$, as we can see from \cref{e:ugradun}, $\grad p \cdot \n = - \uu^\BoldTau \cdot \Cal{A} \uu^\BoldTau$ ($= - \kappa \abs{\uu}^2$ in 2D, where $\kappa$ is the scalar curvature). Hence, when $\Gamma = \Gamma_0$, standard Schauder estimates imply that $\grad p$ and $\uu$ have the same spatial regularity. This is the impermeable boundary case. But for inflow, outflow boundary conditions, the expression for $\grad p \cdot \n$ contains spatial derivatives of $\uu$, as we can see from \cref{e:ugradun}, on which we have no a priori control. (Because $\uu \cdot \n = U^\n$ on all of $\Gamma$, the time derivative in \cref{e:TruePressureEq}$_2$ does not impact the regularity of $p$.)

We circumvent this difficulty using the simple but clever technique in \cite{AKM}: we replace the boundary condition in \cref{e:TruePressureEq}$_2$ using $N[\uu]$ of \cref{e:NuDef}, solving instead, \cref{e:pr} for the pressure $p_r$.
We see from \cref{e:NuDefuComp} that the required compatibility condition coming from $\int_\Gamma \grad p_r \cdot \n = \int_\Omega \Delta p_r = \int_\Omega \dv(-\prt_t \uu - \uu \cdot \grad \uu)$ remains satisfied when using $-\prt_t U^\n - N[\uu]$ in place of $-\prt_t u^\n - (\uu \cdot \grad \uu) \cdot \n$ on $\Gamma$.
For $\uu(t) \in C^{N + 1, \al}(\Omega)$, classical elliptic regularity theory gives a solution to \cref{e:pr}  with $p_r(t)$ in $C^{N + 1, \al}(\Omega)$, unique up to an additive constant. Ultimately, we show, in \cref{P:WeakVorticityIsVelocity}, that $p_r = p$ for the fixed point of the operator $A$.

 Since we are seeking solutions to \Our in \Holder spaces, it would seem natural to use elliptic estimates in \Holder spaces. To obtain the needed control on our pressure estimates in time, however, we will find it necessary to use, instead, elliptic estimates in Sobolev spaces, as given in \cref{L:prLemma}. The reason for this is explained in \cref{R:ReasonForSobolev} following the proof of \cref{P:PressureBounds}.

\begin{lemma}\label{L:prLemma}
    Let $n \ge 0$ and $f \in W^{n, q}(\Omega)$, where $q \in (1, \iny)$.
    If $f$ has mean zero then
	\begin{align*}
		\begin{split}
		\norm{f}_{W^{n + 2, q}(\Omega)}
			&\le C \brac{
				\norm{\Delta f}_{W^{n, q}(\Omega)}
					+ \norm{\grad f \cdot \n}_{W^{n + 1 - \frac{1}{q}, q}(\Gamma)}
				}
		\end{split}
	\end{align*}
    and for any nonempty compact subset of $\Omega'$ of $\Omega \cup \Gamma_+$
    (without assuming $f$ has mean zero)
	\begin{align}\label{e:prLemma}
		\begin{split}
		\norm{f}_{W^{n + 2, q}(\Omega')}
			&\le C \brac{
				\norm{\Delta f}_{W^{n, q}(\Omega)}
					+ \norm{\grad f \cdot \n}_{W^{n + 1 - \frac{1}{q}, q}(\Gamma_+)}
					+ \norm{f}_{L^q(\Omega)}
				}.
		\end{split}
	\end{align}
\end{lemma}
\begin{proof}
	The bounds in \cref{e:prLemma} for $n = 0$ are stated near the bottom of page 174 of \cite{AKM},
	but let us say a few words about them. First, they are derived
	from combining an interior estimate away from all boundaries
	with an estimate that includes only $\Gamma_+$.
	Second, \cite{AKM} treats the $N = 0$ case, and we use
	(15.1.5) of \cite{Agranovich2015} for the $N \ge 1$ case.
\end{proof}

In what follows, we will use $L^q$-based Sobolev spaces on $\Omega$ with $q > 3/(1 - \al)$. This will give us some useful properties, which we summarize in \cref{L:Sobolev}.

\begin{lemma}\label{L:Sobolev}
    Let $q > 3/(1  - \al)$. Then for any integer $k \ge 1$,
    \begin{align}\label{e:SobolevEmbedding}
        W^{k + 1, q}(\Omega) \subseteq C^{k, \al}(\Omega)
    \end{align}
    and for any $r \in [1, \iny]$ and any $f \in W^{2,q}(\Omega)$,
    \begin{align}\label{e:BoundarySobolev}
        \begin{split}
        \norm{f}_{L^r(\Gamma)}
            &\le \norm{f}_{C^\al(\Gamma)}
            \le \norm{f}_{C^\al(\Omega)}
            \le C \norm{f}_{W^{1, q}(\Omega)}, \\
        \norm{\grad f}_{L^q(\Gamma)}
            &\le \norm{\grad f}_{W^{1 - 1/q, q}(\Gamma)}
            \le C \norm{\grad f}_{W^{1, q}(\Omega)}
            \le C \norm{f}_{W^{2, q}(\Omega)}.
        \end{split}
    \end{align}
    For $k \ge 1$, $W^{k, q}(\Omega)$ is an algebra, while for $k \ge 0$ its trace  space $W^{k - 1/q, q}(\Gamma)$ is an algebra.
\end{lemma}
\begin{proof}
    Sobolev embedding gives \cref{e:SobolevEmbedding}.
    The inequality in \cref{e:BoundarySobolev}$_1$ follows from \cref{e:SobolevEmbedding} for any $f \in C^\iny(\Omega) \cap W^{2,q}(\Omega)$, which is dense in $W^{1,q}(\Omega)$; \cref{e:BoundarySobolev}$_2$ follows from the trace inequality.

    For $k \ge 1$,  $kq > 3k \ge 3$ so $W^{k, q}(\Omega)$ is an algebra, and this same condition gives that $W^{k - 1/q, q}(\Gamma)$ is an algebra.
\end{proof}

Not only will we need estimates on $p_r$, but, letting $\uu_1, \uu_2 \in S^{1, \al}$, where $p_{r, 1}, p_{r, 2}$ solve \cref{e:pr} for $\uu_1$, $\uu_2$, respectively, we will need estimates for $N = 0$ on
$
    P
        := p_{1, r} - p_{2, r}
$.
Fixing $t_1, t_2 \in [0, T]$, we will also need to estimate
\begin{align*}
    \ol{p}_r := p_r(t_1) - p_r(t_2), \quad
    \ol{P} := P(t_1) - P(t_2).
\end{align*}

We start in \cref{P:prLqBounds,P:PressureBounds} by controlling only the spatial derivatives of $q$.

\begin{prop}\label{P:prLqBounds}
    Let $q > 3/(1 - \al)$, $t_1, t_2 \in [0, T]$, and let $p_r$ be the unique solution to \cref{e:pr} for some $\uu \in \uSolSpaceN$ normalized so that
    $M_q(p_r(t)) := \int_\Omega p_r \abs{p_r}^{q - 2} = 0$.
    Then
	\begin{align}\label{e:prLqBound}
		\begin{split}
		\norm{p_r(t)}_{L^q(\Omega)}
			&\le C_1,
		\end{split}
	\end{align}
	where
	\begin{align*}
		C_1
			&:=
                        \norm{\UU}_{{L^\iny(Q)}}^2
                        + \norm{\prt_t \UU}_{{L^\iny(Q)}}
					+ \norm{\uu}_{L^\iny(Q)}^2.
	\end{align*}
 
    Fixing $t_1, t_2 \in [0, T]$, normalize $p_r$ so that
    $M_q(p_r(t_1) - p_r(t_2)) = 0$. Then
	\begin{align}\label{e:prt1t2LqBounds}
		\begin{split}
		\norm{p_r(t_1) - p_r(t_2)}_{L^q(\Omega)}
			&\le C_2
				\abs{t_1 - t_2},
		\end{split}
	\end{align} 
	where
    \begin{align*}
		C_2
            &:= C \brac{\norm{\UU}_{S^{2, \al}}
                + \norm{\UU}_{S^{1, \al}} ^2
                + \norm{\uu}_{L^\iny(Q)} \norm{\uu}_{S^{1, \al}}},
	\end{align*}
	the constant $C$ depending only upon $\Omega$ and $q$.
\end{prop}
\begin{proof}
    We adapt the argument on pages 175-176 of \cite{AKM}.
    For now we suppress the time variable.
	
    Let $\beta$ be the unique mean-zero solution to
    \begin{align*}
        \begin{cases}
            \Delta \beta = p_r \abs{p_r}^{q - 2}
                &\text{ in } \Omega, \\
            \grad \beta \cdot \n = 0
                &\text{ on } \Gamma,
        \end{cases}
    \end{align*}
    where the normalization of $p_r$ gives solvability. Letting $q' = q/(q - 1)$, which we note is \Holder conjugate to $q$, \cref{L:prLemma} gives
    \begin{align*}
        \norm{\beta}_{W^{2, q'}(\Omega)}
            &\le C \norm{\abs{p_r}^{q - 1}}_{L^{q'}(\Omega)}
            = C \norm{p_r}_{L^q(\Omega)}^{q - 1}.
    \end{align*}

    From \cref{e:pr,e:NuDef}, the elliptic problem for $p_r$ can be written, for a fixed time, as
    \begin{align}\label{e:prElliptic}
        \begin{split}
    	\begin{cases}
    		\Delta p_r = - \dv (\uu \cdot \grad \uu)
    			&\text{in } \Omega, \\
    		\grad p_r \cdot \n = -\prt_t U^\n
                - (\uu \cdot \grad \uu) \cdot \n
    			- \CharFunc_{\Gamma_+} \dv_\Gamma (U^\n (\uu^\BoldTau - \UU^\BoldTau))
    			&\text{on } \Gamma,
    	\end{cases}
        \end{split}
    \end{align}  
    where $\CharFunc_{\Gamma_+}$ is the characteristic function for $\Gamma_+$.  
    Then, using \cref{e:prElliptic},
    \begin{align*}
        \begin{split}
        \norm{p_r}_{L^q(\Omega)}^q
            &= (\Delta \beta, p_r)
            = - (\grad \beta, \grad p_r) + \int_\Gamma (\grad \beta \cdot \n) p_r
            = (\Delta p_r, \beta) - \int_\Gamma (\grad p_r \cdot \n) \beta \\
             &= - (\dv (\uu \cdot \grad \uu), \beta)
                    - \int_\Gamma (\grad p_r \cdot \n) \beta \\
            &= (\uu \cdot \grad \uu, \grad \beta)
            - \int_\Gamma ((\uu \cdot \grad \uu) \cdot \n) \beta 
                - \int_\Gamma (\grad p_r \cdot \n) \beta \\
            &= (\uu \cdot \grad \uu, \grad \beta)
            - \int_{\Gamma_+} \pr{\prt_t U^n + \dv_\Gamma(U^\n(\uu^\BoldTau - \UU^\BoldTau))} \beta.
        \end{split}
    \end{align*}
    But,
    \begin{align*}
        (\uu \cdot &\grad \uu, \grad \beta)
            = \int_\Omega u^i \prt_i u^j \prt_j \beta
            = \int_\Omega u^i \prt_i (u^j \prt_j \beta)
                - \int_\Omega u^i u^j \prt_i \prt_j \beta \\
            &= (\uu, \grad (\uu \cdot \grad \beta))
                    - (\uu \otimes \uu, \grad \grad \beta)
            = - (\uu \otimes \uu, \grad \grad \beta)
                - \int_\Gamma U^\n (\uu \cdot \grad \beta)
    \end{align*}
    and, applying \cref{L:IBPOnBoundary},
    \begin{align*}
        - \int_{\Gamma_+} \dv_\Gamma(U^\n(\uu^\BoldTau - \UU^\BoldTau)) \beta
            &= \int_{\Gamma_+} U^\n (\uu^\BoldTau - \UU^\BoldTau) \cdot \grad_\Gamma \beta
            = \int_{\Gamma_+} U^\n (\uu - \UU) \cdot \grad \beta.
    \end{align*}
    Here, we used that $\grad \beta \cdot \n = 0$, so $\vv^\BoldTau \cdot \grad_\Gamma \beta = \vv \cdot \grad \beta$. Hence,
    \begin{align*}
        \begin{split}
        \norm{p_r}_{L^q(\Omega)}^q
        &= - (\uu \otimes \uu, \grad \grad \beta)
            - \int_\Gamma U^\n (\uu \cdot \grad \beta)       
            - \int_{\Gamma_+} \prt_t U^n \beta
            + \int_{\Gamma_+} U^\n (\uu - \UU) \cdot \grad \beta \\
        &= - (\uu \otimes \uu, \grad \grad \beta)
            - \int_{\Gamma_-} U^\n (\uu \cdot \grad \beta)       
            - \int_{\Gamma_+} U^\n \, \UU \cdot \grad \beta
            - \int_{\Gamma_+} \prt_t U^n \beta.
        \end{split}
    \end{align*}
    Exploiting \cref{e:BoundarySobolev}, we have the bound,
	\begin{align*}
		\norm{p_r}_{L^q(\Omega)}^q
			&\le \norm{\uu}_{L^\iny(\Omega)} \norm{\uu}_{L^q(\Omega)}
                    \norm{\beta}_{W^{2, q'}(\Omega)}
				+ \norm{\prt_t U^\n}_{L^q(\Gamma)} \norm{\beta}_{L^{q'}(\Gamma)} \\
            &\qquad
                + \norm{U^n}_{L^{q'}(\Gamma)}
                    \brac{\norm{\UU}_{L^\iny(\Gamma)} + \norm{\uu}_{L^\iny(\Gamma)}}
                    \norm{\grad \beta}_{L^q(\Gamma)} \\
            &\le C_1 \norm{\beta}_{W^{2, q'}(\Omega)}
            \le C_1 \norm{p_r}_{L^q(\Omega)}^{q - 1},
	\end{align*}
	from which \cref{e:prLqBound} follows.
 
	To obtain \cref{e:prt1t2LqBounds} we argue the same way,
	bounding $\ol{p}_r := p_r(t_1) - p_r(t_2)$, where now $\beta$ is the
	unique mean-zero solution to
    \begin{align*}
        \begin{cases}
            \Delta \beta = \ol{p}_r \abs{\ol{p}_r}^{q - 2}
                &\text{ in } \Omega, \\
            \grad \beta \cdot \n = 0
                &\text{ on } \Gamma,
        \end{cases}
    \end{align*}

    Letting $\ol{\vv} := \uu(t_1) - \uu(t_2)$, $\ol{\UU} := \UU(t_1) - \UU(t_2)$ we see from \cref{e:prElliptic} that
    \begin{align*}
	\begin{cases}
		\Delta \ol{p}_r = - \dv (\uu(t_1) \cdot \grad \ol{\vv} + \ol{\vv} \cdot \grad \uu(t_2))
			&\text{in } \Omega, \\
		\grad \ol{p}_r \cdot \n = -\prt_t \ol{U}^n
            - (\uu(t_1) \cdot \grad \ol{\vv} + \ol{\vv} \cdot \grad \uu(t_2)) \cdot \n \\
            \qquad\qquad\qquad
            - \CharFunc_{\Gamma_+} \dv_\Gamma
                (U^n(t_1) \ol{\vv}^\BoldTau + \ol{U}^n \cdot \uu^\BoldTau(t_2)
                ) \\
            \qquad\qquad\qquad
            +\CharFunc_{\Gamma_+} \dv_\Gamma
                (U^n(t_1) \ol{\UU}^\BoldTau + \ol{U}^n \cdot \UU^\BoldTau(t_2)
                )
			&\text{on } \Gamma,
	\end{cases}
    \end{align*}

     Then, using \cref{e:prElliptic},
    \begin{align*}
        \begin{split}
        \norm{\ol{p}_r}_{L^q(\Omega)}^q
            &= (\Delta \beta, \ol{p}_r)
            = - (\grad \beta, \grad \ol{p}_r) + \int_\Gamma (\grad \beta \cdot \n) \ol{p}_r
            = (\Delta \ol{p}_r, \beta) - \int_\Gamma (\grad \ol{p}_r \cdot \n) \beta \\
             &= - (\dv (\uu(t_1) \cdot \grad \ol{\vv} + \ol{\vv} \cdot \grad \uu(t_2)), \beta)
                    - \int_\Gamma (\grad \ol{p}_r \cdot \n) \beta \\
            &= (\uu(t_1) \cdot \grad \ol{\vv} + \ol{\vv} \cdot \grad \uu(t_2), \grad \beta)
            - \int_\Gamma ((\uu(t_1) \cdot \grad \ol{\vv} + \ol{\vv} \cdot \grad \uu(t_2)) \cdot \n) \beta \\
            &\qquad
                - \int_\Gamma (\grad \ol{p}_r \cdot \n) \beta \\
            &= (\uu(t_1) \cdot \grad \ol{\vv} + \ol{\vv} \cdot \grad \uu(t_2), \grad \beta)
                + \int_{\Gamma_+}
                    \Spanpr{
                        - \prt_t \ol{U}^n \\
            &\qquad
                    - \dv_\Gamma (U^n(t_1) \ol{\vv}^\BoldTau + \ol{U}^n \cdot \uu^\BoldTau(t_2)
                    + \dv_\Gamma (U^n(t_1) \ol{\UU}^\BoldTau + \ol{U}^n \cdot \UU^\BoldTau(t_2)
                    } \beta.
        \end{split}
    \end{align*}
    
    But, for vector fields $\vv$, $\ww$ with $\dv \vv = 0$,
    \begin{align*}
        (\vv \cdot &\grad \ww, \grad \beta)
            = \int_\Omega v^i \prt_i w^j \prt_j \beta
            = \int_\Omega v^i \prt_i (w^j \prt_j \beta)
                - \int_\Omega v^i w^j \prt_i \prt_j \beta \\
            &= (\vv, \grad (\ww \cdot \grad \beta))
                    - (\vv \otimes \ww, \grad \grad \beta)
            = - (\vv \otimes \ww, \grad \grad \beta)
                - \int_\Gamma v^\n (\ww \cdot \grad \beta),
    \end{align*}
    so
    \begin{align*}
        (\uu&(t_1) \cdot \grad \ol{\vv} + \ol{\vv} \cdot \grad \uu(t_2), \grad \beta) \\
            &= - (\uu(t_1) \otimes \ol{\vv} + \ol{\vv} \otimes \uu(t_2), \grad \grad \beta)
                - \int_\Gamma
                \pr{
                    U^n(t_1) (\ol{\vv} \cdot \grad \beta) + \ol{U}^n (\uu(t_2) \cdot \grad \beta)
                }.
    \end{align*}
    Also, applying \cref{L:IBPOnBoundary},
    \begin{align*}
        - \int_{\Gamma_+} \dv_\Gamma(&U^n(t_1) \ol{\vv}^\BoldTau
                    + \ol{U}^n \cdot \uu^\BoldTau(t_2)) \beta
            = \int_{\Gamma_+} (U^n(t_1) \ol{\vv}^\BoldTau
                    + \ol{U}^n \cdot \uu^\BoldTau(t_2)) \cdot \grad_\Gamma \beta \\
            &
            = \int_{\Gamma_+} (U^n(t_1) \ol{\vv}
                    + \ol{U}^n \cdot \uu(t_2)) \cdot \grad \beta.
    \end{align*}

    Hence,
    \begin{align*}
        \begin{split}
        \norm{\ol{p}_r}_{L^q(\Omega)}^q
            &= - (\uu(t_1) \otimes \ol{\vv} + \ol{\vv} \otimes \uu(t_2),
            \grad \grad \beta)
                - \int_\Gamma
                \pr{
                    U^n(t_1) (\ol{\vv} \cdot \grad \beta) + \ol{U}^n (\uu(t_2) \cdot \grad \beta)
                } \\
            &\qquad
                - \int_\Gamma \prt_t \ol{U}^n \beta
                + \int_{\Gamma_+} (U^n(t_1) \ol{\vv}
                    + \ol{U}^n \cdot \uu(t_2)) \cdot \grad \beta \\
            &\qquad
                + \int_{\Gamma_+} \dv_\Gamma (U^n(t_1) \ol{\UU}^\BoldTau + \ol{U}^n \cdot \UU^\BoldTau(t_2))
                    \, \beta \\
            &= - (\uu(t_1) \otimes \ol{\vv} + \ol{\vv} \otimes \uu(t_2),
                    \grad \grad \beta)
                - \int_\Gamma \prt_t \ol{U}^n \beta \\
            &\qquad
                - \int_{\Gamma_-} (U^n(t_1) \ol{\vv}
                    + \ol{U}^n \cdot \uu(t_2)) \cdot \grad \beta \\
            &\qquad
                - \int_{\Gamma_+} (U^n(t_1) \ol{\UU}^\BoldTau
                    + \ol{U}^n \cdot \UU^\BoldTau(t_2))
                    \cdot \grad \beta
        \end{split}
    \end{align*}

    Thus,
	\begin{align*}
		\norm{\ol{p}_r}_{L^q(\Omega)}^q
			&\le 2 \norm{\uu}_{L^\iny(Q)} \norm{\ol{\vv}}_{L^q(\Omega)}               \norm{\beta}_{W^{2, q'}(\Omega)}
				+ \norm{\prt_t \ol{U}^\n}_{L^q(\Gamma)} \norm{\beta}_{L^{q'}(\Gamma)} \\
            &\qquad
                + \norm{U^n(t_1)}_{L^{q'}(\Gamma)}
                    \brac{\norm{\ol{\vv}}_{L^\iny(\Gamma)} +     \norm{\ol{\UU}}_{L^\iny(\Gamma)}}
                    \norm{\grad \beta}_{L^q(\Gamma)} \\
           &\qquad
                + \norm{\ol{\UU}}_{L^{q'}(\Gamma)}
                    \brac{\norm{\UU}_{L^\iny(\Gamma)} +              \norm{\uu(t_2)}_{L^\iny(\Gamma)}}
                    \norm{\grad \beta}_{L^q(\Gamma)}.
	\end{align*}
    But by \cref{L:LinfTimeDiff},
    \begin{align*}
        \norm{\ol{\vv}}_{L^q(\Omega)} 
            &\le \norm{\ol{\vv}}_{L^\iny(\Omega)}
            \le \norm{\uu}_{\dot{C}^{0, 1}_t(Q)} \abs{t_1 - t_2}
            \le \norm{\uu}_{S^{1, \al}} \abs{t_1 - t_2} \\
        \norm{\prt_t \ol{U}^\n}_{L^q(\Gamma)}
            &\le \norm{\prt_t \ol{U}^\n}_{L^\iny(\Gamma)}
            \le \norm{\prt_t U}_{\dot{C}^{0, 1}_t(Q)}
            \le \norm{\UU}_{S^{2, \al}} \abs{t_1 - t_2}, \\
        \norm{\ol{\UU}}_{L^q(\Gamma)}
            &\le \norm{\UU}_{L^\iny(\Gamma)}
            \le \norm{U}_{\dot{C}^{0, 1}_t(Q)}
            \le \norm{\UU}_{S^{1, \al}} \abs{t_1 - t_2},
    \end{align*}
    so, exploiting \cref{e:BoundarySobolev}, we have
	\begin{align*}
		\norm{\ol{p}_r}_{L^q(\Omega)}^q
			&\le C_2 \norm{\beta}_{W^{2, q'}(\Omega)} \abs{t_1 - t_2}
            \le C_2 \norm{\ol{p}_r}_{L^q(\Omega)}^{q - 1} \abs{t_1 - t_2},
    \end{align*}
	giving \cref{e:prt1t2LqBounds}.
\end{proof}

\begin{prop}\label{P:PressureBounds}
    Assume that the data has regularity $N$ and let $\Omega'$ be as in \cref{L:prLemma}. Let $\uu \in \uSolSpaceN$ and let $p_r$ be the unique mean-zero solution to \cref{e:pr} with $q > 3/(1 - \al)$. Then for any integer $k$, $0 \le k \le N$,
	\begin{align}\label{e:PressureBounds}
		\begin{split}
		\norm{\grad p_r(t_1) - \grad p_r(t_2)}_{W^{k + 1, q}(\Omega')}
			&\le C_3(k)^2 \abs{t_1 - t_2}^\al, \\
		\norm{\grad p_r(t_1) - \grad p_r(t_2)}_{C^{{k, \al}}(\Omega')}
			&\le C_3(k)^2 \abs{t_1 - t_2}^\al
		\end{split}
	\end{align}
	for all $t_1, t_2 \in [0, T]$, where
    \begin{align}\label{e:C3}
        C_3(k)
            &:= C \brac{\norm{\bomega}_{C^{k, \al}(Q)}
                + \norm{\UU}_{C^{k + 2, \al}(Q)}
                + \norm{\uu}_{L^\iny(0, T; H)}}.
    \end{align}
\end{prop}
\begin{proof}
	We first prove \cref{e:PressureBounds}$_1$.
	Defining $\ol{p}_r := p_r(t_1) - p_r(t_2)$ and applying \cref{L:prLemma}, we have
	\begin{align*}
		\norm{\ol{p}_r}_{W^{k + 2, q}(\Omega')}
			\le C \brac{
				\norm{\Delta \ol{p}_r}_{W^{k, q}(\Omega)}
					+ \norm{\grad \ol{p}_r \cdot \n}_{W^{k + 1 - \frac{1}{q}, q}(\Gamma_+)}
					+ \norm{\ol{p}_r}_{L^q(\Omega)}
				}.
	\end{align*}
    This estimate is based upon the specific normalization of $\ol{p}_r$ given in \cref{L:prLemma}, but \cref{e:PressureBounds}$_1$ itself is independent of that normalization, since the gradient eliminates any normalization constant. (But see \cref{R:NormaliationgradpBound}.)
	
	Now,
	\begin{align*}
		\Delta \ol{p}_r
			&= \grad \uu(t_2) \cdot (\grad \uu(t_2))^T - \grad \uu(t_1) \cdot (\grad \uu(t_1))^T \\
			&= \grad (\uu(t_2) - \uu(t_1)) \cdot (\grad \uu(t_2))^T
				+ \grad \uu(t_1) \cdot (\grad (\uu(t_2) - \uu(t_1)))^T.
	\end{align*}
	Thus, for $k = 0$,
	\begin{align*}
		\begin{split}
		\norm{\Delta \ol{p}_r}_{L^q(\Omega)}
			&\le 2 \norm{\grad (\uu(t_1) - \uu(t_2))}_{L^q(\Omega)}
				\brac{\norm{\grad \uu(t_1)}_{L^\iny(\Omega)}
					+ \norm{\grad \uu(t_2)}_{L^\iny(\Omega)}}.
		\end{split}
	\end{align*}
	For $k \ge 1$, $W^{k, q}(\Omega)$ is an algebra, as noted above, so
	\begin{align*}
		\begin{split}
		\norm{\Delta \ol{p}_r}_{W^{k, q}(\Omega)}
			&\le C \norm{\grad (\uu(t_1) - \uu(t_2))}_{W^{k, q}(\Omega)}
				\brac{\norm{\grad \uu(t_1)}_{W^{k, q}(\Omega)}
					+ \norm{\grad \uu(t_2)}_{W^{k, q}(\Omega)}}.
		\end{split}
	\end{align*}
	In either case, we have
	\begin{align*}
		\begin{split}
		\norm{\Delta \ol{p}_r}_{W^{k, q}(\Omega)}
			&\le C \norm{\grad \uu}_{W^{k, \al}(\Omega)}
				\norm{\grad (\uu(t_1) - \uu(t_2))}_{W^{k, q}(\Omega)}.
		\end{split}
	\end{align*}

	But, setting $\bomega = \curl \uu$,
	\begin{align}\label{e:uDecomp}
		\uu(t_1) - \uu(t_2)
			&= K_{U^\n}[\bomega(t_1)] - K_{U^\n}[\bomega(t_2)]
			= K[\bomega(t_1) - \bomega(t_2)] + \ww,
	\end{align}
	where
	\begin{align*}
		\ww
			&= \VV(t_1) - \VV(t_2) + \uu_c(t_1) - \uu_c(t_2).
	\end{align*}
	Hence, applying \cref{L:BSLaw},
	\begin{align}\label{e:gradutqgadut2Diff}
		\begin{split}
		\norm{\grad \uu(t_1) - \grad \uu(t_2)}_{W^{k, q}(\Omega)}
			&\le C \norm{\bomega(t_1) - \bomega(t_2)}_{W^{k, q}(\Omega)}
					+ C \norm{\grad \ww}_{W^{k, q}(\Omega)}.
		\end{split}
	\end{align}
	Applying \cref{L:CNTimeDiffBound}, 
	\begin{align}\label{e:UseOfTrivialEmbedding}
		\norm{\bomega(t_1) - \bomega(t_2)}_{W^{k, q}(\Omega)}
			&\le \norm{\bomega(t_1) - \bomega(t_2)}_{C^k(\Omega)}
			\le \norm{\bomega}_{C^{k, \al}(Q)} \abs{t_1 - t_2}^\al.
	\end{align}
	Using \cref{L:CNTimeDiffBound} again,
	\begin{align*}
		\norm{\grad \ww}_{W^{k, q}(\Omega)}
			&\le C \norm{\grad \ww}_{C^k(\Omega)}
			\le \norm{\grad \ww}_{C^{k, \al}(Q)} \abs{t_1 - t_2}^\al \\
			&\le \norm{\UU}_{C^{k + 1, \al}(Q)} \abs{t_1 - t_2}^\al
				+ \norm{\uu}_{L^\iny(0, T; H)} \abs{t_1 - t_2}^\al.
	\end{align*}
	where we also used \cref{L:TamingHc}.
	Hence,
	\begin{align}\label{e:gradu1u2Bound}
		\begin{split}
		\norm{\grad \uu(t_1) - \grad \uu(t_2)}_{W^{k, q}(\Omega)}
			&\le C_3(k) \abs{t_1 - t_2}^\al.
		\end{split}
	\end{align}
 
    Similarly, $\norm{\grad \uu}_{W^{k, \al}(\Omega)} \le C_3(k)$, so
    \begin{align*}
        \norm{\Delta \ol{p}_r}_{L^q(\Omega)}
            &\le C_3(k)^2 \abs{t_1 - t_2}^\al.
    \end{align*}

    On $\Gamma_+$, we use \cref{e:NuNice}:
    \begin{align*}
        \grad p_r \cdot \n
            &= -\prt_t U_n + U^\n \dv_\Gamma \UU^\BoldTau
			+ (\kappa_1 + \kappa_2) (U^\n)^2
                - \uu^\BoldTau \cdot \grad_\Gamma U^\n
			+ \uu^\BoldTau \cdot \Cal{A} \uu^\BoldTau.
    \end{align*}
    Let us focus on the term $\uu^\BoldTau \cdot \Cal{A} \uu^\BoldTau$, the other terms being similarly, though more simply, bounded.

    By \cref{L:Sobolev}, ${W^{k + 1 - \frac{1}{q}, q}(\Gamma_+)}$ is an algebra. Hence, starting with the trace inequality,
    \begin{align*}
        &\norm{(\uu^\BoldTau \cdot \Cal{A} \uu^\BoldTau)(t_1)
                - \uu^\BoldTau \cdot \Cal{A} \uu^\BoldTau)(t_2)}
                    _{W^{k + 1 - \frac{1}{q}, q}(\Gamma_+)} \\
            &\qquad
            \le C \norm{\uu^\BoldTau(t_1) \cdot
                    \Cal{A} (\uu^\BoldTau(t_1) - \uu^\BoldTau(t_2))}
                        _{W^{k + 1 - \frac{1}{q}, q}(\Gamma_+)}
                    \\
            &\qquad\qquad
                +
            C \norm{(\uu^\BoldTau(t_1) - \uu^\BoldTau(t_2)) \cdot
                    \Cal{A} \uu^\BoldTau(t_2)}_{W^{k + 1 - \frac{1}{q}, q}(\Gamma_+)} \\
            &\qquad
            \le C \norm{\uu^\BoldTau(t_1)}_{W^{k + 1 - \frac{1}{q}, q}(\Gamma_+)}
                    \norm{\Cal{A} (\uu^\BoldTau(t_1) - \uu^\BoldTau(t_2))}
                        _{W^{k + 1 - \frac{1}{q}, q}(\Gamma_+)}
                    \\
            &\qquad\qquad
                +
            C \norm{(\uu^\BoldTau(t_1) - \uu^\BoldTau(t_2))}
                    _{W^{k + 1 - \frac{1}{q}, q}(\Gamma_+)}
                    \norm{\Cal{A} \uu^\BoldTau(t_2)}_{W^{k + 1 - \frac{1}{q}, q}(\Gamma_+)} \\
            &\qquad
            \le C \norm{\uu^\BoldTau(t_1)}_{W^{k + 1, q}(\Omega)}
                    \norm{\Cal{A} (\uu^\BoldTau(t_1) - \uu^\BoldTau(t_2))}
                        _{W^{k + 1, q}(\Omega)}
                    \\
            &\qquad\qquad
                +
            C \norm{(\uu^\BoldTau(t_1) - \uu^\BoldTau(t_2))}
                    _{W^{k + 1, q}(\Omega)}
                    \norm{\Cal{A} \uu^\BoldTau(t_2)}_{W^{k + 1, q}(\Omega)} \\
            &\qquad
            \le C \norm{\uu}_{W^{k + 1, q}(\Omega)}
                \norm{\uu(t_1) - \uu(t_2)}_{W^{k + 1, q}(\Omega)} \\
            &\qquad
            \le C_3(k) \norm{\uu}_{S^{k + 1, \al}} 
            \abs{t_1 - t_2}^\al.
    \end{align*}
    In the last inequality we used that $\norm{\uu}_{W^{k + 1, q}(\Omega)}  \le \norm{\uu}_{L^\iny([0, T]; C^{k + 1, \al}(\Omega))} \le \norm{\uu}_{S^{k + 1, \al}}$.

    Along with similar bounds on the other terms coming from $\grad p \cdot \n$, and using \cref{e:prt1t2LqBounds} of \cref{P:prLqBounds}, noting that the constant $C_2$ can be absorbed into $C_3(k)$, these bounds give \cref{e:PressureBounds}$_1$.
    Then \cref{e:SobolevEmbedding} with \cref{e:PressureBounds}$_1$ gives \cref{e:PressureBounds}$_2$.
\end{proof}

\begin{remark}\label{R:NormaliationgradpBound}
    Suppose that, instead of normalizing $p_r$ so it has mean-zero, we were, for a fixed $t_1, t_2 \in [0, T]$, to normalize it so that $M_q(p_r(t_1) - p_r(t_2)) = 0$. Then we would have
    \begin{align*}
        \norm{p_r(t_1) - p_r(t_2)}_{C^{{k + 1, \al}}(\Omega')}
		&\le C_3(k)^2 \abs{t_1 - t_2}^\al.
    \end{align*}
\end{remark}

\begin{remark}\label{R:ReasonForSobolev}
    In the proof of \cref{P:PressureBounds} we used both the embedding of $W^{k + 1, q}(\Omega)$ in $C^{k, \al}(\Omega)$ of \cref{e:SobolevEmbedding} and, in \cref{e:UseOfTrivialEmbedding}, the simple embedding of $C^k(\Omega)$ in $W^{k, q}(\Omega)$ (using that the domain $\Omega$ is bounded). In each of these inequalities we lost, in a sense, information. It would seem, then, that a more direct estimate using the \Holder space analog of the elliptic estimates in \cref{L:prLemma} would be cleaner. Were we to do that, however, \cref{{e:gradutqgadut2Diff}} would become
    \begin{align*}
        \norm{\grad \uu(t_1) - \grad \uu(t_2)}_{C^{k, \al}(\Omega)}
            &\le C \norm{\bomega(t_1) - \bomega(t_2)}_{C^{k, \al}(\Omega)}
                    + C \norm{\grad \ww}_{C^{k, \al}(\Omega)},
    \end{align*}
    and there would be no way to obtain the needed factor of $\abs{t_1 - t_2}^\al$ in \cref{e:PressureBounds}$_2$ as we obtained in \cref{e:UseOfTrivialEmbedding}.
\end{remark}

To account for time derivatives $\prt_t^j p_r$, $j \le N$, we note that \cref{e:pr} becomes
\begin{align*}
	\begin{cases}
		\Delta \prt_t^j p_r = - \prt_t^j (\grad \uu \cdot (\grad \uu)^T)
			&\text{in } \Omega, \\
		\grad \prt_t^j p_r \cdot \n = -\prt_t^{j + 1} U^\n - \prt_t^j  N[\uu]
			&\text{on } \Gamma,
	\end{cases}
\end{align*}
and the same analysis in \cref{P:PressureBounds,P:prLqBounds} applies to $\prt_t^j p_r$. This yields the following corollaries:

\begin{cor}\label{C:prtjOfpr}
     Let $\uu$ and $p_r$ be as in \cref{P:PressureBounds}
    and let $j$ be an integer with $0 \le j \le N$. 
    Fixing $t_1, t_2 \in [0, T]$, normalize $p_r$ so that
    $M_q(\prt_t^j p_r(t_1) - \prt_t^j p_r(t_2)) = 0$. Then
	\begin{align*}
		\begin{split}
		\norm{\prt_t^j p_r(t_1) - \prt_t^j p_r(t_2)}_{L^q(\Omega)}
			&\le C_2
				\abs{t_1 - t_2}.
		\end{split}
	\end{align*} 
\end{cor}
\begin{proof}
    The proof is the same as that of \cref{e:prt1t2LqBounds}, as we note that the constant $C_2$ already accounts for the presence of $\prt_t^j$.
\end{proof}

\begin{cor}\label{C:PressureBounds}
    Let $\uu$ and $p_r$ be as in \cref{P:PressureBounds}. Then $\grad_\Gamma p_r \in C^{N, \al}([0, T] \times \Gamma_+)$ with, recalling \cref{D:cc},
    \begin{align}\label{e:FinalpBound}
        \begin{split}
        \norm{\grad_\Gamma p_r}_{C^{N, \al}([0, T] \times \Gamma_+)}
            &\le  c_0 + c_X T^b,
        \end{split}
    \end{align}
    where $0 < b \le \al$.  
\end{cor}
\begin{proof}
    The observations above give that for any $j, k \ge 0$ integers with $0 \le j + k \le N$,
    \begin{align}\label{e:PressureBoundsjk}
        \begin{split}
        \norm{\prt_t^j \grad p_r(t_1) - \prt_t^j \grad p_r(t_2)}_{W^{k + 1, q}(\Omega')}
            &\le C_3(N)^2 \abs{t_1 - t_2}^\al, \\
        \norm{\prt_t^j \grad p_r(t_1) - \prt_t^j \grad p_r(t_2)}_{C^{{k, \al}}(\Omega')}
            &\le C_3(N)^2 \abs{t_1 - t_2}^\al
        \end{split}
    \end{align}
    for all $t_1, t_2 \in [0, T]$, where $C_3$ is as in \cref{e:C3}, and we note that
        \begin{align*}
        C_3(N)
            &\le C \brac{\norm{\bomega}_{C^{N, \al}(Q)}
                + \norm{\UU}_{C^{N, \al}(Q)}
                + \norm{\uu}_\uSolSpaceN}
            \le c_X,
    \end{align*}
    using that $\norm{\prt_t^N \uu}_{L^\iny(0, T; H)} \le C  \norm{\uu}_\uSolSpaceN$.

    Letting $\beta$ be any time-space multi-index with $\abs{\beta} \le N$, it follows from \cref{e:PressureBoundsjk}$_2$ that
    \begin{align*}
        \norm{D^\beta \grad p_r(t_1) - D^\beta \grad p_r(t_2)}_{L^\iny(\Omega')}
		\le \norm{D^\beta \grad p_r(t_1) - D^\beta \grad p_r(t_2)}_{C^{j, \al}(\Omega')}
		\le c_X \abs{t_1 - t_2}^\al,
    \end{align*}
    where $j = N - \abs{\beta}$. Then, letting $Q' = [0, T] \times \Omega'$,
    \begin{align*}
        \norm{D^\beta \grad p_r}_{L^\iny(Q')}
            &= \sup_{t \in [0, T]} \norm{D^\beta \grad p_r(t)}_{L^\iny(\Omega')} \\
            &\le \norm{D^\beta \grad p_r(0)}_{L^\iny(\Omega')}
                + \sup_{t \in [0, T]}
                    \norm{D^\beta \grad p_r(t) - D^\beta \grad p_r(0)}_{L^\iny(\Omega')} \\
            &\le \norm{D^\beta \grad p_r(0)}_{L^\iny(\Omega')}
                + c_X T^\al
            \le c_0 + c_X T^\al.
    \end{align*}
    Since $\grad p_r$ is continuous and this bound for all $\abs{\beta} \le N$, we see that
    \begin{align}\label{e:gradprInCN}
        \norm{\grad p_r}_{C^N(Q')} \le c_0 + c_X T^\al.
    \end{align}

    Now suppose that $\abs{\beta} = N$, so $j = N - \abs{\beta} = 0$. Then
    \begin{align}\label{e:gradprInCalx}
        \begin{split}
        \norm{D^\beta \grad p_r}_{\dot{C}^\al_x(Q')}
            &= \sup_{t \in [0, T]} \norm{D^\beta \grad p_r(t)}_{\dot{C}^\al(\Omega')} \\
            &\le \norm{D^\beta \grad p_r(0)}_{\dot{C}^\al(\Omega')}
                + \sup_{t \in [0, T]}
                    \norm{D^\beta \grad p_r(t) - D^\beta \grad p_r(0)}_{\dot{C}^\al(\Omega')} \\
            &\le \norm{D^\beta \grad p_r(0)}_{\dot{C}^\al(\Omega')}
                + c_X T^\al
            \le c_0 + c_X T^\al.
        \end{split}
    \end{align}
    This gives the spatial $C^\al$-regularity of the highest derivatives of $\grad p_r$.
    
    For the time regularity, let $f = D^\beta p_r$ and write \cref{e:PressureBoundsjk}$_2$, noting that $j = 0$, as
    \begin{align*}
		\norm{\grad (f(t_1) - f(t_2))}_{C^\al(\Omega')}
			&\le \norm{\grad (f(t_1) - f(t_2))}_{C^\al(\Omega)}
                \le c_X \abs{t_1 - t_2}^\al.
    \end{align*}
    
    Fix $t_1, t_2 \in [0, T]$ and normalize $f$ as in \cref{C:prtjOfpr}, so that $M_q(f(t_1) - f(t_2)) = 0$. Then by \cref{R:NormaliationgradpBound},
    \begin{align*}
		\norm{f(t_1) - f(t_2)}_{C^{1, \al}(\Omega')}
			&\le c_X \abs{t_1 - t_2}^\al,
    \end{align*}
    and \cref{C:prtjOfpr} gives
    \begin{align*}
        \norm{f(t_1) - f(t_2)}_{L^2(\Omega')}
            \le \norm{f(t_1) - f(t_2)}_{L^2(\Omega)}
            \le C^2 \abs{t_1 - t_2}.
    \end{align*}
    Then, applying the interpolation inequality in \cref{L:PreNirenberg} using \cref{C:prtjOfpr},
    \begin{align*}
    	\begin{split}
    	&\norm{\grad (f(t_1) - f(t_2))}_{C^0(\Gamma_+)} \\
            &\qquad
    		\le C \norm{f(t_1) - f(t_2)}_{C^{1, \al}(\Omega')}^a
    			\norm{f(t_1) - f(t_2)}_{L^2(\Omega')}^{1 - a} \\
    		&\qquad
    		\le C \brac{
    					c_X \abs{t_1 - t_2}^\al
    				}^a 
    			\brac{
    					c_X \abs{t_1 - t_2}
                    }^{1 - a}
    		\le c_X \abs{t_2 - t_1}^{\al'},
    	\end{split}
    \end{align*}
    where $\al < \al' := 1 - a(1 - \al) < 1$ ($a = 2\al/(3 + 2 \al)$ from \cref{L:PreNirenberg}).

    Although this bound was derived using a $(t_1, t_2)$-dependent normalization of the pressure, the bound itself is independent of that normalization and so applies uniformly for all $t_1, t_2 \in [0, T]$.
    It follows that
    \begin{align*}
        \norm{\grad f}_{\dot{C}^\al_t(Q')}
            &= \sup_{\substack{t_1 \ne t_2\\ \x \in \Omega'}}
                \frac{\abs{\grad f(t_1, \x) - \grad f(t_2, \x)}}
                    {\abs{t_1 - t_2}^\al}
            \le c_X \abs{t_1 - t_2}^{\al' - \al}.
    \end{align*}
    
    We conclude that $\norm{D^\beta \grad p_r}_{\dot{C}^\al_t(Q')} \le c_X T^{\al' - \al}$. Combined with \cref{e:gradprInCN,e:gradprInCalx}, \cref{L:HolderEquiv} gives \cref{e:FinalpBound} with $b = \max \set{\al, \al' - \al} > 0$.
\end{proof}

In \cref{P:PressureDiffTwoSolutions}, we obtain estimates on the difference of the pressure gradients for two velocity fields. These estimates will be used in the proof of \cref{P:AContinuous}, which only requires bounding the difference of pressures in $L^\iny([0, T] \times \Gamma_+)$. Hence, we produce the bound in the weakest feasible space, $L^\iny([0, T]; C^\al(\Omega))$.

\begin{prop}\label{P:PressureDiffTwoSolutions}
    Let $\uu_1, \uu_2 \in S^{1, \al}$, where $p_{r, 1}, p_{r, 2}$ solve \cref{e:pr} for $\uu_1$, $\uu_2$, respectively. Then
    \begin{align*}
        \norm{\grad p_{r, 1} - \grad p_{r, 2}}_{L^\iny([0, T]; C^\al(\Omega))}
            &\le C_4,
    \end{align*}
    where
    \begin{align*}
        C_4
            &:= C \sum_{j = 1}^2 \norm{\grad \uu_j}_{L^\iny(Q)}
                \brac{\norm{\uu_1 - \uu_2}_{L^\iny(\Omega)}
                    + \norm{\curl (\uu_1 - \uu_2)}_{L^q(\Omega)}
                    }.
    \end{align*}
\end{prop}  
\begin{proof}
    We parallel the proof of \cref{P:PressureBounds} for the case $k = 0$.

Letting
\begin{align*}
    \ol{\uu} := \uu_1 - \uu_2, \quad P := p_{1, r} - p_{2, r},
\end{align*}
noting that $\ol{\uu} \in H$, the elliptic problem for $P$ can be written,
\begin{align}\label{e:PElliptic}
	\begin{cases}
		\Delta P = - \dv (\uu_1 \cdot \grad \ol{\uu} + \ol{\uu} \cdot \grad \uu_2)
			&\text{in } \Omega, \\
		\grad P \cdot \n = 
            - (\uu_1 \cdot \grad \ol{\uu} + \ol{\uu} \cdot \grad \uu_2) \cdot \n
			- \CharFunc_{\Gamma_+} \dv_\Gamma (U^\n \ol{\uu}^\BoldTau)
			&\text{on } \Gamma.
	\end{cases}
\end{align}
Also, on $\Gamma_+$, from \cref{e:NuNice},
\begin{align}\label{e:POnGammaPlus}
	\begin{split}
		\grad P \cdot \n
			= - \ol{\uu}^\BoldTau \cdot \grad_\Gamma U^\n
			+ (\uu_1^\BoldTau \cdot \Cal{A} \ol{\uu}^\BoldTau
                + \ol{\uu}^\BoldTau \cdot \Cal{A} \uu_2^\BoldTau).
	\end{split}
\end{align}
    We can also write \cref{e:PElliptic} as
    \begin{align}
        \begin{cases}
            \Delta P
                = - (\grad \uu_1)^T \cdot \grad \ol{\uu}
                        + (\grad \ol{\uu})^T \cdot \grad \uu_2)
                &\text{in } \ol{Q}, \\
            \grad P \cdot \n
                = - (\uu_1 \cdot \grad \ol{\uu} + \ol{\uu} \cdot \grad \uu_2)
                &\text{on } [0, T] \times (\Gamma_- \cup \Gamma_0), \\
            \grad P \cdot \n
                = - \ol{\uu}^\BoldTau \cdot \grad_\Gamma U^\n
			+ \uu_1^\BoldTau \cdot \Cal{A} \ol{\uu}^\BoldTau
                + \ol{\uu}^\BoldTau \cdot \Cal{A} \uu_2^\BoldTau
                &\text{on } [0, T] \times \Gamma_+,
        \end{cases}
    \end{align}
    where we used \cref{e:NuNice}. Then applying \cref{L:prLemma} with $q > 3/(1 - \al)$, we have
    \begin{align*}
        \norm{P}_{W^{2, q}(\Omega')}
            \le C \brac{
                \norm{\Delta P}_{L^q(\Omega)}
                    + \norm{\grad P \cdot \n}_{W^{1 - \frac{1}{q}, q}(\Gamma_+)}
                    + \norm{P}_{L^q(\Omega)}
                }.
    \end{align*}
       
    Because $\ol{\uu} \in H$,
    \begin{align*}
        \begin{split}
        \norm{\Delta P}_{L^q(\Omega)}
            &\le 2
                \sum_{j = 1}^2 \norm{\grad \uu_j}_{L^\iny(\Omega)}
                \norm{\grad \ol{\uu}}_{L^q(\Omega)} \\
            &\le C \sum_{j = 1}^2 \norm{\grad \uu_j}_{L^\iny(\Omega)}
                \brac{
                    \norm{\curl \ol{\uu}}_{L^q(\Omega)}
                    +  \norm{\grad P_{H_c} \ol{\uu}}_{L^q(\Omega)}
                } \\
            &\le C \sum_{j = 1}^2 \norm{\grad \uu_j}_{L^\iny(\Omega)}
                \brac{
                    \norm{\curl \ol{\uu}}_{L^q(\Omega)}
                    +  \norm{\ol{\uu}}_H
                },
        \end{split}
    \end{align*}
    where we used \cref{L:TamingHc}, and
    \begin{align*}
        &\norm{\grad P \cdot \n}_{W^{1 - \frac{1}{q}, q}(\Gamma_+)}
            \le \norm{\grad P \cdot \n}_{L^\iny(\Gamma_+)} \\
            &\qquad
            \le \norm{\grad_\Gamma U^\n}_{L^\iny(\Gamma_+)}
                        \norm{\ol{\uu}^\BoldTau}_{L^\iny(\Gamma_+)}
                + C \sum_{j = 1}^2 \norm{\uu_j}_{L^\iny(\Gamma_+)}
			 \norm{\ol{\uu}^\BoldTau}_{L^\iny(\Gamma_+)} \\
            &\qquad
            \le \norm{\UU}_{C^1(Q)}
                        \norm{\ol{\uu}}_{L^\iny(\Omega)}
                + C \sum_{j = 1}^2 \norm{\uu_j}_{L^\iny(\Omega)}
			 \norm{\ol{\uu}}_{L^\iny(\Omega)}.
    \end{align*}

    It remains to bound $\norm{P}_{L^q(\Omega)}$. We follow the proof of \cref{e:prLqBound} of \cref{P:prLqBounds}, letting
    $\beta$ solve
	\begin{align*}
		\begin{cases}
			\Delta \beta = P \abs{P}^{q - 2}
				&\text{ in } \Omega, \\
			\grad \beta \cdot \n = 0
				&\text{ on } \Gamma,
		\end{cases}
	\end{align*}
    where $P$ is normalized so that $M_q(P) = 0$.
    We find that
	\begin{align}\label{e:prnormInLrt1t2Alt}
		\begin{split}
			\norm{P}_{L^q(\Omega)}^q
				&= - (\uu_1 \otimes \uu_1 - \uu_2 \otimes \uu_2, \grad \grad \beta)
					- \int_\Gamma U^\n \ol{\uu} \cdot \grad \beta.
		\end{split}	
	\end{align}

    For the first term on the right-hand side of \cref{e:prnormInLrt1t2Alt}, we use that
	\begin{align*}
		\norm{\uu_1 \otimes \uu_1 - \uu_2 \otimes \uu_2}_{L^q(\Omega)}
			&\le \sum_{j = 1}^2 \norm{\uu_j}_{L^\iny(\Omega)}
				\norm{\ol{\uu}}_{L^{q'}(\Omega)},
	\end{align*}
    so
    \begin{align*}
        - (\uu_1 \otimes \uu_1 - \uu_2 \otimes \uu_2, \grad \grad \beta)
            &\le C \sum_{j = 1}^2 \norm{\uu_j}_{L^\iny(\Omega)}
				\norm{\ol{\uu}}_{L^{q'}(\Omega)}
                    \norm{\beta}_{W^{2, q}(\Omega)}.
    \end{align*}
	
	For the boundary integral in \cref{e:prnormInLrt1t2Alt}, we have
	\begin{align*}
		&- \int_\Gamma U^\n \ol{\uu} \cdot \grad \beta
			\le C \norm{\UU}_{L^q(Q)} \norm{\ol{\uu}}_{L^\iny(\Omega)}
			\norm{\beta}_{W^{2, q'}(\Omega)}
			\le C \norm{\UU}_{S^{1, \al}} \norm{\ol{\uu}}_{L^\iny(\Omega)}
			\norm{\beta}_{W^{2, q'}(\Omega)}.
	\end{align*}
 
    Combining these bounds, we have
    \begin{align*}
        \norm{P}_{L^q(\Omega)}^q
            &\le C_4 \norm{\beta}_{W^{2, q'}(\Omega)}
            \le C_4 \norm{P}_{L^q(\Omega)}^{q - 1}.
    \end{align*}
    Since this bound holds uniformly over time, we can use \cref{e:SobolevEmbedding} to conclude  $\norm{P}_{L^\iny([0, T]; C^\al(\Omega))} \le C_4$, completing the proof.
\end{proof}

%
%
\Part{Part III: Estimates on the Operator $A$}

\subsection*{Organization of Part III}

In \cref{S:InvariantSet} we give the proof of \cref{P:InvariantSet} by first obtaining sufficient estimates on the operator $A$ using (primarily) the pressure estimates from \cref{S:Pressure} along with the estimates on the flow map from \cref{S:FlowMap}. In \cref{S:AContinuous}, we use these estimates on $A$ to prove \cref{P:AContinuous}. In \cref{S:FullBCs}, we give the proof of \cref{P:WeakVorticityIsVelocity}. Finally, in \cref{S:VorticityBCs}, we prove \cref{T:VorticityBCs}.

%
%
\section{An invariant set}\label{S:InvariantSet}

\noindent
We now make a series of estimates leading, in \cref{P:InvariantSet}, to the existence of an invariant set for the operator $A$. Recall that $\Dom_N(A)$ is given in \cref{D:DomNA}.

\begin{prop}\label{P:HBounds}
    Assume that $\uu \in \Dom_N(A)$. Then
    \begin{align}\label{e:HBound}
        \norm{\H}_{C^{N, \al}([0, T] \times \Gamma_+)}
                &\le c_0 + c_X T^b,
    \end{align}
    where $0 < b \le \al$ and $\H$ is given by \cref{e:H}.
    Suppose that $\uu_1$, $\uu_2$ both lie in $\Dom_0(A)$ with
    $\norm{\uu_1}_{S^{1, \al}}, \norm{\uu_1}_{S^{1, \al}} \le M$,
    and let $\H_i$, $i = 1, 2$, be given by \cref{e:H} with $\uu = \uu_i$.
    Then, letting $q > 3/(1 - \al)$,
    \begin{align}\label{e:HBoundDiff}
        \norm{\H_1 - \H_2}_{L^\iny([0, T] \times \Gamma_+)}
                &\le c_0 M
                \brac{\norm{\uu_1 - \uu_2}_{L^\iny(\Omega)}
                    + \norm{\curl (\uu_1 - \uu_2)}_{L^q(\Omega)}
                    }.
    \end{align}
\end{prop}
\begin{proof}
    Let $X := C^{N, \al}([0, T] \times \Gamma_+)$. From \cref{e:H}, we can write,
	\begin{align*}
		\H^\BoldTau
			&= \delta_1 + \delta_2
				- \grad_\Gamma p_r, \quad
		H^\n
			= \curl_\Gamma \UU^\BoldTau,
	\end{align*}
	where
	\begin{align*}
		\delta_1
			&:= \frac{1}{U^\n}
				\brac{
					-\prt_t \UU^\BoldTau - \grad_\Gamma \pr{\frac{1}{2} \abs{\UU}^2}
				 	 	+ \f}^\perp, \quad
		\delta_2
			:= \frac{1}{U^\n}
					\curl_\Gamma \UU^\BoldTau \uu^\BoldTau.
	\end{align*}
    Hence, $\norm{H^\n}_X \le \norm{\UU}_{C^{N + 1, \al}(Q)} + \norm{\f}_{C^{N, \al}(Q)} \le c_0$, $\norm{\delta_1} \le c_0$, and $\norm{\grad_\Gamma p_r}_X \le c_0 + c_X T^b$ by \cref{C:PressureBounds}.
    Then, letting $\phi = (U^\n)^{-1} \curl_\Gamma \UU^\BoldTau$ and applying \cref{P:SConsequences},
    \begin{align*}
        \norm{\delta_2}_X
            &= \norm{\phi \uu^\BoldTau}_X
            \le \norm{\phi}_X \norm{\uu}_X
            \le c_0 + (c_0 + M) \max \set{T^{1 - \al}, T^\al}.
    \end{align*}
    Together, these bounds yield \cref{e:HBound}.

    Now suppose that $\uu_1$, $\uu_2$ both lie in $\Dom_0(A)$; hence, they then have the same initial data, and the same $\UU$ and $\f$. So reviewing the estimates that led to \cref{e:HBound}, we see that $H_1^n - H_2^n = 0$ and many terms in $\H^\BoldTau$ cancel, leaving
    \begin{align*}
	\H_1^\BoldTau - \H_2^\BoldTau
		&:= \frac{1}{U^\n}
				\brac
                    {\grad_\Gamma^\perp (p_{r, 1} - p_{r, 2})
				+ \curl_\Gamma \UU^\BoldTau
                    (\uu_1^\BoldTau - \uu_2^\BoldTau)
                    },
    \end{align*}
    where $p_{r, j}$ is the pressure corresponding to $u_j$. The bound in \cref{e:HBoundDiff}
    follows from this observation and \cref{P:PressureDiffTwoSolutions}.
\end{proof}    

\begin{prop}\label{P:Koch22Sub}
	Assume that $\uu \in \Dom_N(A)$. With $\Lambda$ as in \cref{e:Lambda},
	\begin{align*}
		\norm{\Lambda \uu}_{C^{N, \al}(Q)}
                &\le \cN(c_0, M, T),
	\end{align*}
	where $\cN$ is as in \cref{D:cc} and $M$ bounds $\norm{\uu}_\uSolSpaceN$,
	as in \cref{e:M}.
\end{prop}
\begin{proof}
	First assume no forcing.
	Let $\bomega_0 = \bomega(0)$
	and recall the definition of $\pos_0$ in \cref{e:pos0}.
	From \cref{e:LagrangianForm}, we can write,
	$\ol{\bomega} := \Lambda \uu = \ol{\bomega}_\pm$ on $U_\pm$,
	where
	\begin{align}\label{e:omegaplusminus}
		\begin{split}
		\ol{\bomega}_-(t, \x)
			&= \grad \eta(0, t; \pos_0) \bomega_0(\pos_0)
				\text{ on } U_-, \\
		\ol{\bomega}_+(t, \x)
			&= \grad \eta(\time(t, \x), t; \pos(t, \x)) \H(\time(t, \x), \pos(t, \x))
				\text{ on } U_+.
		\end{split}
	\end{align}
	
	It follows, using \cref{L:QuickFlow}, \cref{R:cNFlowBounds}, \cref{P:HBounds} that
	\begin{align*}
		\norm{\ol{\bomega}_-(t, \x)}_{L^\iny(U_-)}
			&\le \norm{\grad \eta}_{L^\iny(Q)} \norm{\bomega_0}_{L^\iny(\Omega)}
			\le \norm{\bomega_0}_{L^\iny(\Omega)} e^{M T}, \\
		\norm{\ol{\bomega}_+(t, \x)}_{L^\iny(U_+)}
			&\le \norm{\grad \eta}_{L^\iny(Q)} \norm{\H}_{L^\iny([0, T] \times \Gamma_+)}
			\le [\norm{\bomega_0}_{L^\iny(\Gamma_+)} + M T^\al] e^{M T},
	\end{align*}
	which shows that $\norm{\Lambda \uu}_{L^\iny(Q)} \le \cN$.
	
	Let us now first treat the case $N = 0$, to get a better understanding
	of the estimates involved.
	Using \cref{L:QuickFlow,R:cNFlowBounds} along with \cref{L:HolderProd,L:HolderComp},
	we see that
	\begin{align*}
		\norm{\ol{\bomega}_-}_{C^\al(U_-)}
			&\le \norm{\grad \eta(0, t; \pos_0)}_{C^\al(U_-)}
					\norm{\bomega_0(\pos_0)}_{C^\al(U_-)}
						\\
			&\le \norm{\grad \eta(0, t; \cdot)}_{C^\al(Q)}
						[\norm{\grad \pos_0}_{L^\iny(U_-)}^\al]^2
					\norm{\bomega_0}_{C^\al(\Omega)} \\
			&\le \norm{\bomega_0}_{C^\al(\Omega)} [1 + M e^{(1 + 2 \al) M T} T^{1 - \al}]
				e^{2 M T}
                \le \cN(c_0, M, T).
	\end{align*}
	Note that $\cN(c_0, M, 0) = \norm{\bomega_0}_{C^\al(\Omega)}$, giving \cref{e:cN0}.
		
	Similarly,
	\begin{align*}
		\norm{\ol{\bomega}_+(t, \x)}_{C^\al(U_+)}
			&\le \norm{\grad \eta(\time(t, \x), t; \pos(t, \x))}_{C^\al(U_+)}
					\norm{\H(\time(t, \x), \pos(t, \x))}_{C^\al(U_+)}.
	\end{align*}
	Using \cref{L:QuickFlow,L:HolderComp},
	\begin{align*}
		&\norm{\grad \eta(\time(t, \x), t; \pos(t, \x))}_{C^\al(U_+)}
			\le \norm{\grad \eta(t_1, t_2; \x)}_{C^\al([0, T]^2 \times \Omega)}
				[1 + \norm{D \mu}_{L^\iny(U_+)}]^\al \\
			&\qquad
			\le [e^{MT} + e^{(1 + 2 \al) MT} M T^{1 - a}]
				[(1 + M^2) e^{M T}]^\al.
	\end{align*}
	Then, using \cref{L:QuickFlow} with \cref{R:cNFlowBounds}, and \cref{P:HBounds},
	\begin{align*}
		&\norm{\H(\time(t, \x), \pos(t, \x))}_{C^\al(U_+)}
			\le \norm{\H}_{C^\al([0, T] \times \Gamma_+)} [1 + \norm{D \mu}_{L^\iny(U_+)}]^\al \\
			&\qquad
			\le [(1 + M^2) e^{MT}]^\al
                [\cN + C M(M + 1)].
	\end{align*}
	These bounds lead to
	\begin{align*}
		&\norm{\ol{\bomega}_+(t, \x)}_{C^\al(U_+)}
			\le \cN
                e^{(1 + 4 \al)MT}
                (1 + M^2)^{1 + \al}
                (1 + M T^{1 - \al})
                \le \cN(c_0, M, T),
	\end{align*}
	and we can see also that \cref{e:cN0} holds.
				
    But we know from \cref{T:LinearExistence} that $\ol{\bomega} \in \vortSpace$, because we assumed $\cond_0$: hence, taking the maximum of the bounds for $\ol{\bomega}_\pm$ on $U_\pm$, and using that $\norm{\bomega_0}_{L^\iny(\Omega)} \le M$ yields the bound on  $\norm{\Lambda \uu}_{C^\al(Q)}$.
	
	Now consider $N \ge 1$.
	The expressions for $\ol{\bomega}_\pm$ in \cref{e:omegaplusminus} each
	consist of two factors.
	We first apply Leibniz's product rule to these expressions then apply
	the chain rule to each term. For $\ol{\bomega}_+$, if $\beta$
	is a time-space multi-index with $\abs{\beta} = N$, then $D^\beta \ol{\bomega}_+$
	consists of a finite sum of terms of the form,
	\begin{align*}
		D^{\beta_1}
			\grad \eta(\time(t, \x), t; \pos(t, \x))
		D^{\beta_2} \H(\time(t, \x), \pos(t, \x))
		\prod_{\ell = 1}^n D^{\beta_3^\ell} \mu(t, \x)
				\text{ on } U_+,
	\end{align*}
	where $\beta_1 + \beta_2 = \beta$ and $\sum_{\ell = 1}^n \abs{\beta_3^\ell} = \abs{\beta}$.
	The factors can be controlled by
	\cref{P:HBounds,L:QuickFlow} with \cref{R:cNFlowBounds}.
	Following the similar process for $D^\beta \ol{\bomega}_-$
	leads to an estimate for
	$\norm{\Lambda \uu}_{C^{N, \al}(Q)}$ of the same general form as for
	$\norm{\Lambda \uu}_{C^\al(Q)}$.

    This gives the bound $\norm{\Lambda \uu}_{C^{N, \al}(Q)} \le \cN(c_0, M, T)$ in the absence of forcing. For forcing, we must bound $\G_\pm$ of \cref{e:LagrangianForm}. Now, as noted in
\cref{R:Paper1},
$\ol{\bomega}_\pm$ need not be continuous across $\Cal{S}$ and $\G_\pm$ need not be continuous across $\Cal{S}$; rather, $\ol{\bomega}_\pm + \G_\pm$ is $C^{N, \al}$-continuous across $\Cal{S}$. Nonetheless, adding separate bounds on $\ol{\bomega}_+$, $\ol{\bomega}_-$, $\G_+$, and $\G_-$ in $C^{N, \al}(U_\pm)$ give the bound on $\Lambda \uu = \ol{\bomega}$ in $C^{N, \al}(Q)$.

    It follows from \cref{R:cNFlowBounds} and repeated applications of the chain rule that, in fact, 
    \begin{align*}
        \norm{\G_\pm}_{C^{N, \al}(Q))} \le \cN,
    \end{align*}
    which completes the proof.
\end{proof}

\begin{prop}\label{P:AuCNalBound}
    For any $\uu \in \Dom_N(A)$,
    \begin{align*}
        \norm{A \uu}_{C^{N, \al}(Q)}
            &\le \cN(c_0, M, T).
    \end{align*}
\end{prop}
\begin{proof}
    Let $\uu \in \Dom_N(A)$ and $\vv = A \uu$, which we note satisfies \cref{e:LinearVelEq}. Our goal is to bound $\vv$ in the space $R^{N + 1, \al}$ of \cref{D:Rkal} and then apply \cref{L:Rkal}.
    
    From \cref{e:LinearVelEq} we have $\prt_t \vv = \f - P_H[\uu \cdot \ol{\bOmega}]$. But,
    \begin{align*}     
        P_H \vv
            = P_H (\vv - \UU) + P_H \UU
            = \vv  - \ww, \quad
        \ww := \UU - P_H \UU,
    \end{align*}
    so
    \begin{align}\label{e:prtvvBase}
        \prt_t \vv = \f + \ww - P_H[\uu \cdot \ol{\bOmega}],
    \end{align}
    and, for any $1 \le j \le N + 1$,
    \begin{align*}
        \prt_t^j \vv = \prt_t^{j - 1} \f + \prt_t^{j - 1} \ww - P_H[\prt_t^{j - 1} (\uu \cdot \ol{\bOmega})].
    \end{align*}
    Hence, letting $X^j := L^\iny([0, T]; C^{N + 1 - j, \al}(\Omega))$,
    \begin{align*}
        &\norm{\prt_t^j \vv}_{L^\iny([0, T]; C^{N + 1 - j, \al}(\Omega))}
            = \norm{\prt_t^j \vv}_{X^j} \\
            &\qquad
            \le \norm{\prt_t^{j - 1} \f}_{X^j} + \norm{\prt_t^{j - 1} \ww}_{X^j}
                + C \norm{\prt_t^{j - 1} (\uu \cdot \ol{\bOmega})}_{X^j} \\
            &\qquad
            \le \norm{\f}_{C^{N, \al}(Q)} + \norm{\ww}_{C^{N, \al}(Q)}
                + C \norm{(\uu \cdot \ol{\bOmega})}_{C^{N, \al}(Q)} \\
            &\qquad
            \le \norm{\f}_{C^{N, \al}(Q)} + \norm{\ww}_{C^{N, \al}(Q)}
                + C \norm{\uu}_{C^{N, \al}(Q)} \norm{\ol{\bOmega}}_{C^{N, \al}(Q)} \\
            &\qquad
            \le c_0 + \cN
            \le \cN.
    \end{align*}
    We used the continuity of $P_H$ in the algebra $C^{N, \al}(\Omega)$ (though not in $C^{N, \al}(Q)$)
    and in the last inequality we used that $\norm{\uu}_{C^{N, \al}(Q)} \le \cN$ by \cref{P:SConsequences} and, because $\ol{\bomega} = \Lambda \uu$, $\norm{\ol{\bomega}}_{C^{N, \al}(Q)} \le \cN$ by \cref{P:Koch22Sub}. This same bound follows for $j = 0$ by applying \cref{L:BSLaw},
	which completes the demonstration that $\vv \in R^{N + 1, \al}$.

    We conclude by \cref{L:Rkal} that $\norm{A \uu}_{C^{N, \al}(Q)} = \norm{\vv}_{C^{N, \al}(Q)} \le \cN$.
\end{proof}   

\begin{cor}\label{C:AuSolspaceBound}
    For any $\uu \in \Dom_N(A)$, $\norm{A \uu}_\uSolSpaceN \le \cN(c_0, M, T)$.
\end{cor}
\begin{proof}
    In light of \cref{P:AuCNalBound,P:Koch22Sub}, it remains only to bound
    $\prt_t^{N + 1} A \uu$ in the space $L^\iny([0, T]; C^\al(\Omega))$.
    For this, we apply the Leray projector to \cref{e:LinearVelEq}. Writing
    $A \uu = \vv = \ol{\vv} + \VV$, where $\vv \in H$, we have
   	$
		\prt_t \ol{\vv}
			= - P_H (\uu \cdot \ol{\bOmega}) + \f,
	$
	using that $P_H \VV = 0$ because $\VV$ is a gradient.
	Hence,
   	\begin{align*}
		\prt_t^{N + 1} \ol{\vv}
			&= - \prt_t^N P_H (\uu \cdot \ol{\bOmega}) + \prt_t^N \f
			= - P_H \pr{\prt_t^N (\uu \cdot \ol{\bOmega})} + \prt_t^N \f.
	\end{align*}
	Using that $P_H$ is continuous in $C^\al(\Omega)$, we see that
	\begin{align*}
		\norm{\prt_t^{N + 1} \vv}_{L^\iny([0, T]; C^\al(\Omega))}
			&\le C \norm{\uu \cdot \ol{\bOmega}}_{C^{N, \al}(Q)}
				+ \norm{f}_{C^{N, \al}(Q)} \\
			&\le \norm{\uu}_{C^{N, \al}(Q_T)}
				\norm{\ol{\bomega}}_{C^{N, \al}(Q_T)} + \norm{f}_{C^{N, \al}(Q)} \\
			&\le \cN(c_0, M, T) \cN(c_0, M, T) + \cN(c_0, M, T)
			\le \cN(c_0, M, T),
	\end{align*}
	where we applied \cref{P:SConsequences,P:Koch22Sub}.
	This gives the required bound on $\prt_t^{N + 1} A \uu$ in $L^\iny([0, T]; C^\al(\Omega))$.
\end{proof}

Having established our many estimates, we can now give the proof of \cref{P:InvariantSet}.

\begin{proof}[\textbf{Proof of \cref{P:InvariantSet}}]
	For an arbitrary $T > 0$, recall that we set
	\begin{align*}
		\Cal{K}
			&= \Cal{K}_{M, T}
            := \set{\uu \in \Dom_N(A) \colon
                \norm{\uu}_\uSolSpaceN \le M}.
	\end{align*}

	By \cref{C:AuSolspaceBound}, for \textit{any} $T, M > 0$,
	\begin{align*}
		\uu \in \Cal{K}_{M, T}
			\implies \norm{A \uu}_\uSolSpaceN \le \cN(c_0, M, T).
	\end{align*}
	
	Since $\cN(c_0, M, T) = c_0$, which is independent of $M$,
	we can now choose a specific $M > c_0$. Then the continuity of
	$\cN(c_0, M, T)$ allows us to choose
	$T > 0$ for which $\cN(c_0, M, T) < M$. But this means that
	\begin{align*}
		\uu \in \Cal{K}_{M, T}
			\implies \norm{A \uu}_\uSolSpaceN \le M
			\implies A \uu \in \Cal{K}_{M, T}.
	\end{align*}
	That is, $\Cal{K}_{M, T}$ is invariant under the operator $A$.
\end{proof}

%
%
\section{Continuity of the operator $A$}\label{S:AContinuous}

\noindent 
Throughout this section, we let $M$, $T$, and $\Cal{K} = \Cal{K}_{M, T}$ be fixed, as given by \cref{P:InvariantSet}. We also fix $\beta \in (0, \al)$ arbitrarily. 

Before giving the proof of \cref{P:AContinuous}, we establish a series of estimates on the difference between two velocity fields in $\Cal{K}$ and the difference of their corresponding flow maps.
We assume that 
\begin{align*}
	\uu_1, \uu_2 \text{ are two vector fields in } \Cal{K},
\end{align*}
and define the following:
\begin{itemize}
	\item
		$\bomega_j := \curl \uu_j$ for $j = 1, 2$,
		
	\item
		$\eta_j$, $\time_j$, $\pos_j$, $U_\pm^j$, and hypersurface $\Cal{S}_j$
		are defined as in \cref{S:FlowMap} for $\uu_j$, $j = 1, 2$,
		
	\item
		$V_\pm := U_\pm^1 \cap U_\pm^2$,
		
	\item
		$W := Q \setminus (V_+ \cup V_-)$.
\end{itemize}

We define  $\mu_j \colon U_+ \to [0, T] \times \Gamma_+$ by
\begin{align*}
    \mu_j(t, \x) := (\time_j(t, \x), \pos_j(t, \x))
\end{align*}
and
\begin{align*}
	\ww := \uu_1 - \uu_2, \quad \mu := \mu_1 - \mu_2.
\end{align*}
We set
\begin{align}\label{e:thetabeta}
	\theta_\beta
		:= \norm{\ww}_{\WeakSArg{0}}
		= \norm{\ww}_{C^\beta(Q)} + \norm{\curl \ww}_{C^\beta(Q)}.
\end{align}

\begin{remark}
    In this section, we use the convention that $\FMOC$ stands for a continuous function from $[0, \iny)$ to $[0, \iny)$ with $\FMOC(0) = 0$. Its precise values will be unimportant, and may vary from occurrence to occurrence in expressions.  
 \end{remark}

\cref{L:WeakSBound} gives two interpolation inequalities between $L^\iny(Q)$-based spaces and $\WeakSTwoArgs{N + 1}{\al}$. When applied to $\ww = \uu_1 - \uu_2$, we have $\norm{\ww}_{\WeakSTwoArgs{N + 1}{\al}} \le \norm{\ww}_{\uSolSpaceN} \le 2 M$, so \cref{L:WeakSBound} will allow us to control the size of higher norms of $\ww$ and $\curl \ww$ by the size of their $L^\iny(Q)$ norms. This will greatly simplify our arguments, since estimating $\ww$ and $\curl \ww$ in $L^\iny$ norms is much easier than in higher norms, and because of this, all of the estimates we obtain following \cref{L:WeakSBound} will be in $L^\iny$.

\begin{lemma}\label{L:WeakSBound} 
    For any $\uu \in \WeakSTwoArgs{N + 1}{\al}$,
    \begin{align*}
        \norm{\curl \uu}_{C^{N, \beta}(Q)}
            &\le \FMOC \pr{\norm{\curl \uu}_{L^\iny(Q)}}, \\
        \norm{\uu}_{\WeakS}
            &\le \FMOC \pr{\norm{\uu}_{L^\iny(Q)} + \norm{\curl \uu}_{L^\iny(Q)}},
    \end{align*}
    where the function $\FMOC$ depends on $\norm{\uu}_{\WeakSTwoArgs{N + 1}{\al}} \le M$.
\end{lemma}
\begin{proof}
Let $f \in C^{N, \al}(Q)$. By \cref{L:Interpolation},
\begin{align}\label{e:KeyfBound}
	\norm{f}_{C^{N, \beta}(Q)}
		&\le \norm{f}_{L^\iny(Q)}
			+ F_c(\norm{f}_{C^{N, \al}(Q)})
			\norm{f}_{L^2(Q)}^{1 - a},
\end{align}
where $F_c(x) = x^{a_1} + x^{a_N} + x^{a'}$, where $a_1$ and $a_N$ are given in \cref{L:PreNirenberg}, and $a'$ is given in \cref{L:Interpolation}, and each of $a_1, a_N, a'$ lies in $(0, 1)$. The exponent $a$, which also lies in $(0, 1)$, depends upon whether $\norm{f}_{L^2(Q)}$ is greater or less than $1$.

Because $0 < a_1, a_N, a' < 1$ and for any $0 < b < 1$, $(x + y)^b < x^b + y^b \le 2^b (x + y)^b$, we see that $F_c(x + y) \le C (F_c(x) + F_c(y))$. Applying this inequality with \cref{e:KeyfBound} to $f = \uu$ and $f = \curl \uu$, and using that $F_c(x + y) \le C (F_c(x) + F_c(y))$ gives the result.
\end{proof}

By \cref{L:QuickFlow}, we have, for $j = 1, 2$,
\begin{align}\label{e:uletal}
	\begin{array}{ll}
		\norm{\eta_j(0, \cdot; \cdot)}_{C^{N + 1, \al}(Q} \le C(T, M).
	\end{array}
\end{align}
We generally do not state the dependence of constants on $T$ and $M$, which are fixed and hence  have no impact on the proof of \cref{P:AContinuous}. We do state such dependence explicitly, however, when it makes the nature of the bound being derived clearer.

\begin{lemma}\label{L:muDiff}
	We have,
	\begin{align*}
		\norm{\mu}_{L^\iny(V_+)}
			\le C(T, M) T \theta_\beta.
	\end{align*}
\end{lemma}
\begin{proof}
	We know from \LinearPaper{Lemma 3.5} that
	$\mu_j$ is transported by the flow map for $\uu_j$; that is,
	\begin{align*}
		&\prt_t \mu_1 + \uu_1 \cdot \grad \mu_1 = 0, \\
		&\prt_t \mu_2 + \uu_2 \cdot \grad \mu_2 = 0.
	\end{align*}
	Hence,
	\begin{align*}
		\prt_t \mu + \uu_1 \cdot \grad \mu
			= - \ww \cdot \grad \mu_2,
	\end{align*}
	or,
	\begin{align*}
		\diff{}{t} &\mu(t, \eta_1(0, t; \x))
			= - (\ww \cdot \grad \mu_2)(t, \eta_1(0, t; \x)).
	\end{align*}
	Integrating in time, using that
	$\mu(t, \eta_1(0, t; \x))|_{t = 0} = 0$,
	and employing \cref{L:QuickFlow} gives
	\begin{align*}
		\mu(t, \eta_1(0, t; \x))
			&= - \int_0^t (\ww \cdot \grad \mu_2)(s, \eta_1(0, s; \x))
			\le \norm{\ww}_{L^\iny(Q)} \norm{\grad \mu_2}_{L^\iny(Q)} \\
			&\le C(T, M) \theta_\beta.
			\qedhere
	\end{align*}
\end{proof}

\begin{lemma}\label{L:FlowDiff}
		We have
		\begin{align*}
			\norm{\eta_1  - \eta_2}_{L^\iny([0, T]^2 \times \Omega)}
				&\le C(T, M) T \theta_\beta, \\
			\norm{\grad \eta_1 - \grad \eta_2}_{L^\iny([0, T]^2 \times \Omega)}
				&\le C(T, M) T [\theta_\beta + \theta_\beta^{\al}].
		\end{align*}
\end{lemma}
\begin{proof}
	We have,
	\begin{align*}
		\eta_1(t_1, t_2; \x) - \eta_2(t_1, t_2; \x)
			&= \int_{t_1}^{t_2} \brac{\uu_1(s, \eta_1(t_1, s; \x))
				- \uu_2(s, \eta_2(t_1, s; \x))} \, ds.
	\end{align*}
	Fixing $t_1$, using \cref{e:uletal}, \cref{L:HolderComp}, \cref{L:HolderDiff}, and applying
	Minkowski's integral inequality gives
	\begin{align*}
		&\abs{\eta_1(t_1, t; \x) - \eta_2(t_1, t; \x)} \\
			&\qquad
			\le \int_{t_1}^t
					\abs{\uu_1(s, \eta_2(t_1, s; \x))
						- \uu_2(s, \eta_2(t_1, s; \x))} \, ds
				\\
			&\qquad\qquad
			+ \int_{t_1}^t
					\abs{\uu_1(s, \eta_1(t_1, s; \x))
						- \uu_1(s, \eta_2(t_1, s; \x))} \, ds
				\\
			&\qquad
			\le \int_{t_1}^t
					\norm{\uu_1(s) - \uu_2(s)}_{L^\iny(\Omega)} \, ds
				+ \int_{t_1}^t
					\norm{\uu_1(s)}_{\dot{C}^1(\Omega)}
							\abs{\eta_1(t_1, s; \x) - \eta_2(t_1, s; \x)} \, ds
				\\
			&\qquad
			\le T \theta_\beta
				+ C(T, M) \int_{t_1}^t \norm{\eta_1(t_1, s; \cdot)
						- \eta_2(t_1, s; \cdot)}_{L^\iny(\Omega)} \, ds.
	\end{align*}
	Taking the supremum over $\x$
	and applying \Gronwalls Lemma gives	\begin{align*}
		\norm{\eta_1(t_1, t; \x) - \eta_2(t_1, t; \x)}_{C_t([0, T]; L^\iny_{\x}(\Omega))} 
			&\le T e^{C(M, T) T} \theta_\beta.
	\end{align*}
	Since this holds uniformly for all $t_1 \in [0, T]$, we obtain the first bound.

	Similarly, starting from
	\begin{align*}
		\grad \eta_1&(t_1, t; \x) - \grad \eta_2(t_1, t; \x)
			= \int_{t_1}^t \brac{\grad_\x (\uu_1(s, \eta_1(t_1, s; \x)))
				- \grad_\x(\uu_2(s, \eta_2(t_1, s; \x)))} \, ds \\
			&= \int_{t_1}^t
					\brac{\grad \uu_1(s, \eta_1(t_1, s; \x)) \grad \eta_1(t_1, s; \x)
						- \grad \uu_2(s, \eta_2(t_1, s; \x))
							\grad \eta_2(t_1, s; \x)} \, ds,
	\end{align*}
	we find
	\begin{align*}
		&\abs{\grad \eta_1(t_1, t; \x) - \grad \eta_2(t_1, t; \x)} \\
			&\qquad
			\le \int_{t_1}^t
					\abs{\grad \uu_1(s, \eta_1(t_1, s; \x)) \grad \eta_1(t_1, s; \x)
						- \grad \uu_1(s, \eta_2(t_1, s; \x))
							\grad \eta_1(t_1, s; \x)} \, ds
				\\
			&\qquad\qquad
			+ \int_{t_1}^t
					\abs{(\grad \uu_1(s, \eta_2(t_1, s; \x))
						- \grad \uu_2(s, \eta_2(t_1, s; \x)))
						\grad \eta_2(t_1, s; \x)} \, ds
				\\
			&\qquad\qquad
			+ \int_{t_1}^t
					\abs{\grad \uu_1(s, \eta_2(t_1, s; \x))
							(\grad \eta_1(t_1, s; \x)
								-\grad \eta_2(t_1, s; \x))} \, ds
				\\
			&\qquad
			\le \int_{t_1}^t
					\norm{\uu_1(s)}_{\dot{C}^\al}
							\norm{\eta_1(t_1, s; \cdot) - \eta_2(t_1, s; \cdot)}_{L^\iny(\Omega)}^\al
							\norm{\grad \eta_1(t_1, s; \cdot)}_{L^\iny(\Omega)} \, ds \\
			&\qquad\qquad
			+ \int_{t_1}^t
					\norm{\grad \uu_1(s) - \grad \uu_2(s)}_{L^\iny(\Omega)}
						\norm{\grad \eta_2(t_1, s; \cdot)}_{L^\iny(\Omega)}\, ds
				\\
			&\qquad\qquad
			+ \int_{t_1}^t
					\norm{\uu_1(s)}_{\dot{C}^1}
							\norm{\grad \eta_1(t_1, s; \cdot)
								- \grad \eta_2(t_1, s; \cdot)}_{L^\iny(\Omega)} \, ds
				\\
			&\qquad
			\le C(T, M) [T e^{(C(T, M) T} \theta_\beta]^\al T + C(M, T) T \theta_\beta \\
			&\qquad\qquad
				+ C(M, T) \int_{t_1}^t
					\norm{\grad \eta_1(t_1, s; \cdot)
						- \grad \eta_2(t_1, s; \cdot)}_{L^\iny(\Omega)} \, ds.
	\end{align*}
	In the last inequality, we used \cref{L:BSLaw} to conclude that
	$\norm{\grad \uu_1(s) - \grad \uu_2(s)}_{L^\iny(\Omega)}
	\le \norm{\ww(s)}_{C^{1, \beta}(\Omega)}
	\le C \norm{\curl \ww(s)}_{C^\beta(\Omega)} + C \norm{\ww(s)}_H
	\le C \theta_\beta$.
	Taking the supremum over $\x$ and applying \Gronwalls Lemma as before
	gives the second bound.
\end{proof}

\begin{lemma}\label{L:WWidth}
	Letting $\abs{W}$ be the Lebesgue measure of $W := Q \setminus (V_+ \cup V_-)$,
	we have
	\begin{align*}
		\abs{W} \le C(T, M) T^2 \theta_\beta.
	\end{align*}
\end{lemma}
\begin{proof}
	The set $W(t) := \set{\x \in \Omega \colon (t, \x) \in W}$ consists of all points lying
	between the surfaces $\Cal{S}_1(t)$ and $\Cal{S}_2(t)$. Any $\x_1 \in \Cal{S}_1(t)$ is of the form $\x_1 = \eta_1(0, t; \y)$
	for some $\y \in \Gamma_+$, and by \cref{L:FlowDiff}, the point
	$\x_2 = \eta_2(0, t; \y)$ is within a distance $\delta = C(T, M) T \theta_\beta$ of $\x_1$.
	That is, any point in $\Cal{S}_1(t)$ is within a distance $\delta$ of $\Cal{S}_2(t)$ and
	the relation is symmetric. So
	\begin{align*}
		W(t)
			\subseteq W_\delta(t) := \set{x \in \Omega \colon \dist(x, \Cal{S}_1(t)) \le \delta}.
	\end{align*}

	As we observed in \cref{S:FlowMap},
	$\Cal{S}_1(t)$ is at least $C^{1, \al}$ regular as a surface in $\Omega$,
	and so has finite Hausdorff measure;
	hence, we can see that
	$\abs{W_\delta(t)} \le C \delta$. Moreover, this constant can depend upon $T$
	and $M$, but is bounded over $[0, T]$, for as also observed in \cref{S:FlowMap},
	$\Cal{S}_1$ is at least $C^{1, \al}$ regular as a hypersurface in $Q$. Thus,
	$\abs{W} \le T \abs{W_\delta(t)} \le C(T, M) T^2 \theta_\beta$.
\end{proof}

\begin{proof}[\textbf{Proof of \cref{P:AContinuous}}]
We will show that
\begin{align}\label{e:AContinuousGoal}
	\norm{A \uu_1 - A \uu_2}_{L^\iny(Q)}
		&\le \FMOC(\theta_\beta).
\end{align}
Once we obtain \cref{e:AContinuousGoal}, we will have
\begin{align*}
	\norm{A \uu_1 - A \uu_2}_{L^\iny(Q)}
		&\le \FMOC(\theta_\beta)
		= \norm{\uu_1 - \uu_2}_{\WeakSArg{0}}
		\le \norm{\uu_1 - \uu_2}_{\WeakS},
\end{align*}
and continuity of $A \colon \Cal{K} \to \Cal{K}$ in $\WeakS$ will follow from \cref{L:WeakSBound}.

We will obtain the bound on $\norm{A \uu_1 - A \uu_2}_{L^\iny(Q)}$ in \cref{e:AContinuousGoal} by the following three steps:
\begin{enumerate}[label=(\emph{\Alph*})]
    \item
        Bound the difference in vorticities, $\Lambda \uu_1 - \Lambda \uu_2$, in $L^\iny(Q)$ assuming  zero forcing.
        
    \item
        Account for forcing in the bound on $\Lambda \uu_1 - \Lambda \uu_2$ in $L^\iny(Q)$.
        
    \item
        Bound  $A \uu_1 - A \uu_2$ in $L^\iny(Q)$ with the help of (B).
    \end{enumerate}

\medskip\noindent\textbf{(A) Vorticity}:
Letting $(t, \x) \in Q$, we must estimate $\abs{\Lambda \uu_1(t, \x) - \Lambda \uu_2(t, \x)}$. This involves three cases: \textbf{(1)} $(t, \x) \in V_-$, \textbf{(2)} $(t, \x) \in V_+$, \textbf{(3)} $(t, \x) \in W$, which we consider separately. We argue first without forcing.

\medskip
\noindent\textbf{(1)} Define, for $(t, \x) \in V_-$, $j = 1, 2$,
\begin{align}\label{e:pos0j}
	\pos_0^j
		= \pos_0^j(t, \x)
		:= \eta_j(t, 0; \x).
\end{align}
From \cref{e:LagrangianForm}, we can write,
\begin{align*}
	&\Lambda \uu_1(t, \x) - \Lambda \uu_2(t, \x)
		= \grad \eta_1(0, t; \pos_0^1) \bomega_0(\pos_0^1)
			- \grad \eta_2(0, t; \pos_0^2) \bomega_0(\pos_0^2)
		= I_1 + I_2,
\end{align*}
where
\begin{align*}
	I_1
		&:= \bomega_0(\pos_0^1)
				\cdot (\grad \eta_1(0, t; \pos_0^1)
						- \grad \eta_2(0, t; \pos_0^2)), \\
	I_2
		&:= (\bomega_0(\pos_0^1) - \bomega_0(\pos_0^2))
				\cdot \grad \eta_2(0, t; \pos_0^2).
\end{align*}
We also make the decomposition, $I_1 = \bomega_0(\pos_0^1) \cdot (I_1^1 + I_1^2)$, where
\begin{align*}
	I_1^1
		&:= \grad \eta_1(0, t; \pos_0^1)
						- \grad \eta_1(0, t; \pos_0^2), \\
	I_1^2
		&:= \grad \eta_1(0, t; \pos_0^2)
						- \grad \eta_2(0, t; \pos_0^2).
\end{align*}
Then,
\begin{align*}
	\norm{I_1}_{L^\iny(V_-)}
		&\le \norm{\bomega_0}_{L^\iny(\Omega)}
			\pr{\norm{I_1^1}_{L^\iny(V_-)} + \norm{I_1^2}_{L^\iny(V_-)}},
\end{align*}
with
\begin{align*}
	\norm{I_1^1}_{L^\iny(V_-)}
		&\le \norm{\grad \eta_1(0, t; \cdot)}_{\dot{C}^\al(\Omega)}
			\norm{\eta_1(t, 0; \cdot) - \eta_2(t, 0; \cdot)}_{L^\iny(\Omega)}^\al \\
		&\le C(T, M) T [T \theta_\beta]^\al
		\le C(T, M) T^{1 + \al} \theta_\beta^\al, \\
	\norm{I_1^2}_{L^\iny(V_-)}
		&\le \norm{\grad \eta_1(0, t; \cdot) - \grad \eta_2(0, t; \cdot)}_{L^\iny(\Omega)}
		\le C(T, M) T [\theta_\beta + \theta_\beta^{\al}],
\end{align*}
where we applied \cref{L:FlowDiff}.
Similarly, applying \cref{L:HolderDiff,L:FlowDiff},
\begin{align*}
	\norm{I_2}_{L^\iny(V_-)}
		&\le \norm{\bomega_0}_{\dot{C}^\al(\Omega)}
			\norm{\eta_1(t, 0; \cdot) - \eta_2(t, 0; \cdot)}_{L^\iny(V_-)}^\al
			\norm{\grad \eta_2(0, t, \cdot)}_{L^\iny(V_-)} \\
		&\le C(T, M) M [C(T, M) T \theta_\beta]^\al.
\end{align*}

Dropping the dependence upon $M$ or the initial data, which play no role here, we conclude
\begin{align*}
	\norm{\Lambda \uu_1(t, \x) - \Lambda \uu_2(t, \x)}_{L^\iny(V_-)}
		&\le C(T) [\theta_\beta + \theta_\beta^{\al}].
\end{align*}

\medskip
\noindent\textbf{(2)} For $(t, \x) \in V_+$, we have
\begin{align*}
	\Lambda \uu_1(t, \x) - \Lambda \uu_2(t, \x)
		&= \H_1(\mu_1(t, \x))
			\cdot \grad \eta_1 (\time_1(t, \x), t; \pos_1(t, \x)) \\
		&\qquad
			- \H_2(\mu_2(t, \x))
			\cdot \grad \eta_2 (\time_2(t, \x), t; \pos_2(t, \x)) \\
		&= J_1 + J_2 + J_3,
\end{align*}
where $\H_j(t, \x)$ is defined in \cref{e:H} for $\uu_j$, and
\begin{align*}
	J_1
		&:= \H_1(\mu_1(t, \x))
			\cdot (\grad \eta_1 (\time_1(t, \x), t; \pos_1(t, \x))
				- \grad \eta_2 (\time_1(t, \x), t; \pos_1(t, \x))), \\
	J_2
		&:= \H_1(\mu_1(t, \x))
			\cdot (\grad \eta_2 (\time_1(t, \x), t; \pos_1(t, \x))
				- \grad \eta_2 (\time_2(t, \x), t; \pos_2(t, \x))), \\
	J_3
		&:= (\H_1(\mu_1(t, \x)) - \H_2(\mu_2(t, \x))
			\cdot \grad \eta_2 (\time_2(t, \x), t; \pos_2(t, \x)).
\end{align*}

Now, since $\H_j(s, \y) = \bomega_j(s, \y)$ for $(s, \y) \in [0, T] \times \Gamma_+$, we have, using \cref{L:FlowDiff},
\begin{align*}
	\norm{J_1}_{L^\iny(V_+)}
		&\le \norm{\bomega_1}_{L^\iny(Q)}
			\norm{\grad \eta_1(\cdot, t; \cdot) - \grad \eta_2(\cdot, t; \cdot))}
				_{L^\iny(Q)}
		\le C(T, M) [\theta_\beta + \theta_\beta^\al],
\end{align*}
where we also used $\cond_0$.
For $J_2$, we have, using \cref{L:HolderDiff,L:muDiff},
\begin{align*}
	\norm{J_2}_{L^\iny(V_+)}
		&\le \norm{\bomega_1}_{L^\iny(Q)}
			\norm{\grad \eta_2}_{\dot{C}^\al(Q)}
			\norm{(\time_1(t, \x), \pos_1(t, \x))
				- (\time_2(t, \x), \pos_2(t, \x))}
				_{L^\iny(Q)}^\al \\
		&\le C(T, M) \norm{\mu}_{L^\iny(U_+)}^\al
		\le C(T, M) \theta_\beta^\al.
\end{align*}

For $J_3$, we have
\begin{align*}
	J_3
		&\le \norm{\H_1(\mu_1(t, \x)) - \H_2(\mu_2(t, \x)}_{L^\iny(U_+)}
			\norm{\grad \eta_2}_{L^\iny(Q)}.
\end{align*}
But, $\norm{\grad \eta_2}_{L^\iny(Q)} \le C(T, M)$ by \cref{L:QuickFlow}, and, using \cref{L:HolderDiff},
\begin{align*}
	&\norm{\H_1(\mu_1(t, \x)) - \H_2(\mu_2(t, \x)}_{L^\iny(U_+)} \\
		&\qquad
		\le \norm{\H_1(\mu_1(t, \x)) - \H_2(\mu_1(t, \x)}_{L^\iny(U_+)}
			+ \norm{\H_2(\mu_1(t, \x)) - \H_2(\mu_2(t, \x)}_{L^\iny(U_+)} \\
		&\qquad
		\le \norm{\H_1 - \H_2}_{L^\iny([0, T] \times \Gamma_+)}
			+ \norm{\H_2}_{\dot{C}^\al([0, T] \times \Gamma_+)}
				\norm{\mu}_{L^\iny}^\al \\
		&\qquad
		\le \norm{\bomega_1 - \bomega_2}_{L^\iny([0, T] \times \Gamma_+)}
			+ C(T, M)  \theta_\beta^\al
		\le C(T, M) [\theta_\beta + \theta_\beta^\al],	
\end{align*}
where in the second-to-last inequality we used the bounds on $\H_1$ and $\H_2$ from \cref{P:HBounds} and use that $\H_j = \bomega_j$ on $[0, T] \times \Gamma_+$, since $\Lambda \uu$ solves \OurLinear.

Combined, we see that
\begin{align*}
	&\norm{\Lambda \uu_1(t, \x) - \Lambda \uu_2(t, \x)}_{L^\iny(V_+)}
		\le C(T, M) [\theta_\beta + \theta_\beta^\al].
\end{align*}

\medskip
\noindent\textbf{(3)} Now assume $(t, \x) \in W$. Applying \cref{L:MOCBound} with the Lipschitz \MOC, $r \mapsto \norm{\Lambda \uu_1 - \Lambda \uu_2}_{\dot{C}^\al} r \le M r$, 
\begin{align*}
	\norm{\Lambda \uu_1 - \Lambda \uu_2}_{L^\iny(W)}
		&\le \FMOC \pr{\norm{\Lambda \uu_1 - \Lambda \uu_2}_{L^2(W)}}.
\end{align*}
From \cref{L:WWidth},
\begin{align*}
    \norm{\Lambda \uu_1 - \Lambda \uu_2}_{L^2(W)}
        &\le \norm{\Lambda \uu_1 - \Lambda \uu_2}_{L^\iny(W)} \abs{W}^{\frac{1}{2}}
        \le C M \abs{W}^{\frac{1}{2}}
        \le C(T, M) \theta_\beta,
\end{align*}
which then gives $\norm{\Lambda \uu_1 - \Lambda \uu_2}_{L^\iny(W)} \le \FMOC(\theta_\beta)$. We conclude that $\norm{\Lambda \uu_1 - \Lambda \uu_2}_{L^\iny(Q)} \le \FMOC(\theta_\beta)$.

\medskip\noindent\textbf{(B) Accounting for forcing}: 
To treat forcing, let $\G_\pm^j$ be given by \cref{e:LagrangianForm} for $\eta_j$. Then
\begin{align*}
	&\norm{\G_\pm^1 - \G_\pm^2}_{L^\iny(V_\pm)} \\
		&\qquad
		\le \int_0^T \norm{
						\grad \eta_1(s, t; \eta_1(t, s; \x)) \g(s, \eta_1(t, s; \x))
						- \grad \eta_2(s, t; \eta_2(t, s; \x)) \g(s, \eta_2(t, s; \x))
						}_{L^\iny(\Omega)}
					\, ds.
\end{align*}
But,
\begin{align*}
	&\norm{
		\grad \eta_1(s, t; \eta_1(t, s; \x)) \g(s, \eta_1(t, s; \x))
		- \grad \eta_2(s, t; \eta_2(t, s; \x)) \g(s, \eta_2(t, s; \x))
		}_{L^\iny(\Omega)} \\
		&\quad
		\le
			\norm{\grad \eta_1(s, t; \eta_1(t, s; \x)) \g(s, \eta_1(t, s; \x))
			- \grad \eta_2(s, t; \eta_1(t, s; \x)) \g(s, \eta_1(t, s; \x))
			}_{L^\iny(\Omega)} \\
		&\qquad
			+ \norm{\grad \eta_2(s, t; \eta_1(t, s; \x)) \g(s, \eta_1(t, s; \x))
			- \grad \eta_2(s, t; \eta_2(t, s; \x)) \g(s, \eta_1(t, s; \x))
			}_{L^\iny(\Omega)} \\
		&\qquad
			+ \norm{\grad \eta_2(s, t; \eta_2(t, s; \x)) \g(s, \eta_1(t, s; \x))
			- \grad \eta_2(s, t; \eta_2(t, s; \x)) \g(s, \eta_2(t, s; \x))
			}_{L^\iny(\Omega)} \\
		&\quad
		\le
			\norm{\grad \eta_1 - \grad \eta_2}_{L^\iny([0, T]^2 \times \Omega)}
				\norm{\g}_{L^\iny(Q)} \\
		&\qquad
			+ \norm{\grad \eta_2}_{\dot{C}^\al([0, T]^2 \times \Omega)}
				\norm{\grad \eta_1 - \grad \eta_2}_{L^\iny([0, T]^2 \times \Omega)}^\al
				\norm{\g}_{L^\iny(Q)}  \\
		&\qquad
			+ \norm{\grad \eta_2}_{L^\iny([0, T]^2 \times \Omega))}
				\norm{\g}_{\dot{C}^\al} \norm{\eta_1 - \eta_2}_{L^\iny(Q)}^\al,
\end{align*}
where we used \cref{L:HolderComp,L:HolderDiff}.

Since $\g \in L^\iny(Q)$, while $\grad \eta_1$ and $\grad \eta_2$ are bounded in $\dot{C}^\al([0, T]^2 \times \Omega)$, by \cref{L:FlowDiff} we see that
\begin{align*}
	&\norm{\G_\pm^1 - \G_\pm^2}_{L^\iny(V_\pm)}
		\le C T [\theta_\beta + \theta_\beta^{\al}].
\end{align*}
Hence, the inclusion of forcing does not change our bounds on $\norm{\Lambda \uu_1(t, \x) - \Lambda \uu_2(t, \x)}_{L^\iny(V_\pm)}$ in \textbf{(1)}, \textbf{(2)}. And $\G_\pm^1$, $\G_\pm^2$ are bounded on $Q$, so the estimate on $\norm{\Lambda \uu_1 - \Lambda \uu_2}_{L^2(W)}$ in \textbf{(3)} is also unchanged.

In summary, what we have done so far is to show that
\begin{align}\label{e:LambdaDiffBound}
    \norm{\Lambda \uu_1 - \Lambda \uu_2}_{L^\iny(Q)}
        &\le \FMOC(\theta_\beta).
\end{align}

\medskip\noindent\textbf{(C) Velocity}:
From \cref{e:LinearVelEq}, we have,
\begin{align*}
    P_H A \uu_j(t, \x)
        &= \uu_0(\x) + \int_0^t \f(s, \x) \, ds
            - \int_0^t P_H[\uu_j(s) \cdot \ol{\bOmega}_j(s)](\x) \, ds.
\end{align*}
Then because $A \uu_1 - A \uu_2 = P_H (A \uu_1 - A \uu_2)$ and $P_H$ is continuous in $C^\beta(\Omega)$, which is an algebra,
\begin{align*}
    &\norm{A \uu_1(t) - A \uu_2(t)}_{L^\iny(\Omega)}
    	\le \norm{A \uu_1(t) - A \uu_2(t)}_{C^\beta(\Omega)} \\
        &\qquad
        \le C \int_0^t \norm{\ww}_{C^\beta(\Omega)} \norm{\ol{\bOmega}_1}_{C^\beta(\Omega)}
            + C \int_0^t \norm{\uu_1}_{C^\beta(\Omega)}
            	\norm{\ol{\bOmega}_1 - \ol{\bOmega}_2}_{C^\beta(\Omega)} \\
        &\qquad
        = C \int_0^t \norm{\ww}_{C^\beta(\Omega)} \norm{\ol{\bomega}_1}_{C^\al(\Omega)}
            + C \int_0^t \norm{\uu_1}_{C^\al(\Omega)}
            	\norm{\Lambda \uu_1 - \Lambda \uu_2}_{C^\beta(\Omega)} \\
        &\qquad
        \le C T \theta_\beta
            + C T \FMOC \pr{\norm{\Lambda \uu_1 - \Lambda \uu_2}_{L^\iny(Q)}}
        \le \FMOC (\theta_\beta),
\end{align*}
where we used \cref{P:AuCNalBound}, \cref{P:Koch22Sub}, \cref{e:LambdaDiffBound}, and \cref{L:WeakSBound}.

This gives \cref{e:AContinuousGoal}, which completes the proof.
\end{proof}

%
%
\section{Full inflow boundary condition satisfied}\label{S:FullBCs}

\noindent We now prove \cref{P:WeakVorticityIsVelocity}, which shows that a solution satisfying \OurWeakBCs also satisfies \OurInflow, and hence satisfies the full inflow boundary conditions. This can be done by defining $\H$ by \cref{e:H} and recovering the pressure using $N[\uu]$ of \cref{e:NuDef}, as already observed in \cite{AKM}.

\begin{proof}[\textbf{Proof of \cref{P:WeakVorticityIsVelocity}}]
	Our proof is inspired by the proof of Lemma 4.2.1 pages 156-159 of \cite{AKM}.
	Let
	\begin{align*}
		\ww = \uu^\BoldTau - \UU^\BoldTau, \quad
		P := p - p_r,
	\end{align*}
	where $p_r$ is the regularized pressure given by \cref{e:pr}.
	By \cref{P:omegaOnGammaPlus}, $\bomega = \W[\uu, p]$ on $[0, T] \times \Gamma_+$,
	where we recall that $\W[\uu, p]$ is defined in \cref{e:W}.
	From \cref{e:TruePressureEq}, \cref{e:pr}, and \cref{e:NuDef}, we see that
	on $\Gamma_+$,
	$
		\grad P \cdot \n
			= \dv_\Gamma (U^\n \ww)
	$.
	Hence, $P$ satisfies
	\begin{align*}
		\begin{cases}
			\Delta P = 0
				&\text{in } \Omega, \\
			\grad P \cdot \n = 0
				&\text{on } \Gamma_- \cup \Gamma_0, \\
			\grad P \cdot \n
				= \dv_\Gamma (U^\n \ww)
				&\text{on } \Gamma_+.
		\end{cases}
	\end{align*}
	Multiplying by $P$, integrating over $\Omega$, and integrating by parts over $\Gamma_+$ gives
	\begin{align}\label{e:gradPSquared}
		\norm{\grad P}_{L^2(\Omega)}^2
			&= - (\Delta P, P) + \int_{\Gamma_+} (\grad P \cdot \n) P
			= \int_{\Gamma_+} \dv_\Gamma (U^\n \ww) P
			= - \int_{\Gamma_+} U^\n \w \cdot \grad_\Gamma P.
	\end{align}
	
	By \cref{e:omegaIsW} and the assumption that $\H = \bomega := \curl \uu$
        on $\Gamma_+$, we know that
	$U^\n [\H^\BoldTau]^\perp = U^\n [\W^\BoldTau[\uu, p]]^\perp$.
	Using also that $(\vv^\perp)^\perp = - \vv$, 
	we have, from \cref{e:W,e:H}, that on $\Gamma_+$,
	\begin{align*}
			&\prt_t \UU^\BoldTau + \grad_\Gamma \pr{p_r + \frac{1}{2} \abs{\UU}^2}
			 	 	- \f^\BoldTau
				+ \curl_\Gamma \UU^\BoldTau [\uu^\BoldTau]^\perp
				= \H \\
				&\qquad
				= \bomega
				= \prt_t \uu^\BoldTau + \grad_\Gamma \pr{p + \frac{1}{2} \abs{\uu}^2}
			 	 	- \f^\BoldTau 
					+ \curl_\Gamma \uu^\BoldTau [\uu^\BoldTau]^\perp.
	\end{align*}
	Subtracting the left hand side from the right hand side, we have
	\begin{align*}
		0
			&= \grad_\Gamma P
				+ \frac{1}{2} \grad_\Gamma (\abs{\uu}^2 - \abs{\UU}^2)
				+ \prt_t \ww + \curl_\Gamma \ww [\uu^\BoldTau]^\perp.
	\end{align*}
	But, 
	$\omega^\n = H^\n$ on $\Gamma_+$, which gives
	$
		\curl_\Gamma \UU^\BoldTau = \curl_\Gamma \uu^\BoldTau
	$.
	Hence, $\curl_\Gamma \ww = 0$, so
	\begin{align*}
		\grad_\Gamma P
			&= - \prt_t \ww
				- \frac{1}{2} \grad_\Gamma (\abs{\uu}^2 - \abs{\UU}^2).
	\end{align*}
	
	Returning to \cref{e:gradPSquared},
	we thus have
	\begin{align*}
		\norm{\grad P}_{L^2(\Omega)}^2
			&= \int_{\Gamma_+} U^\n \ww \cdot \prt_t \ww 
				+ \frac{1}{2} \int_{\Gamma_+} U^\n \ww \cdot
					\grad_\Gamma (\abs{\uu}^2 - \abs{\UU}^2).
	\end{align*}
	Now,
	\begin{align*}
		\int_{\Gamma_+} U^\n \ww \cdot \prt_t \ww 
			&= \frac{1}{2} \int_{\Gamma_+} U^\n \prt_t \abs{\ww}^2
			= \frac{1}{2} \int_{\Gamma_+} \prt_t [U^\n \abs{\ww}^2]
				- \frac{1}{2} \int_{\Gamma_+} \prt_t U^\n \, \abs{\ww}^2 \\
			&= \frac{1}{2} \diff{}{t} \int_{\Gamma_+} U^\n \abs{\ww}^2
				- \frac{1}{2} \int_{\Gamma_+} \prt_t U^\n \, \abs{\ww}^2,
	\end{align*}
	so
	\begin{align}\label{e:ddtUnw2}
		\diff{}{t} \int_{\Gamma_+} U^\n \abs{\ww}^2
			&= \int_{\Gamma_+} \prt_t U^\n \abs{\ww}^2
				-  \int_{\Gamma_+} U^\n \ww \cdot
					\grad_\Gamma (\abs{\uu}^2 - \abs{\UU}^2)
				+ 2 \norm{\grad P}_{L^2(\Omega)}^2.
	\end{align}
	
	Writing
	$\abs{\UU}^2 - \abs{\uu}^2 = \abs{\uu^\BoldTau}^2 - \abs{\UU^\BoldTau}^2
	= \ww \cdot \vv$ on $\Gamma_+$, since $U^\n = u^\n$, 
	where $\vv := \UU^\BoldTau + \uu^\BoldTau$,
	we have
	\begin{align*}
		\int_{\Gamma_+} U^\n \ww &\cdot \grad_\Gamma (\abs{\uu}^2 - \abs{\UU}^2)
			= \int_{\Gamma_+} U^\n \ww \cdot \grad_\Gamma (\ww \cdot \vv) \\
			&= \int_{\Gamma_+} U^\n (\ww \cdot \grad_\Gamma \vv) \cdot \ww
				+ \int_{\Gamma_+} U^\n (\ww \cdot \grad_\Gamma \ww) \cdot \vv \\
			&= \int_{\Gamma_+} U^\n (\ww \cdot \grad_\Gamma \vv) \cdot \ww
				- \frac{1}{2} \int_{\Gamma_+} \abs{\ww}^2 \dv_\Gamma (U^\n \vv).
	\end{align*}
	For the last term above, we used that
	$U^\n (\ww \cdot \grad_\Gamma \ww) \cdot \vv  = (1/2) U^\n \vv \cdot \grad_\Gamma \abs{\ww}^2 $
	and integrated by parts via \cref{L:IBPOnBoundary}.
	Then because $\vv$ and $U^\n$ are sufficiently regular, we have
	\begin{align*}
		\abs[\bigg]{\int_{\Gamma_+} U^\n \ww \cdot \grad_\Gamma (\abs{\uu}^2 - \abs{\UU}^2)}
			\le C \int_{\Gamma_+} \abs{\ww}^2.
	\end{align*}
	
	Changing sign in \cref{e:ddtUnw2} and integrating in time, we see that
	\begin{align*}
		\int_{\Gamma_+} &\abs{U^\n(t)} \abs{\ww(t)}^2
			= -\int_{\Gamma_+} U^\n(t) \abs{\ww(t)}^2 \\
			&\le -\int_0^t \int_{\Gamma_+} \prt_t U^\n \abs{\ww}^2
				+  \int_0^t \int_{\Gamma_+} U^\n \ww \cdot
					\grad_\Gamma (\abs{\uu}^2 - \abs{\UU}^2)
				- 2 \int_0^t \norm{\grad P}_{L^2(\Omega)}^2
				\\
			&\le C \int_0^t \int_{\Gamma_+} \abs{\ww(s)}^2 \, ds
				- 2 \int_0^t \norm{\grad P}_{L^2(\Omega)}^2
			\le C \int_0^t \int_{\Gamma_+} \abs{\ww(s)}^2 \, ds.
	\end{align*}
	In the first equality we used that $U^\n < 0$ on $\Gamma_+$,
	in the second equality we used that $\ww(0) = 0$, and in the third
	equality we used that $\prt_t U^\n$ is bounded.
	
	Now since $\abs{U^\n}$ is bounded away from zero, we have
	\begin{align*}
		\int_{\Gamma_+} \abs{\ww(t)}^2
			\le C \int_0^t \int_{\Gamma_+} \abs{\ww(s)}^2 \, ds,
	\end{align*}
	and we conclude from \Gronwalls Lemma that $\ww \equiv 0$. This means that
	$\uu^\BoldTau = \UU^\BoldTau$, so \OurInflow holds. 
\end{proof}

\begin{remark}\label{R:ImpermeableOnly}
	If $\Gamma_0 = \Gamma$, the classical setting of impermeable boundary
	conditions on the whole boundary, our proof of existence and uniqueness still applies,
	though a number
	of things trivialize. First, no vorticity is transported off of the boundary,
	so there is no need for the pressure estimates in \cref{S:Pressure},
	and $U_-$ is all of $Q$, so many
	of the flow map constructs, such as $S$, $\time$, and $\pos$ are unnecessary.
	And, of course, none of the estimates involving $U_+$ are needed.
\end{remark}

%
%
\section{Vorticity boundary conditions}\label{S:VorticityBCs}

\noindent

\begin{proof}[\textbf{Proof of \cref{T:VorticityBCs}}]
	The proof of existence is the same as that for \cref{T:WellPosedness},
	though with substantial simplifications.
	Because $\H$ is given with sufficient regularity, it satisfies
	\begin{align*}
		\norm{\H}_{L^\iny([0, T] \times \Gamma_+)}
			\le c_0, \quad
		\norm{\H}_{C^{N, \al}([0, T] \times \Gamma_+)}
			\le c_0.
	\end{align*}
	Hence, there are no pressure estimates involved, so
	the condition in \cref{e:RangeOfCurlCondSpecial} immediately gives
	\cref{e:RangeOfCurlCond}, and there is no need to appeal
	to \cref{P:H0CondSatisfied}. Since we only require $\uu \cdot \n = U^\n$
	on $\Gamma_+$, we simplify the definition of $\Dom_N(A)$ in \cref{e:DomNA} to
	\begin{align*}
	\Dom_N(A)
		:= \set{&\uu \in \uSolSpaceN \colon
			\uu(0) = \uu_0},
	\end{align*}
	and there is no need to invoke \cref{P:WeakVorticityIsVelocity} or \cref{L:DomNA}.
	Otherwise, the remainder of the proof of existence proceeds unchanged.
	
For uniqueness when $N \ge 1$, let $\bomega_j = \curl \uu_j$, $j = 1, 2$, and let $\ww = \uu_1 - \uu_2$. Then $\ww \in H_0$, since $\uu_1$, $\uu_2$ have the same prescribed harmonic component, $\uu_c$. Let
\begin{align*}
	\bmu := \curl \ww = \bomega_1 - \bomega_2.
\end{align*}
Since $N \ge 1$, we have enough regularity to write $\prt_t \bomega_j + \uu_j \cdot \grad \bomega_j = \bomega_j \cdot \grad \uu_j + \curl \f$, and subtracting this relation for $j = 2$ from that for $j = 1$ gives
\begin{align}\label{e:VortDiffEq}
	\prt_t \bmu
		&+ \uu_1 \cdot \grad \bmu + \ww \cdot \grad \bomega_2
		= \bomega_1 \cdot \grad \ww + \bmu \cdot \grad \uu_2.
\end{align}
Multiplying by $\bmu$, integrating over $\Omega$, and using that $(\uu_1 \cdot \grad \bmu, \bmu) = (1/2) (\uu_1, \grad \abs{\bmu}^2)$, gives
\begin{align}\label{e:VortEnergyPairing}
	\begin{split}
	\frac{1}{2} &\diff{}{t} \norm{\bmu}^2
			+ \frac{1}{2} \int_\Omega \uu_1 \cdot \grad \abs{\bmu}^2
		= - (\ww \cdot \grad \bomega_2, \bmu)
			+ (\bomega_1 \cdot \grad \ww, \bmu) + (\bmu \cdot \grad \uu_2, \bmu) \\
		&\le \frac{1}{2} \norm{\grad \bomega_2}_{L^\iny} \norm{\ww}^2 + \frac{1}{2} \norm{\bmu}^2
			+ \frac{1}{2} \norm{\bomega_1}_{L^\iny} \norm{\grad \ww}^2
			+ \frac{1}{2} \norm{\bmu}^2
			+ \norm{\grad \uu_2}_{L^\iny} \norm{\bmu}^2,
	\end{split}
\end{align}
where $\norm{\cdot} := \norm{\cdot}_{L^2(\Omega)}$ here. Elements of $H$ have mean zero, so by \Poincares inequality, $\norm{\ww} \le C \norm{\grad \ww}$. Moreover, since $\ww \in H_0$, we have $\norm{\grad \ww} \le C \norm{\bmu}$ and so obtain
\begin{align*}
	\diff{}{t} \norm{\bmu}^2
		&\le - \int_\Omega \uu_1 \cdot \grad \abs{\bmu}^2 + C \norm{\bmu}^2.
\end{align*}
We note that $\grad \bomega_2 \in L^\iny([0, T] \times \Omega)$ by the $N = 1$ existence result. But,
\begin{align*}
	- \int_\Omega \uu_1 \cdot \grad \abs{\bmu}^2
		&= \int_\Omega \dv \uu_1 \, \abs{\bmu}^2
			- \int_\Gamma U^\n \abs{\bmu}^2
		= - \int_{\Gamma_-} U^\n \abs{\bmu}^2
		\le 0,
\end{align*}
so we conclude from Gronwall's lemma, since $\bmu(0) = 0$, that $\bmu \equiv 0$. That is, $\uu_1 = \uu_2$.

Finally, from \cref{e:OurVorticity}$_1$, we have
\begin{align*}
	\prt_t \uu^\BoldTau + (\uu \cdot \grad \uu)^\BoldTau
		= (\f - \grad p )^\BoldTau + \zz^\BoldTau.
\end{align*}
From $\cond_0$, then, we see that $\zz^\BoldTau(0) = 0$. Since also $z^\n(0) = 0$, we know that $\zz(0) = 0$.
\end{proof}

\newcommand{\Acknowledgments}{\noindent Gie was partially supported by a Simons Foundation Collaboration Grant for Mathematicians; Research R-II Grant, Office of EVPRI, University of Louisville; Brain Pool Program through the National Research Foundation of Korea (NRF) (grant number: 2020H1D3A2A01110658). Mazzucato was partially supported by the US National Science Foundation Grant DMS-1909103 and Simons Foundation Grant 1036502. Part of this work was prepared while Kelliher and Mazzucato were participating in a program hosted by the Simons Laufer Mathematical Sciences Research Institute (SLMath) in Berkeley, California, in Spring 2021 and again in Summer 2023, supported by the National Science Foundation under Grant No. DMS-1928930.
Mazzucato would like to thank the Isaac Newton Institute for
Mathematical Sciences, Cambridge, for support and hospitality during the
programme, \textit{Mathematical aspects of turbulence: where do we stand?},
where work on this paper was partially undertaken.
The work of the Institute is supported by EPSRC grant no EP/R014604/1. The authors thank Igor Kukavica, Helena Nussenzveig Lopes, Marco Sammartino, and Wojciech S. O\.za\'nski for stimulating discussions.
}

\section*{Acknowledgements}
\Acknowledgments

%
%

%
%
\addtocontents{toc}{\protect\vspace{0.65em}}

\appendix

%
%
\section{\Holder space lemmas}\label{A:HolderLemmas}

\noindent We collect here a number of estimates in \Holder spaces, defined in \cref{S:BSLaw}, which we use throughout much of this paper. We include proofs only of the less standard ones.

\begin{lemma}\label{L:HolderProd}
	Let $f, g \in C^\al(U)$. Then
	\begin{align*}
		\norm{fg}_{C^\al}
			&\le \norm{f}_{C^\al} \norm{g}_{C^\al}, \\
		\norm{fg}_{\dot{C}^\al}
			&\le \norm{f}_{L^\iny} \norm{g}_{\dot{C}^\al}
				+ \norm{g}_{L^\iny} \norm{f}_{\dot{C}^\al}, \\
		\norm{fg}_{C^\al}
			&\le \norm{f}_{L^\iny} \norm{g}_{L^\iny} 
				+ \norm{f}_{L^\iny} \norm{g}_{\dot{C}^\al}
				+ \norm{g}_{L^\iny} \norm{f}_{\dot{C}^\al}, \\
			&\le
				\norm{f}_{L^\iny} \norm{g}_{C^\al}
				+ \norm{g}_{L^\iny} \norm{f}_{C^\al}, \\
		\norm{fg}_{C^\al}
			&\le
				\norm{f}_{L^\iny} \norm{g}_{\dot{C}^\al}
				+ \norm{g}_{L^\iny} \norm{f}_{C^\al},
	\end{align*}
	where $\dot{C}^\al$ is the \Holder semi-norm, as in \cref{e:HolderNormDef}.
	Also, for any $\beta \in (0, \al)$, allowing $\al = 1$, we have the interpolation inequality,
	\begin{align*}
		\norm{f}_{\dot{C}^\beta}
			&\le 2 \norm{f}_{\dot{C}^\al}^{\frac{\beta}{\al}}
				\norm{f}_{L^\iny}^{1 - \frac{\beta}{\al}}.
	\end{align*}
\end{lemma}

\begin{lemma}\label{L:HolderComp}
	Let $U, V$ be open subsets of $\R^n$, $\al \in (0, 1]$, and $k \ge 1$ an integer. 
	If $f \in C^{k, \al}(U)$ and $g \in C^{k + 1}(V)$ with $g(V) \subseteq U$
	then
	\begin{align}\label{e:HolderComp}
		\begin{split}
		\norm{f \circ g}_{\dot{C}^\al(V)}
			&\le \norm{f}_{\dot{C}^\al(U)} \norm{g}_{Lip(V)}^\al, \\
		\norm{f \circ g}_{C^\al(V)}
			&\le \norm{f}_{L^\iny(U)}
				+ \norm{f}_{\dot{C}^\al(U)} \norm{g}_{Lip(V)}^\al
			\le \norm{f}_{C^\al(U)} \brac{1 + \norm{g}_{Lip(V)}^\al}, \\
		\norm{f \circ g}_{C^{k, \al}(V)}
			&\le C(k) \norm{f}_{C^{k, \al}(U)} \brac{1 + \norm{g}_{C^{k + 1}(V)}}^{k + 1},
		\end{split}
	\end{align}
	where $Lip$ is the homogeneous Lipschitz semi-norm
	and $\dot{C}^\al$ is the homogeneous \Holder norm.
\end{lemma}

\begin{lemma}\label{L:HolderDiff}
	Let $U, V$ be open subsets of $\R^d$, $d \ge 1$, and let $\al \in (0, 1]$.
	Assume that the domain of $f$ is $U$ and the domains of $g$ and $h$ are $V$,
	with  $g(V), h(V) \subseteq U$. Then
	\begin{align*}
		\norm{f \circ g - f \circ h}_{L^\iny(V)}
			&\le \norm{f}_{\dot{C}^\al(U)} \norm{g - h}_{L^\iny(V)}^\al.
	\end{align*}
\end{lemma}

We also have the following interpolation-like inequality:

\begin{lemma}\label{L:PreNirenberg}
	Let $U$ be a bounded open subset of $\R^d$, $d \ge 1$, let $n \ge 1$,
	and $\grad^n f \in C^\al(U)$.
	Then 
	\begin{align*}
		\norm{\grad^n f}_{L^\iny(U)}
			&\le C \norm{f}_{C^{n, \al}(U)}^a
				\norm{f}_{L^2(U)}^{1 - a},
	\end{align*}
	where
	\begin{align*}
		a = a_n = \frac{2n + d}{2n  + d + 2 \al} < 1.
	\end{align*}
\end{lemma}
\begin{proof}
	First extend $f$ continuously to all of $\R^d$  in all
	\Holder spaces,
	as can be done using the extension operator
	in Theorem $5'$, chapter VI of \cite{Stein1970}. Applying a cutoff function,
	we can insure that the extension, which we continue to call $f$,
	has support with a diameter no more than twice $\diam(U)$.

	Then
	\begin{align*}
		\begin{split}
		\norm{\grad^n  f}_{L^\iny(U)}
			&= \sup_{\x \in \supp f}
				\abs{\grad^n  f(\x)}
			= \sup_{\x \in \supp f}
				\abs{\grad^n  f(\x) - \grad^n  f(\x_0)}
			\le R,
		\end{split}
	\end{align*}
	where $\x_0$ is a fixed point in $(\supp f)^C$ (so $\grad \f(x_0) = 0$) and
	\begin{align*}
		R
			&= \sup_{\x \in \supp f} \abs{\x - \x_0}^\al
				\sup_{\x \in \supp f} 
				\frac{\abs{\grad^n  f(\x) - \grad^n  f(\x_0)}}
					{\abs{\x - \x_0}^\al}
			= \sup_{\x \in \supp f}  \abs{\x - \x_0}^\al
				\norm{\grad^n (f(s \cdot))}_{\dot{C}^\al(\R^d)}.
	\end{align*}
	In particular, 
	\begin{align}\label{e:gradfnBound}
		\norm{\grad^n f}_{L^\iny(\R^d)}
			&\le R
				+ \norm{f}_{L^2(\R^d)}
	\end{align}
	for all $f \in C_0^\iny(\R^d)$.

	Following the scaling argument in the proof of Proposition 13.3.4 of \cite{Taylor2011VolIII},
	we write \cref{e:gradfnBound}
	schematically in the form
	$Q \le R + P$.
	Replacing $f(\cdot)$ with $f(s \cdot)$, we have
	$
		\grad^n (f(s \x))
			= s^n \grad f(s \x).
	$
	This gives
	$
		\norm{\grad^n (f(s \cdot))}_{L^\iny(\R^d)}
			= s^n \norm{\grad f}_{L^\iny(\R^d)}
	$
	and
	$
		\norm{f(s \cdot)}_{L^2(\R^d)}
			= s^{-\frac{d}{2}} \norm{f}_{L^2(\R^d)}.
	$
	Also, $R$ becomes
	\begin{align*}
		\sup_{\x \in \supp f} \abs{s \x - s \x_0}^\al
				\sup_{\x \in \supp f}
				s^n \frac{\abs{\grad^n  f(s \x) - \grad^n  f(s \x_0)}}
					{\abs{s \x - s \x_0}^\al}
			= s^{n + \al} R.
	\end{align*}
	Thus, $Q \le R + P$ becomes
	\begin{align*}
		s^n Q \le s^{n + \al} R + s^{-\frac{d}{2}} P
			\implies
			Q \le s^\al R + s^{-(n + \frac{d}{2})} P.
	\end{align*}
	
	As in \cite{Taylor2011VolIII}, we conclude that
	\begin{align*}
		\norm{\grad^n f}_{L^\iny(\R^d)}
			\le \norm{\grad^n f}_{\dot{C}^\al(\R^d)}^a
					\norm{f}_{L^2(\R^d)}^{1 - a}
			\le C \norm{\grad^n f}_{\dot{C}^\al(U)}^a
					\norm{f}_{L^2(U)}^{1 - a}
	\end{align*}
	as long as $\al a = (n + \frac{d}{2})(1 - a)$, which gives the
	stated value of $a$ and the stated estimate, using the continuity of the
	extension operator.	
\end{proof}

The inequality in \cref{L:PreNirenberg} is similar to that in the lemma on page 126 of \cite{Nirenberg1959}, used by the authors of \cite{AKM} (for $N = 0$).

\begin{lemma}\label{L:Interpolation}
	Let $U$ be a bounded open subset of $\R^d$, $d \ge 1$, let $n \ge 1$,
	and suppose that $f \in C^{n, \al}(U)$. Let $a_n$ be as in \cref{L:PreNirenberg}.
	For any $\beta \in (0, \al)$,
	\begin{align*}
		\norm{f}_{C^{n, \beta}(U)}
			\le \norm{f}_{L^\iny(U)}
				+ C &\brac{\norm{f}_{C^{n, \al}(U)}^{a_1} + \norm{f}_{C^{n, \al}(U)}^{a_n}}
					\brac{\norm{f}_{L^2(U)}^{1 - {a_1}} + \norm{f}_{L^2(U)}^{1 - {a_n}}} \\
				&+  C \norm{f}_{C^{n, \al}(U)}^{a'}
				\norm{f}_{L^2(U)}^{1 - a'},
	\end{align*}
	where
	\begin{align*}
		a'
			= (\beta/\al) + a_n (1 - \beta/\al)
			< 1.
	\end{align*}
\end{lemma}

\begin{lemma}\label{L:LinfTimeDiff}
    Recalling the definition of $\dot{C}^{0, \beta}_t(Q)$ in \cref{e:Crtx},
    for any $\beta \in (0, 1]$ and $f \in \dot{C}^{0, \beta}(Q)$,
    \begin{align*}
        \norm{f(t_1) - f(t_2)}_{L^\iny(\Omega)}
            &\le \norm{f}_{\dot{C}^{0, \beta}_t(Q)} \abs{t_1 - t_2}^\beta.   
    \end{align*}
\end{lemma}
\begin{proof}
    For any $t_1, t_2 \in [0, T]$ with $t_1 \ne t_2$,
    \begin{align*}
        &\norm{f(t_1) - f(t_2)}_{L^\iny(\Omega)}
            = \sup_{\x \in \Omega}
                \frac{\abs{f(t_1, \x) - f(t_2, \x)}}{\abs{t_1 - t_2}^\beta}
                    \abs{t_1 - t_2}^\beta
            \le \norm{f}_{\dot{C}_t(Q)} \abs{t_1 - t_2}^\beta.
            \qedhere
    \end{align*}
\end{proof}

\begin{lemma}\label{L:CNTimeDiffBound}
	Let $f \in \vortSpaceN$ for some $N \ge 0$ and $\al \in (0, 1]$. Then
	for any $t_1, t_2 \in [0, T]$,
	\begin{align*}
		&\norm{f(t_1) - f(t_2)}_{C^N(\Omega)}
			\le \norm{f}_\vortSpaceN \abs{t_1 - t_2}^\al.
	\end{align*}
\end{lemma}
\begin{proof}
    We have, applying \cref{L:LinfTimeDiff},
    \begin{align*}
        \norm{f(t_1) - f(t_2)}_{C^N(\Omega)}
            &= \sum_{k = 0}^N \norm{\grad^k (f(t_1) - f(t_2))}_{L^\iny(\Omega)} \\
            &\le \sum_{k = 0}^{N - 1} \norm{\grad^k f}_{\dot{C}^{0, 1}_t(Q)} \abs{t_1 - t_2}
                + \norm{\grad^N f}_{\dot{C}^{0, \al}_t(Q)} \abs{t_1 - t_2}^\al \\
            &\le \sum_{k = 0}^N \norm{\grad^k f}_{\dot{C}^\al_t(Q)} \abs{t_1 - t_2}^\al
            \le \norm{f}_\vortSpaceN \abs{t_1 - t_2}^\al.
            \qedhere
    \end{align*}
\end{proof}

\begin{cor}\label{C:CNTimeDiffBound}
	If $f \in \vortSpaceN$ for some $N \ge 0$ and $\al \in (0, 1]$ then
	\begin{align*}
		&\norm{f(t) - f(0)}_{C^N(Q)}
			\le C \norm{f}_\vortSpaceN T^\al.
	\end{align*}
\end{cor}

\cref{L:MOCBound} is adapted from Lemma 8.3 of \cite{KExpandingDomain}.
\begin{lemma}\label{L:MOCBound}
	Suppose that $f_j \colon \R^d \to \R$,
	$j = 1, 2$, each have the modulus of continuity $\Theta$, with
	$\Theta \colon [0, \iny) \to [0, \iny)$ continuous and increasing with $\Theta(0) = 0$.
	There exists a continuous increasing function $\FMOC \colon [0, \iny) \to \iny$,
	depending on $\Theta$, with $\FMOC(0) = 0$ for which
	\begin{align*}
		\norm{f_1 - f_2}_{L^\iny(\R^d)}
			\le \FMOC(\norm{f_1 - f_2}_{L^2(\R^d)}).
	\end{align*}
\end{lemma}
\begin{proof}
	Fix $x \in \R^d$ arbitrarily and suppose that $\delta = \abs{f_1(x) - f_2(x)} > 0$.
	Let $y$ be in the ball $B$ of radius $a = \Theta^{-1}(\delta/4)$ about $x$,
	so that $\abs{f_1(x) - f_1(y)}, \abs{f_2(x) - f_2(y)} \le \delta/4$.
	Then
	\begin{align*}
		\abs{f_1(y) - f_2(y)}
			&\ge \delta - \abs{f_1(x) - f_1(y)} - \abs{f_2(x) - f_2(y)} 
			= \frac{\delta}{2}.
	\end{align*}
	Hence,
	\begin{align*}
		\norm{f_1 - f_2}_{L^2(\R^d)}
			&\ge \norm{f_1 - f_2}_{L^2(B)}
			\ge \pr{\int_B \pr{\frac{\delta}{2}}^2}^{\frac{1}{2}}
			= \frac{\delta}{2} \sqrt{\pi} a,
	\end{align*}
	or,
	\begin{align*}
		h(\delta)
			:= \frac{\sqrt{\pi}}{2} \delta \Theta^{-1} \pr{\delta/4}
			\le \norm{f_1 - f_2}_{L^2(\R^d)}.
	\end{align*}
	Since $\Theta^{-1}$ must be increasing, so must $h$, so setting $\FMOC = h^{-1}$
	(noting that $\FMOC(0) = 0$) we have
	\begin{align*}
		\abs{f_1(x) - f_2(x)}
			= \delta
			\le \FMOC(\norm{f_1 - f_2}_{L^2(\R^d)}).
	\end{align*}
	This inequality applies for all $x$ even when
	$\delta = \abs{f_1(x) - f_2(x)} = 0$, giving the result. 
\end{proof}

%
%
\section{Boundary differential operators}\label{A:BoundaryCalculations}
\noindent We can define differential operators up to order two on $\prt \Omega$ by treating it as a manifold having at least $C^2$ regularity. In this appendix, we describe the properties that we need of the first-order differential operators, $\grad_\Gamma$, $\dv_\Gamma$, and $\curl_\Gamma$. We refer the reader to standard references for such operators (for instance, Section 2.2 of \cite{Taylor2011VolI}).

We will also have the need to calculate $\grad$, $\dv$, and $\curl$ in 3-space, but restricted to the boundary. This can be done by introducing a convenient coordinate system in a tubular neighborhood of the boundary in such a way that on the boundary itself, the coordinates reduce to a convenient coordinate system on the boundary. This is as done, for instance, in \cite{GK11}, drawing upon \cite{Klingenberg1978}, and we refer the reader to those references for details.

We can define  $\grad_\Gamma$---and then from it, $\dv_\Gamma$ and $\curl_\Gamma$---in a coordinate-free manner by requiring that for any $f \in C^\iny(\Gamma)$ and any smooth curve $\x(s)$ on $\Gamma$ parameterized by arc length,
\begin{align*}
	\grad_\Gamma f \cdot \x'(0)
		&= \lim_{s \to 0} \frac{f(\x(s)) - f(\x(0))}{s}.
\end{align*}
We then define $\dv_\Gamma$ as the adjoint of $\grad_\Gamma$, in the sense of \cref{L:IBPOnBoundary}:

\begin{lemma}\label{L:IBPOnBoundary}
	Let $f \in C^1(\Gamma)$, $\vv \in (C^1(\Gamma))^d$.
	Then
	\begin{align*}
		\int_\Gamma \vv \cdot \grad_\Gamma f
			&= - \int_\Gamma \dv_\Gamma \vv \, f.
	\end{align*}
	Moreover,
	\begin{align}\label{e:dvfv}
		\dv_\Gamma (f \vv)
			&= f \dv_\Gamma \vv + \vv \cdot \grad_\Gamma f.
	\end{align}
\end{lemma}
\begin{proof}
	This is classical for smooth functions (see, for instance, Proposition 2.2.2 of
	\cite{Taylor2011VolI}), and follows in the same way for $C^1$ functions,
	integrating by parts on the boundary in charts.
\end{proof}

Finally, we define (with the $\perp$ operator as in \cref{D:perp})
\begin{align*}
	\curl_\Gamma \vv
		&:= -\dv_\Gamma \vv^\perp.
\end{align*}

We collect now a few useful facts.

For $\uu$, $\vv$ tangent vectors,
\begin{align*}
	(\uu \cdot \grad_\Gamma \vv) \cdot \vv
		= \frac{1}{a_j} u^j \prt_j v^i \, v^i
		= \frac{1}{2 a_j} u^i \prt_j \abs{\vv}^2
		= \frac{1}{2} \uu \cdot \grad \abs{\vv}^2,
\end{align*}
so for any component $\Gamma_n$ of the boundary,
\begin{align*}
	\int_{\Gamma_n} (\uu \cdot \grad_\Gamma \vv) \cdot \vv
		&= \frac{1}{2} \int_{\Gamma_n} \uu \cdot \grad_\Gamma \abs{\vv}^2.
\end{align*}
For a vector field $\vv$ on $\ol{\Omega}$,
\begin{align}\label{e:curlRestricted}
	\curl_\Gamma \vv^\BoldTau
		= (\curl \vv) \cdot \n
\end{align}
and
\begin{align}\label{e:dvCompared}
	\dv \vv
		&= \dv_\Gamma \vv^\BoldTau
			+ \prt_\n v^\n
			+ (\kappa_1 + \kappa_2) v^\n
			\text{ on } \Gamma,
\end{align}
where $\kappa_1, \kappa_2$ are the principal curvatures on $\Gamma$.

\begin{lemma}\label{L:CrossProduct}
	Let $\uu, \vv$ be vector fields on $\ol{\Omega} \subseteq \R^3$. Then
	\begin{align*}
		[\uu \times \vv]^\BoldTau
			&= u^\n [\vv^\BoldTau]^\perp - v^\n [\uu^\BoldTau]^\perp, \quad
		u^\n \vv^\BoldTau - v^\n \uu^\BoldTau
			= [(\vv \times \uu)^\BoldTau]^\perp.
	\end{align*}
\end{lemma}
\begin{proof}
	We have,
	\begin{align*}
		\uu \times \vv
			&= (\uu^\n + \uu^\BoldTau) \times (\vv^\n + \vv^\BoldTau)
			= \uu^\n \times \vv^\BoldTau - \vv^\n \times \uu^\BoldTau
				+ \uu^\BoldTau \times \vv^\BoldTau,
	\end{align*}
	since $\uu^\n \times \vv^\n = 0$. Now, $\uu^\BoldTau \times \vv^\BoldTau$
	is parallel to $\n$, so we see that
	\begin{align*}
		[\uu \times \vv]^\BoldTau
			&= \uu^\n \times \vv^\BoldTau - \vv^\n \times \uu^\BoldTau.
	\end{align*}
	But, $\uu^\n$ is perpendicular to $\vv^\BoldTau$, so we see that
	$\uu^\n \times \vv^\BoldTau = u^\n [\vv^\BoldTau]^\perp$,
	and similarly, $\vv^\n \times \uu^\BoldTau = v^\n [\uu^\BoldTau]^\perp$.
	Hence,
	$
		[\uu \times \vv]^\BoldTau
			= u^\n [\vv^\BoldTau]^\perp - v^\n [\uu^\BoldTau]^\perp
	$,
	giving also
	$u^\n \vv^\BoldTau - v^\n \uu^\BoldTau = [(\vv \times \uu)^\BoldTau]^\perp$.
\end{proof}

\begin{proof}[\textbf{Proof of \cref{P:ugradunOnGammaPlus}}]
	All the following calculations are on $\Gamma$.
	We start with a short calculation 
	\textit{in rectangular coordinates},
	using that $\dv \uu = \prt_i u^i = 0$:
	\begin{align*}
		(\uu \cdot \grad \uu) \cdot \n
			&= u^i \prt_i u^j n^j
			= \prt_i (u^i u^j n^j)
				- u^j u^i \prt_i n^j
			= \dv (u^\n \uu)
				- \uu \cdot (\uu \cdot \grad \n) \\
			&= \dv (u^\n \uu)
				- \uu^\BoldTau \cdot \Cal{A} \uu^\BoldTau.
	\end{align*}
	In the last equality, we used that because $\n$ does not change in the direction of $\n$,
	\begin{align*}
		\uu \cdot \grad \n
			&= (\uu^\n \cdot \grad) \n + \uu^\BoldTau \cdot \grad \n
			= \Cal{A} \uu^\BoldTau,
	\end{align*}
	which is a tangent vector.

	From \cref{e:dvCompared} followed by \cref{e:dvfv}, then,
	\begin{align*}
		\dv (u^\n \uu)
		&= \dv_\Gamma (u^n \uu^\BoldTau)
			+ \prt_\n  (u^\n)^2
			+ (\kappa_1 + \kappa_2) (u^\n)^2 \\
		&= u^\n \dv_\Gamma \uu^\BoldTau + \uu^\BoldTau \cdot \grad_\Gamma u^\n
			+ \prt_\n  (u^\n)^2
			+ (\kappa_1 + \kappa_2) (u^\n)^2.
	\end{align*}
	
	Using \cref{e:dvCompared} again,
	\begin{align*}
		0
			&= (\dv \uu) u^\n
			= (\dv_\Gamma \uu^\BoldTau + \prt_\n u^\n + (\kappa_1 + \kappa_2) u^\n) u^\n,
	\end{align*}
	so
	\begin{align*}
		\prt_\n (u^\n)^2
			&= 2 u^\n \prt_\n u^\n
			= - 2 u^\n \dv_\Gamma \uu^\BoldTau - 2 (\kappa_1 + \kappa_2) (u^\n)^2.
	\end{align*}
	Hence,
	\begin{align*}
		(\uu \cdot \grad \uu) \cdot \n
		&= -u^\n \dv_\Gamma \uu^\BoldTau + \uu^\BoldTau \cdot \grad_\Gamma u^\n
			- (\kappa_1 + \kappa_2) (u^\n)^2
			- \uu^\BoldTau \cdot \Cal{A} \uu^\BoldTau.
		\qedhere
	\end{align*}
\end{proof}

%
%
\section{Compatibility conditions: special case}\label{A:CompCondsSpecialCase}

\noindent 
In \cite{TemamWang2002}, Temam and Wang consider a periodic domain with $\UU = (0, 0, -1)$, so $\UU^\BoldTau = 0$ for all time. More generally, the authors of \cite{GieHamoudaJungTemam2018} consider $\UU = -U^I \n$, where $U^I > 0$ is constant, so $\UU^\BoldTau = 0$ on $\Gamma_+$ for all time. The compatibility conditions simplify in these settings.

\begin{prop}\label{P:TWCompConds}
	Assume that $\UU^\BoldTau \equiv 0$ and $U^\n$ is spatially constant along $\Gamma_+$
	($U^\n$ need not be constant in time).
	Then the compatibility condition $\cond_N$ for $N \ge 0$ is
	\begin{align}\label{e:TWLikeCompCond}
		\prt_t^j \f^\BoldTau|_{t = 0}
			= \prt_t^j \grad_\Gamma p|_{t = 0}
				- U_0^\n (\prt_t^j \bomega^\BoldTau)^\perp|_{t = 0}
			\text{ for all } 0 \le j \le N,
	\end{align}
	where $\prt_t^j \grad_\Gamma p|_{t = 0}$ and
	$\prt_t^j \bomega|_{t = 0}$ must be treated as explained following
	\cref{e:condN}.
\end{prop}
\begin{proof}
	Since $\uu^\BoldTau = \UU^\BoldTau = 0$, \cref{e:curlRestricted} gives that on $\Gamma_+$,
	\begin{align*}
		\omega^\n
			= \bomega \cdot \n
			= \curl_\Gamma \uu^\BoldTau
			= 0.
	\end{align*}
	In particular, this holds at time zero.
	Both $\prt_t \UU^\BoldTau = 0$ and $\curl_\Gamma \UU^\BoldTau = 0$, while
	$\abs{\UU}^2 = (U^\n)^2$ is constant on $\Gamma_+$, so also
	$\grad_\Gamma \abs{\UU}^2 = 0$. We see, then, that $\H^\BoldTau$ simplifies to
	$
		\H^\BoldTau
			= (U^\n)^{-1}
				\brac{\f^\BoldTau - \grad_\Gamma p}^\perp
	$,
	so $\lincond_0$ (which follows from $\cond_0$ by \cref{P:condImplieslincond0}) becomes
	\begin{align*}
		\brac{\f^\BoldTau  - \grad_\Gamma p}^\perp_{t = 0}
			= U_0^\n \bomega_0^\BoldTau,
	\end{align*}
	which is \cref{e:TWLikeCompCond} for $N = 0$. The inductive extension of this
	to higher $N$ follows readily, leading to \cref{e:TWLikeCompCond} for $N \ge 0$.
\end{proof}

The condition in \cref{e:TWLikeCompCond} for $N = 0$ also follows from $\cond_0$ with slightly more work, though the inductive extension to higher $N$ is not so transparent as it is starting from $\cond_0'$.

Because $\dv \f = 0$ with $\f \cdot \n = 0$ on $\Gamma$, $\f$ plays no role in the calculation of $\grad_\Gamma p$ for $N = 0$. By writing the condition in \cref{e:TWLikeCompCond} as we do, we are stressing that, given initial data one can always choose a forcing at time zero so that $\cond_0$ is satisfied.  

For all $N \ge 1$, though, forcing enters into the calculation of $\prt_t \grad_\Gamma p$, when $\prt_t \uu_0$ is replaced by $\f(0) - \uu_0 \cdot \grad \uu_0 - \grad p^0$: even though $\f \cdot \n = 0$, the forcing still does not, in general, vanish from even the $N = 1$ condition. Because of this fact, the forcing is intimately entwined in $\cond_N$ for $N \ge 1$, appearing on both sides of the condition, even for the simplest nontrivial case considered in \cite{TemamWang2002}. These same comments hold in the general setting, but are more transparent in this simplified setting.

\def\cprime{$'$} \def\polhk#1{\setbox0=\hbox{#1}{\ooalign{\hidewidth
  \lower1.5ex\hbox{`}\hidewidth\crcr\unhbox0}}}

\end{document}